\documentclass[a4paper,12pt]{amsart}

\usepackage[margin=2cm]{geometry}

\usepackage{graphicx}
\usepackage{tikz}
\usepackage{tikz-cd}
\usepackage{amssymb, amsmath, mathtools}
\usetikzlibrary{decorations.markings}

\usepackage{hyperref}
\usepackage{floatrow}
\newtheorem{theorem}{Theorem}[section]
\newtheorem{proposition}[theorem]{Proposition}
\newtheorem{lemma}[theorem]{Lemma}
\newtheorem{corollary}[theorem]{Corollary}

\theoremstyle{definition}

\theoremstyle{remark}
\newtheorem{remark}[theorem]{Remark}

\numberwithin{equation}{section}
\numberwithin{figure}{section}

\newcommand{\Kk}{\mathcal{K}}
\newcommand{\Z}{\mathbb{Z}}

\newcommand{\N}{\mathbb{N}}
\newcommand{\T}{\mathbb{T}}
\newcommand{\lb}{\lambda}
\newcommand{\Lb}{\Lambda}

\newcommand{\cantor}{\mathbf{K}}

\newcommand{\image}{\operatorname{im}}
\newcommand{\lt}{\operatorname{lt}}
\newcommand{\spann}{\operatorname{span}}
\newcommand{\Aut}{\operatorname{Aut}}
\newcommand{\coker}{\operatorname{coker}}
\newcommand{\MCE}{\operatorname{MCE}}
\newcommand{\clsp}{\operatorname{\overline{\mathrm{span}}}}

\thanks{This research was supported by grant DP180100595 of the Australian Research Council  and  Marsden Grant 18-VUW-056 of the Royal Society of New Zealand. We thank Jamie Gabe for helpful discussions.}

\begin{document}
	
\title{Stably finite extensions of rank-two graph $C^*$-algebras}

\author[Astrid an Huef]{Astrid an Huef}
\address[A.\ an Huef]{School of Mathematics and Statistics, Victoria University of Wellington, P.O. Box 600, Wellington 6140, New Zealand.}
\email{astrid.anhuef@vuw.ac.nz}
\author[Abraham Ng]{Abraham C.S.\ Ng}
\author[Aidan Sims]{Aidan Sims}
\address[A.C.S.\ Ng and A.\ Sims]{School of Mathematics and Applied Statistics,
University of Wollongong,
NSW 2522, Australia.}

\email{abrahamn and asims@uow.edu.au}

	
\date{\today}
	
	
\begin{abstract}
We study stable finiteness of extensions of $2$-graph $C^*$-algebras determined by saturated hereditary sets of vertices. We use two iterations of the Pimsner--Voiculescu sequence to calculate the map in $K$-theory induced by the inclusion of a hereditary subgraph into the larger $2$-graph it lives in. We then apply a theorem of Spielberg about stable finiteness of extensions to provide a sufficient condition for the $C^*$-algebra of the larger $2$-graph to be stably finite. We illustrate our results with examples.
\end{abstract}
	
\maketitle
	
\section{Introduction}
In this paper we investigate when the $C^*$-algebra of a row-finite $2$-graph with no sources is stably finite by applying a theorem of Spielberg \cite{Spielberg1988} about stable finiteness of extensions of $C^*$-algebras.

A  $C^*$-algebra $A$ is stably finite if no matrix algebra over $A$ contains a partial isometry whose initial projection strictly dominates its final projection. This yields an easy sufficient condition for stable finiteness:  if $A$ admits a faithful trace, then it is stably finite. This has been exploited in a number of analyses of $C^*$-algebras arising from dynamical systems that are minimal in a suitable sense \cite{BL, CaHS, RS}---the idea is to identify an order-theoretic property in $K$-theory that guarantees the existence of infinite projections, then show that in the absence of this condition the $C^*$-algebra admits a nonzero tracial state, and then use minimality to show that this trace must be faithful.

A motivating example for us is the class of $C^*$-algebras associated to higher-rank graphs or $k$-graphs. A $k$-graph $\Lambda$ is like a $k$-dimensional directed graph; it is made up of  $k$ ordinary directed graphs (distinguished by $k$ different colours) sharing a common vertex set $\Lambda^0$, together with factorisation rules that organise pairs of bi-coloured paths into commuting squares. The associated $C^*$-algebra $C^*(\Lambda)$ is generated by copies of the graph $C^*$-algebras of each of the $k$ singly-coloured subgraphs subject to commutation rules between generators that correspond to the factorisation rules. The Cuntz--Krieger relations for the $k$ subgraphs determine equivalences between projections in matrix algebras over the $k$-graph $C^*$-algebra $C^*(\Lambda)$. This leads to a natural necessary condition on the adjacency matrices of the subgraphs, for stable finiteness of $C^*(\Lambda)$ \cite[Proposition~3.1]{CaHS}. The main result of \cite{CaHS} is that if $\Lambda$ is cofinal (which is equivalent to minimality of an underlying groupoid \cite{KP2000}) and satisfies the matrix condition mentioned above, then there exists a faithful  graph trace \cite{PRS2005}  on $\Lambda$, and therefore a faithful trace on $C^*(\Lambda)$. It follows that the matrix condition is necessary and sufficient for stable finiteness of $C^*(\Lambda)$ when $\Lambda$ is cofinal. A similar program has since been implemented for $C^*$-algebras of minimal ample groupoids \cite{BL, RS}.

This program typically fails for dynamical systems that are not minimal, because $C^*$-algebras associated to such systems need not admit a faithful trace even when they are stably finite (we describe one example of this phenomenon in Section~\ref{subsec:spine example}). However, for the $C^*$-algebras of $2$-graphs, another approach is available. In this case, if  one of the two singly-coloured subgraphs of $\Lambda$ has no cycles, then $C^*(\Lambda)$ is Morita equivalent to a crossed product of an AF algebra $A$ by an automorphism $\alpha$ of $A$ that implements an action of the integers $\Z$. Brown addressed stable finiteness of such crossed products in \cite{Brown}. He observed that there is a very natural obstruction to stable finiteness involving the interaction of the induced automorphism $\alpha_*$ of $K_0(A)$ with the positive cone $K_0(A)_+$ of the $K_0$-group. Specifically, if the crossed product is stably finite, then the range of $1 - \alpha_*$ includes no nonzero elements of $K_0(A)_+$. He then employed a deep analysis, using powerful results on classification of AF algebras and their automorphisms to prove the converse:  if the range of $1 - \alpha_*$ contains no nonzero positive elements, then the crossed product embeds in an AF algebra, and hence is stably finite. It turns out that when $A$ and $\alpha$ arise from a $2$-graph with no cycles in one of its singly-coloured subgraphs, the range of $1 - \alpha_*$ contains no nonzero elements of $K_0(A)_+$ if and only if the $2$-graph satisfies the matrix condition of \cite{CaHS}. So for any cofinal $2$-graph or  any $2$-graph that contains no cycles of at least one colour, the matrix condition characterises stable finiteness.

 But it is easy to construct an example of a $2$-graph which is not cofinal and contains cycles of both colours. To tackle such $2$-graphs, it is natural to try to address stable finiteness in terms of extensions of $C^*$-algebras. Stable finiteness has a number of nice permanence properties (it is preserved under passing to subalgebras, under Morita equivalence and under direct limits), but in general it is poorly behaved with respect to extensions.
 However, a result of Spielberg \cite{Spielberg1988} characterises the obstruction to stable finiteness of an extension of one stably finite $C^*$-algebra by another. If $0 \to I \to A \to B \to 0$ is an exact sequence in which $I$ and $B$ are both stably finite, then $A$ is stably finite exactly when the homomorphism of $K_0$-groups induced by the inclusion $I \hookrightarrow A$ satisfies the positivity condition that its kernel contains no nonzero elements of $K_0(I)_+$.

In this paper we investigate how this applies to extensions of $2$-graph $C^*$-algebras. The ideal structure of such $C^*$-algebras is fairly well understood. In particular, the gauge-invariant ideals are indexed by so-called saturated hereditary sets $H$ of vertices. Moreover, the ideal $I_H$ corresponding to a given $H$ is Morita equivalent to the $C^*$-algebra of a sub-2-graph $H\Lambda$ while the quotient is isomorphic to the $C^*$-algebra of a second sub-2-graph $\Lambda \setminus \Lambda H$. But to apply Spielberg's positivity condition, we need to understand the kernel of the map in $K$-theory induced by the inclusion $i_H: C^*(H\Lambda) \hookrightarrow C^*(\Lambda)$.

This is difficult for two reasons. The first difficulty is that the existing computations of $K$-theory for $2$-graph $C^*$-algebras \cite{RobSte, APS, Evans2008} employ Kasparov's spectral sequence to give a description of the $K_0$-group as an extension (indeed, a direct sum) of the quotient of $\Z\Lambda^0$ determined  by the Cuntz--Krieger relations, and a second more mysterious summand. We have not found a way to extract from these results any statement about the map $(\iota_H)_*$ between these extensions induced by the inclusion. However, intuition and informal calculations suggest that the first summand in the $K_0$-group is just the subgroup generated by the classes of the vertex projections, that $(\iota_H)_*$ is a homomorphism of extensions, that the kernel of $(\iota_H)_*$ is just the kernel of its restriction to the first summand, and that on this summand $\iota_H$ is induced by the canonical inclusion of $\Z H$ into $\Z\Lambda^0$. To prove this, we re-calculate the $K$-theory of $C^*(\Lambda)$ by iterating the Pimsner--Voiculescu sequence---this allows us to use naturality of the sequence and the concreteness of the connecting maps appearing in it. We pay a price for using our approach instead of Evans' in \cite{Evans2008}: his calculation gives a direct formula for $K_1(C^*(\Lambda))$, whereas our approach only fits it into a short exact sequence of abelian groups.

The second difficulty is that existing calculations of the $K$-theory of a $2$-graph algebra provide no information about the positive cone in $K_0$, which we need in order to apply Spielberg's theorem. Our new calculation using the Pimsner--Voiculescu sequence does not completely address this issue. However, we are able to show that Spielberg's positivity condition follows from the matrix condition for $\Lambda$ provided that the subgraph $H\Lambda$ satisfies a condition~\eqref{nicecond} which asks that the intersection of the positive cone of $K_0(C^*(H\Lambda))$ with the subgroup generated by the classes $[p_v]$ of vertex projections is precisely the collection of nonnegative-integer combinations of classes of vertex projections. We also relate condition~\eqref{nicecond} to the following natural, and open, question for crossed products by $\Z$: if $i : A \to A \rtimes \Z$ is the inclusion of a $C^*$-algebra in a crossed product by $\Z$, when is $i_*(K_0(A)_+)$ equal to $i_*(K_0(A)) \cap K_0(A \rtimes \Z)_+$? We provide an example in Section~\ref{subsec:handlebar example}  that shows our condition~\eqref{nicecond} is checkable at least for some concrete $2$-graphs that are not accessible to the existing characterisations of stable finiteness in \cite{CaHS}.

\smallskip

The paper is organised as follows. We start with the necessary preliminaries and notation in Section~\ref{sec:prelims}. In Section~\ref{sec:graph algebras} we revisit the calculation of $K$-theory for graph $C^*$-algebras due to Raeburn and Szyma\'nski, because we need an explicit description of the homomorphisms  in $K$-theory that are induced by naturality of the Pimsner--Voiculescu exact sequence applied to the inclusion of the $C^*$-algebra of a hereditary subgraph into a graph $C^*$-algebra.

We then go on in Section~\ref{sec: k theory main} to the main bulk of the work: the computation of $K$-theory, and induced maps between $K$-groups, for $2$-graph algebras. This is long and complex, so we break it up into a number of steps. We consider a $2$-graph $\Lambda$ and a hereditary subset $H$ of $\Lambda^0$. Starting from Section~\ref{subsec:1-dim skew}, we consider the skew products $\Lambda \times_{d_1} \Z$ and $H\Lambda \times_{d_1} \Z$ of $\Lambda$ and of $H\Lambda$ by the first coordinate of the degree maps. By realising the $C^*$-algebra of the skew-product as a direct limit of subalgebras that are Morita equivalent to the second coordinate graph of $\Lambda$ and then applying our results from Section~\ref{sec:graph algebras}, we compute the $K$-theory of $C^*(\Lambda \times_{d_1} \Z)$ and the map in $K$-theory induced by the inclusion $H\Lambda \times_{d_1} \Z \hookrightarrow \Lambda \times_{d_1} \Z$. This computation occupies Sections \ref{subsec:1-dim skew}--\ref{subsec:step3}. With these key technical results in hand, we can apply the Pimsner--Voiculescu sequence again, together with Takai duality, to obtain our first main theorem, Theorem~\ref{thm:inc}: a description of the $K$-theory of $C^*(\Lambda)$ and of the map in $K$-theory induced by the inclusion $H\Lambda \hookrightarrow \Lambda$. We give a second, slightly cleaner, description of the $K_0$-groups and the maps between them in Theorem~\ref{thm:inc nice}. Specifically let $A_1$ and $A_2$ be the $\Lambda^0 \times \Lambda^0$ adjacency matrices whose $v,w$ entries give the number of edges of the two different colours from $w$ to $v$. Then the block matrices
\[
(1 - A_1^t, 1 - A_2^t) \quad\text{ and }\quad \begin{pmatrix}1 - A_1^t\\1 - A_2^1\end{pmatrix}
\]
determine homomorphisms between $\Z\Lambda^0 \oplus \Z\Lambda^0$ and $\Z\Lambda^0$, and $K_0(C^*(\Lambda))$ is an extension of the kernel of the second of these homomorphisms by the cokernel of the first. Moreover, the map $(\iota_H)_*$ in $K$-theory respects this decomposition and intertwines with the maps $\tilde\iota_H$ between cokernels and $\iota_H|$ between kernels induced by the inclusion $\Z H \hookrightarrow \Z\Lambda^0$.

In Section~\ref{sec:stablyfiniteextensions} we use our $K$-theory calculations from the preceding section to describe how Spielberg's positivity condition applies to an extension of the form
\[
0 \to I_H \to C^*(\Lambda) \to C^*(\Lambda \setminus \Lambda H) \to 0
\]
arising from a saturated hereditary subset $H$ of $\Lambda$. We first show in Lemma~\ref{lem:scl1} that if $\Lambda$ satisfies the matrix condition of \cite{CaHS}, then so does the subgraph $H\Lambda$, and the kernel of the map $\tilde\iota_H$ between cokernels described in the previous paragraph has trivial intersection with the image of $\N H$ in $(1 - A_{H,1}^t, 1 - A_{H,2}^t)$.
We then introduce a positive-elements condition on a $2$-graph $\Lambda$, that the image of $\N \Lambda^0 + \image(1 - A_{H,1}^t, 1 - A_{H,2}^t)$ in $K_0(C^*(\Lambda))$ agrees with the intersection of $K_0(C^*(\Lambda))_+$ with the subgroup generated by the vertex projectons. We show in Proposition~\ref{prop:equivalent conditions} that this condition is equvalent to the condition~\eqref{nicecond} mentioned earlier. We use these two results to obtain our second main theorem, Theorem~\ref{thm:simplified stably finite restated}, which says that if $\Lambda$ satisfies the matrix condition, the subgraph $H\Lambda$ satisfies condition~\eqref{nicecond} and both $C^*(H\Lambda)$ and $C^*(\Lambda \setminus \Lambda H)$ are stably finite, then so is $C^*(\Lambda)$. We end the section with Corollary~\ref{cor: maximal chain} which gives a sufficient condition for stable finiteness of the $C^*$-algebra of a $2$-graph that admits a finite maximal chain of saturated hereditary sets.

We close the paper in Section~\ref{sec:examples} with two examples. The first is an example of a $2$-graph that does not satisfy the hypotheses of the theorems of \cite{CaHS}, but whose $C^*$-algebra we can prove is stably finite using Theorem~\ref{thm:simplified stably finite restated} for an appropriate choice of saturated hereditary set.

The second example is a $2$-graph $\Lambda$ that first appeared in \cite{EvansPhD} and was considered again in \cite{ES}. Its $C^*$-algebra  is already known to be stably finite by the theorems of \cite{CaHS}, but it provides a good illustration of the conditions on our main results. It also provides an example in which $C^*(\Lambda)$ satisfies the matrix condition even though it has a subgraph of the form $\Lambda \setminus \Lambda H$ that does not. Further $C^*(\Lambda)$ is stably finite (indeed, AF embeddable) even though $\Lambda$ admits no faithful graph trace in the sense of \cite{PRS2005}, and hence $C^*(\Lambda)$ admits no faithful trace.

	\section{Preliminaries}\label{sec:prelims}
	
	Throughout this paper, we write $\N$ for the additive abelian monoid $\{0,1,2, \dots\}$. For $k \ge 1$, $\N^k$ then denotes the monoid of $k$-tuples of elements of $\N$ under pointwise addition. The usual generators of $\N^k$ are denoted $e_1, \dots, e_k$, so for $n \in \N^k$, we have $n = \sum^k_{i=1} n_i e_i$. For $m,n \in \N^k$ we write $m \vee n$ for the coordinatewise maximum of $m,n$, which is the least common upper bound of $m,n$ in the usual algebraic order. We often regard monoids as categories with a single object.
	
	The next two subsections summarise some key definitions and notation for $k$-graphs, mostly taken from the original paper \cite{KP2000} on the subject. For more details, see \cite{AidanAbeNotes}.
	
	\subsection{\texorpdfstring{$k$}{k}-graphs}\label{ss:kgraphs}
	A $k$-graph is a countable category $\Lambda$ endowed with a functor $d : \Lambda \to \N^k$ that has the \emph{factorisation property}: for all $\lambda \in \Lambda$ and all $m,n \in \N^k$ such that $m+n = d(\lambda)$ there are unique elements $\mu \in d^{-1}(m)$ and $\nu \in d^{-1}(n)$ such that $\lambda = \mu\nu$. We write $\Lambda^n$ for $d^{-1}(n)$. We think of $\Lambda$ as the path category of a $k$-dimensional directed graph, so we often refer to its morphisms as \emph{paths}. The factorisation property guarantees that $\Lambda^0$ is precisely the collection of identity morphisms. We identify the objects of $\Lambda$ with $\Lambda^0$, whose elements we call \emph{vertices}, and we write $r,s : \Lambda \to \Lambda^0$ for codomain and domain maps respectively. For $S, T \subseteq \Lambda$, we write $ST \coloneqq  \{\lambda\mu : \lambda \in S,\mu \in T\}$. If $S = \{\lambda\}$, we write $\lambda T$ rather than $\{\lambda\}T$ and similarly if $T$ is a singleton. In particular, for $v \in \Lb^0$, we have $vS = S \cap r^{-1}(v)$ and $Sv = S \cap s^{-1}(v)$.
	
	The $k$-graph $\Lambda$ has \emph{no sources} if $v\Lambda^n \not= \emptyset$ for all $v \in \Lambda^0$ and $n \in \N^k$, and it is \emph{row-finite} if $v\Lambda^n$ is finite for all $v \in \Lambda^0$ and $n \in \N^k$. In this paper we work exclusively with $k$-graphs that are row-finite and have no sources. Given a $k$-graph $\Lambda$, and given $\mu,\nu \in \Lambda$, we write $\MCE(\mu,\nu) \coloneqq  \mu\Lambda \cap \nu\Lambda \cap \Lambda^{d(\mu) \vee d(\nu)}$ for the collection of \textit{minimal common extensions} of $\lambda$ and $\mu$.
	
	Given an abelian group $(G, +)$ and a functor $c : \Lambda \to G$, we can form the \emph{skew product} $\Lambda \times_c G$ which is equal as a set to $\Lambda \times G$, has objects $\Lambda^0 \times G$, has range and source maps $r(\lambda, g) = (r(\lambda), g)$ and $s(\lambda, g) = (s(\lambda), g + c(\lambda))$, and has composition $(\lambda, g)(\mu, g + c(\lambda)) = (\lambda\mu, g)$. This skew-product is itself a $k$-graph with degree map $d(\lambda,g) = d(\lambda)$.
	
	The \emph{adjacency matrices} of the $k$-graph $\Lambda$ are the $\Lambda^0 \times \Lambda^0$ nonnegative integer matrices $A_1, \dots, A_k$ with entries given by $A_i(v,w) = |v\Lambda^{e_i}w|$; as usual, $A_i^t$ then denotes the transposed matrix $A_i^t(v,w) = |w\Lambda^{e_i} v|$. The factorisation property implies that
\[
A_iA_j(v,w)
    = |v\Lb^{e_i}\Lb^{e_j}w|
    = |v\Lb^{e_i+e_j}w|
    = |v\Lb^{e_j}\Lb^{e_i}w|
    = A_jA_i(v,w)
\]
for all $v,w \in \Lb^0$ so that the matrices $A_1,\dots, A_k$ pairwise commute. We regard these matrices as maps from the free abelian group $\Z \Lb^0$ to itself.
	
%
	
%

\subsection{\texorpdfstring{$k$}{k}-graph \texorpdfstring{$C^*$}{C*}-algebras} \label{ss:cuntzkrieger}
	
Let $\Lambda$ be a row-finite $k$-graph with no sources. A \emph{Cuntz--Krieger $\Lambda$-family} in a $C^*$-algebra $A$ is a map $\lambda \mapsto s_\lambda$ from $\Lambda \to A$ such that
\begin{enumerate}
	\item $\{s_v : v \in \Lambda^0\}$ is a family of mutually orthogonal projections;
	\item $s_\mu s_\nu = s_{\mu\nu}$ whenever $s(\mu) = r(\nu)$;
	\item $s_\mu^*s_\mu = s_{s(\mu)}$ for all $\mu \in \Lambda$; and
	\item $\sum_{\lambda \in v\Lambda^n} s_\lambda s^*_\lambda = s_v$ for all $v \in \Lambda^0$ and $n \in \N^k$.
\end{enumerate}
We often write $p_v$ rather than $s_v$ for $v \in \Lambda^0$ to emphasise that these are
projections. However when looking at Cuntz--Krieger families for graphs skewed by $\Z^2$, we use
$s_{\alpha}$ even when $\alpha$ is a vertex, as $p$ in that setting is often an integer. It
follows from these relations that for each $n \in \N^k$, the collection $\{s_\lambda : \lambda \in
\Lambda^n\}$ is a family of partial isometries with mutually orthogonal ranges, and that for any
$\mu,\nu \in \Lambda$ and $n \ge d(\mu) \vee d(\nu)$, we have $s_\mu s^*_\mu s_\nu s^*_\nu =
\sum_{\lambda \in \mu\Lambda \cap \nu\Lambda \cap \Lambda^n} s_\lambda s^*_\lambda$, and $s_\mu^*
s_\nu = \sum_{\mu\alpha=\nu\beta \in \Lambda^n} s_\alpha s_\beta^*$.
	
The $k$-graph $C^*$-algebra, denoted $C^*(\Lambda)$ is the universal $C^*$-algebra generated by a Cuntz--Krieger $\Lambda$-family. It follows from the preceding paragraph that $C^*(\Lambda) = \clsp\{s_\mu s^*_\nu : s(\mu) = s(\nu)\}$. The universal property ensures that there is a strongly continuous action $\gamma : \T^k \to \Aut(C^*(\Lambda))$, called the \emph{gauge action} such that $\gamma_z(s_\mu) = z^{d(\mu)} s_\mu$ for all $\mu \in \Lambda$.

\subsection{Saturated hereditary sets and ideals}\label{subsec:H<->I}
	
Let $\Lambda$ be a row-finite $k$-graph with no sources. A subset $H \subseteq \Lambda^0$ is \emph{hereditary} if $s(H\Lambda) \subseteq H$. It is \emph{saturated} if whenever $s(v\Lambda^n) \subseteq H$ we also have $v \in H$. If $H$ is a hereditary subset of $\Lambda^0$ then the smallest saturated hereditary set $\overline{H}$ containing $H$ is called the \emph{saturation} of $H$. The following description of gauge-invariant ideals of $C^*(\Lambda)$ is proved in \cite{RaeburnSimsYeend2003}. Given a hereditary set $H \subseteq \Lambda^0$, the space $I_H \coloneqq  \clsp\{s_\mu s^*_\nu : \mu,\nu \in \Lambda H\}$ is a gauge-invariant ideal of $C^*(\Lambda)$. The set $\{v \in \Lambda^0 : p_v \in I_H\}$ is equal to $\overline{H}$. Given an ideal $I$ of $C^*(\Lambda)$, the set $H_I \coloneqq  \{v \in \Lambda^0 : p_v \in I\}$ is a saturated hereditary set. The map $H \mapsto I_H$ from the collection of \emph{saturated} hereditary subsets of $\Lambda^0$ to the collection of gauge-invariant ideals of $C^*(\Lambda)$ is a lattice isomorphism with inverse $I \mapsto H_I$.
	
Given any hereditary set $H \subseteq \Lambda^0$, the series $\sum_{v \in H} p_v$ converges strictly to a projection $P_H$ in the multiplier algebra of $I_H$. The corner $P_H I_H P_H$ is full and is equal to $\clsp\{s_\mu s^*_\nu : \mu,\nu  \in H\Lambda\}$. The set $H \Lambda$ is a sub-$k$-graph of $\Lambda$, and the inclusion $H\Lambda \hookrightarrow \Lambda$ induces an inclusion $\iota_H : C^*(H \Lambda) \to C^*(\Lambda)$ sending generators to generators, which is an isomorphism of $C^*(H\Lambda)$ onto $P_H I_H P_H$. In particular, $\iota_H$ induces an isomorphism $K_*(C^*(H\Lambda)) \cong K_*(I_H)$.
	
	If $\Lambda$ is a row-finite $2$-graph with no sources and $H$ is a hereditary subset of $\Lambda^0$, then for any $v \in H$, any $i \le k$ and any $w$ such that $A_i(v,w) \not= 0$, we have $v\Lambda^{e_i}w \not= \emptyset$ and hence $w \in H$. It follows that partitioning $\Lb^0=(\Lb^0\setminus H)\sqcup H$ yields block decompositions
	$$A_i =  \begin{pmatrix}
		A_{\Lb^0\setminus H,i} & * \\
		0 & A_{H,i}
	\end{pmatrix} \quad \text{ and } \quad A_i^t = \begin{pmatrix}
		A_{\Lb^0\setminus H,i}^t & 0 \\
		* & A_{H,i}^t
	\end{pmatrix}.$$
	Hence we have the following result that we will use (often implicitly) throughout the paper.
	
	\begin{lemma}\label{lem:fundamental}
		Let $\Lb$ be a row-finite $2$-graph with no sources and let $H$ be a hereditary subset of $\Lb$. For $i=1,2$, the inclusion $\iota : \Z H \to \Z\Lb^0$ restricts to inclusions $\iota : \ker(1-A_{H,i}^t)\subset \Z H \to \ker(1-A_i^t) \subset \Z \Lb^0$ and $\iota : \image(1-A_{H,i}^t) \subset \Z H \to \image(1-A_i^t) \subset \Z \Lb^0$.
	\end{lemma}

\subsection{Skew-graphs and their \texorpdfstring{$C^*$}{C*}-algebras}

Let $\Lb$ be a row-finite $2$-graph with no sources and degree map $d$. For $i = 1,2$, define $d_i : \Lambda \to \N$ by $d_i(\lambda) = d(\lambda)_i$, the $i$-th coordinate of $d(\lambda)$, and let $\Lambda \times_{d_i}\Z$ be the skew-product $2$-graph defined in Subsection \ref{ss:kgraphs}.
We also consider $\Lb \times_d \Z^2$, the skew-product of $\Lb$ by the degree map $d:\Lb \to \N^2$. Observe that $((\lb,m),n) \mapsto (\lb,(m,n))$ is an isomorphism of $(\Lb \times_{d_1}\Z)\times_{d_2}\Z$ onto $\Lb \times_d \Z^2$.

The automorphism $(\lb,n) \mapsto (\lb,n-1)$ of $\Lb\times_{d_1}\Z$ induces an automorphism $\lt_1$ of $C^*(\Lb\times_{d_1}\Z)$ such that $\lt_1(s_{(\lb,n)}) = s_{(\lb,n-1)}$.

Similarly for the skew-graph $\Lambda\times_d\Z^2$, we let $\lt_1$ denote the automorphism of
$C^*(\Lambda\times_d\Z^2)$ such that $\lt_1(s_{(\lambda,n,m)}) = s_{(\lambda,n-1,m)}$, and let $\lt_2$ denote the automorphism such that $\lt_2(s_{(\lambda,n,m)}) = s_{(\lambda,n,m-1)}$.


The following lemma deals with hereditary saturated subsets and skew-graphs.

	\begin{lemma}\label{lem:saturated product}
		Let $\Lambda$ be a row-finite $2$-graph with no sources and let $H$ be a saturated hereditary subset of $\Lambda^0$. Then the subset  $H \times \Z$  of $(\Lambda \times_{d_1}\Z)^0$ is  hereditary and saturated, and induces  an inclusion $\iota_{H \times \Z} : C^*(H\Lambda \times_{d_1} \Z) \to
		C^*(\Lambda \times_{d_1} \Z)$.
	\end{lemma}

	\begin{proof}
		If $(\mu,n) \in (H\times\Z)(\Lb\times_{d_1}\Z)$, then $r(\mu,n) = (r(\mu),n) \in H\times \Z$ so that $r(\mu) \in H$. Since $H$ is hereditary, $s(\mu) \in H$ so that $s(\mu,n) = (s(\mu),n + d_1(\mu)) \in H\times \Z$. Thus $H\times\Z$ is a hereditary subset of $(\Lb\times_{d_1}\Z)^0$. Suppose now that $s\left((v,m)(\Lb\times_{d_1}\Z)\right) \subset H\times\Z$. Fix $\mu \in v\Lb^n$. Then $(\mu,m) \in (v,m)(\Lb\times_{d_1}\Z)$ so that $(s(\mu),m+d_1(\mu)) = s(\mu,m) \in H\times\Z$. In particular, $s(\mu) \in H$. Hence $s(v\Lb^n) \subset H$. As $H$ is saturated, $v\in H$, and so $(v,m) \in H\times\Z$. Thus $H\times\Z$ is saturated. Since $(H\times\Z)(\Lb\times_{d_1}\Z) = H\Lb \times_{d_1}\Z$, the result follow.
	\end{proof}

	\subsection{Stable finiteness and extensions of \texorpdfstring{$C^*$}{C*}-algebras}
	
	A $C^*$-algebra $A$ is \emph{finite} if it contains no partial isometry $S$ such that $SS^* \lneq S^*S$. It is \emph{stably finite} if its stabilisation $A \otimes \Kk$ is finite.
	
	As mentioned in the introduction, stable finiteness passes to subalgebras, but has poor permanence properties with respect to quotients and extensions. However, Spielberg proved the following theorem that describes when an extension of one stably finite $C^*$-algebra by another is again stably finite. We will use this theorem to prove our main result.
	
	\begin{theorem}[{\cite[Lemma 1.5]{Spielberg1988}}] \label{thm:spielberg}
		Suppose that $0 \to I \stackrel{i}{\to} A \to B \to 0$ is an exact sequence of $C^*$-algebras and that both $I$ and $B$ are stably finite. Then $A$ is stably finite if and only if $\ker(i_*) \cap K_0(I)_+ = \{0\}$.
	\end{theorem}
	
	\subsection{Direct limits}
	
	Here we establish our notation for direct limits of groups (for details, see for example \cite[Chapter 6]{bluebook}). Suppose that $(G_n)_n$ is a sequence of groups and $\alpha_n : G_n \to G_{n+1}$ are homomorphisms. We denote the direct limit by $G_\infty \coloneqq  \varinjlim(G_n,\alpha_n)$ and we denote the inclusion of each $G_n$ into $G_\infty$ by $\alpha_{n,\infty}$. If $G_n = G$ for all $n \in \N$, then for $x \in G$ we often write $[x]_n$ for $\alpha_{n,\infty}(x)$. 
	
	
	We will use standard results about direct limits of groups without proof, including the following lemma showing that direct limits preserve exactness and commutativity of diagrams.
	
	\begin{lemma}\label{lem:direct limits}Let $(G_n)$ and $(H_n)$ be sequences of groups. For each $n\in \Z$, suppose that $\phi_n : G_n \to H_n, \alpha_n : G_n \to G_{n+1}$ and $\beta_n : H_n \to H_{n+1}$ are homomorphisms, and that for each $n$, $\phi_{n+1} \circ \alpha_n = \beta_n \circ \phi_n$. Then each $\alpha_n(\ker(\phi_n)) \subseteq \ker(\phi_{n+1})$, and each $\beta_n(\image(\phi_n)) \subseteq \image(\phi_{n+1})$. Let $G_\infty \coloneqq  \varinjlim(G_n, \alpha_n)$, let $H_\infty := \varinjlim(H_n, \beta_n)$, and let $\alpha_{n,\infty} : G_n \to G_\infty$ and $\beta_{n,\infty} : H_n \to H_\infty$ be the canonical maps. Then there is a unique homomorphism $\phi_\infty : G_\infty \to H_\infty$ such that $\phi_\infty \circ \alpha_{n,\infty} = \beta_{n,\infty} \circ \phi_n$ for all $n$, we have $\ker(\phi_{\infty}) = \bigcup_n \alpha_{n,\infty}(\ker(\phi_n))$, and we have $\image(\phi_\infty) = \bigcup_n(\beta_{n,\infty}(\image(\phi_n))$.
	\end{lemma}

	\section{\texorpdfstring{$K$}{K}-theory for graph algebras}\label{sec:graph algebras}
	Evans computed the $K$-theory of a $2$-graph $C^*$-algebra in \cite{Evans2008}, using Kasparov's spectral sequence for an action of $\Z^k$ on a $C^*$-algebra. Evans' formula shows that the $K_0$-group decomposes as a direct sum in which one summand is isomorphic to the quotient of $\Z\Lambda^0$ by the subgroup generated by the images of $1 - A_1^t$ and $1 - A_2^t$, and the second summand is the intersection of the kernels of the same two maps. The first summand is the subgroup of $K_0$ generated by the classes of the vertex projections, the identification given by $[p_v] \mapsto \delta_v + \image(1 - A_1^t) + \image(1 - A_2^t)$, but the second one is more mysterious. To prove our main theorem, we will need to know that, for a saturated hereditary subset $H$ of $\Lambda^0$, the map in $K$-theory induced by the canonical inclusion of $C^*(H\Lambda)$ into $C^*(\Lambda)$ respects the above direct-sum decomposition. Unfortunately, the direct sum decomposition obtained by Evans is not natural with respect to $C^*$-homomorphisms, because it is obtained from the general fact that homomorphisms from abelian groups to free abelian groups always split.
	
	Consequently, in order to prove our main theorem, we need to re-prove Evans' theorem by iterating the Pimsner--Voiculescu exact sequence, in such a way that we can fit $K_0(C^*(\Lambda))$ into a short-exact sequence with the two subgroups described above, and guarantee that this short exact sequence is natural with respect to inclusions of saturated hereditary subgraphs.
	
	To do this, we first need to spell out the corresponding naturality property of the well-known exact sequence in $K$-theory for graph $C^*$-algebras proved in \cite{PaskRaeburn1996, RaeburnSzymanski2004}; we do this in this section. We will need, both here, and again later in the paper, the following generalisation of a result by Pask and Raeburn.
	
	\begin{lemma}[c.f.\ {\cite[Theorem 4.2.4]{PaskRaeburn1996}}]\label{lem:paskraeburn}
		Let $G$ be an abelian group and $\phi:G \to G$ a homomorphism. Let
		$\phi_\infty: \varinjlim(G,\phi) \to \varinjlim(G,\phi)$ be the homomorphism such that the diagrams
		\[
		\begin{tikzcd}
			G \arrow[r,"\phi_{n,\infty}"] \arrow[d,"\phi"] & \varinjlim(G,\phi) \arrow[d,"\phi_\infty"] \\
			G \arrow[r,"\phi_{n,\infty}"] & \varinjlim(G,\phi)
		\end{tikzcd}
		\]
        commute. Then $\phi_{1,\infty}$ restricts to an isomorphism of $\ker(1-\phi)$ onto $\ker(1-\phi_\infty)$, and there is an isomorphism $\tilde{\phi}_{1,\infty}: \coker(1-\phi) \to \coker(1-\phi_\infty)$ such that
        \[
			\tilde{\phi}_{1,\infty}(g + \image(1 - \phi)) = \phi_{1,\infty}(g) + \image(1 - \phi_\infty)\qquad\text{ for all $g \in G$.}
        \]
		For each $n \in \N$,
		\[
		\tilde{\phi}_{1,\infty}^{-1}(\phi_{n,\infty}(g) + \image(1-\phi_\infty)) = g + \image(1-\phi)
		\qquad\text{ and }\qquad
		\phi_{1,\infty}^{-1}(\phi_{n,\infty}(g)) = g.
		\]
	\end{lemma}
	\begin{proof}
		As usual, we define an equivalence relation $\sim$ on $\prod^\infty_{n=1} G$ by $(g_n) \sim (h_n)$ if $g_n = h_n$ for large $n$, and identify $\varinjlim(G,\phi)$ with
		\[
		\Big\{(g_n)_{n=1}^\infty \in \prod_{n=1}^{\infty}G : \phi(g_n) = g_{n+1} \text{ for large } n\Big\}/{\sim}.
		\]
		We write $e$ for the identity element of $G$. Then
		\[
		\phi_{n,\infty}(g) = [(\underbrace{e,\dots,e}_{n-1},g,\phi(g),\phi^2(g),\dots)] \qquad \text{ for all $n\in \N, g \in G$.}
		\]
		For $(g_n) \in \varinjlim(G,\phi)$, we have
		\[
		\phi_\infty([(g_n)]) = [(\phi(g_n))].
		\]
		
		Since $\phi_{1,\infty} \circ \phi = \phi_\infty \circ \phi_{1,\infty}$ we have $\phi_{1,\infty}(\ker(1 - \phi)) \subseteq \ker(1 - \phi_\infty)$, so $\phi_{1,\infty}$ restricts to a homomorphism $\ker(1 - \phi) \to \ker(1 - \phi_\infty)$. Moreover, $\phi_{1,\infty}((1-\phi)(G)) \subseteq (1 - \phi_\infty)(\varinjlim(G, \phi))$, and so $\phi_{1,\infty}$ induces a homomorphism $\tilde{\phi}_{1,\infty} : \coker(1 - \phi) \to \coker(1 - \phi_\infty)$ that satisfies $\tilde{\phi}_{1,\infty}(g + \image(1 - \phi)) = \phi_{1,\infty}(g) + \image(1 - \phi_\infty)$.
		
		To see that $\tilde{\phi}_{1,\infty}$ is surjective, fix $[(g_n)] \in \varinjlim(G,\phi)$. Then
		\begin{equation}\label{PRstar}
			\begin{split}
				[(g_n)] + \image(1-\phi_\infty) &= [(g_n)] - (1-\phi_\infty)[(g_n)] + \image(1-\phi_\infty) \\ & = \phi_\infty([(g_n)]) + \image(1-\phi_\infty) = [(\phi(g_n))]+ \image(1-\phi_\infty).
			\end{split}
		\end{equation}
		Thus, for $g\in G$ and $n\in \N$,
		\begin{align*}
			\phi_{n,\infty}(g) + \image(1-\phi_\infty) & = [(\underbrace{e,e,\dots,e}_{n-1\text{ terms}},g,\phi(g),\phi^2(g),\dots)] + \image(1-\phi_\infty) \\
			& = [(\underbrace{\phi^n(e),\phi^n(e),\dots,\phi^n(e)}_{n-1 \text{ terms}},\phi^n(g),\phi^{n+1}(g),\phi^{n+2}(g),\dots)] + \image(1-\phi_\infty) \\
			& = [(g,\phi(g),\dots,\phi^{n-1}(g),\phi^n(g),\phi^{n+1}(g),\phi^{n+2}(g),\dots)] + \image(1-\phi_\infty) \\
			& = \tilde{\phi}_{1,\infty}(g+ \image(1-\phi)).
		\end{align*}
		So $\tilde{\phi}_{1,\infty}$ is surjective.
		
		To see that $\phi_{1,\infty}$ is injective, suppose that $\tilde{\phi}_{1,\infty}(g+\image(1-\phi)) = 0$. Then $\phi_{1,\infty}(g) \in \image(1-\phi_\infty)$, so there exists $(h_n) \in \varinjlim(G,\phi)$ such that
		\[
            [(g,\phi(g),\phi^2(g),\dots)] = [(h_n)] - [\phi(h_n))].
        \]
        Hence $\phi^n(g) = h_n -\phi(h_n) \in \image(1-\phi)$ for large $n$. Now
		\[
            g=(g-\phi(g)) + (\phi(g) -\phi^2(g)) + \dots + (\phi^{n-1}(g)-\phi^n(g)) + \phi^n(g) \in \phi^n(g) + \image(1-\phi) 
        \] belongs to $\image(1-\phi)$.
        Hence $g+ \image(1-\phi) = 0$, proving injectivity of $\tilde{\phi}_{1,\infty}$.
		
		To see that $\phi_{1,\infty}$ maps $\ker(1 - \phi)$ to $\ker(1 - \phi_\infty)$, fix $g\in \ker(1-\phi)$. Then
		\[
            (1-\phi_\infty)(\phi_{1,\infty}(g)) = [(\phi^n(g)-\phi^{n+1}(g))] = [\phi^n((1-\phi)(g))] = 0.
        \]
		
		Now to see that $\phi_{1,\infty}(\ker(1 - \phi)) = \ker(1 - \phi_\infty)$, fix $[(g_n)] \in \ker(1-\phi_\infty)$. We have $g_n = \phi(g_n) = g_{n+1}$ for, say, $n \ge N$. In particular, $\phi(g_n) = g_N$ for all $n\ge N$ and thus
		\[
            \phi_{1,\infty}(g_N) = [(g_N,\phi(g_N),\phi^2(g_N),\dots)] = [(g_N,g_N,g_N,\dots)] = [(g_1,\dots,g_N,g_N,g_N,\dots)]= [(g_n)].
        \]
		
		Finally, to see that $\phi_{1,\infty}|_{\ker(1 - \phi)}$ is injective, fix $g\in \ker(1-\phi) \cap \ker(\phi_{1,\infty})$. Then
		\[
            [(g,\phi(g),\phi^2(g),\dots)] = [(e,e,e,\dots)],
        \]
        so that $\phi^n(g)= e$ for some $n$. Thus,
		\[
		g = \phi(g) = \phi^2(g) = \dots = \phi^n(g) = e.\qedhere
		\]
	\end{proof}
	
	The following is a fairly straightforward consequence of Lemma~\ref{lem:paskraeburn}, but we record it separately for repeated use later.
	
	\begin{lemma}\label{lem:PR naturality}
        For $i = 1,2$, let $H^i, G^i, K^i$ be abelian groups and let $\phi^i : G^i \to G^i$ be a homomorphism. Let $\phi^i_\infty : \varinjlim(G^i, \phi^i) \to \varinjlim(G^i, \phi^i)$, and $\phi^i_{n,\infty} : G^i \to \varinjlim(G^i, \phi^i)$ be as in Lemma~\ref{lem:paskraeburn}. Let $\psi : G^1 \to G^2$ be a homomorphism and let $\psi_\infty : \varinjlim(G^1, \phi^1) \to \varinjlim(G^2,\phi^2)$ be the homomorphism such that $\psi_\infty \circ \phi^1_{n,\infty} = \phi^2_{n,\infty} \circ \psi$ for all $n$. Suppose that
		\[\begin{tikzcd}
			0 \arrow[r] & H^1 \arrow[r, "\iota^1"] \arrow[d, "\eta"] & {\varinjlim(G^1,\phi^1)} \arrow[r,"1-\phi^1_\infty"] \arrow[d, "\psi_\infty"] & {\varinjlim(G^1,\phi^1)} \arrow[r, "\pi^1"] \arrow[d, "\psi_\infty"] & K^1 \arrow[r] \arrow[d, "\theta"] & 0 \\
			0 \arrow[r] & H^2 \arrow[r, "\iota^2"]                   & {\varinjlim(G^2,\phi^2)} \arrow[r,"1-\phi^2_\infty"]                   & {\varinjlim(G^2,\phi^2)} \arrow[r, "\pi^2"]                   & K^2 \arrow[r]                     & 0
		\end{tikzcd}\]
		is a commuting diagram of group homomorphisms with exact rows.
		
		For each $i$, let $\phi^i_{1,\infty}| : \ker(1 - \phi^i) \to \ker(1 - \phi^i_\infty)$, and $\tilde{\phi}^i_{1,\infty} : \coker(1 - \phi^i) \to \coker(1 - \phi^i_\infty)$ be the isomorphisms of Lemma~\ref{lem:paskraeburn}. Let $\tilde{\pi}^i : \coker(1 - \phi^i_\infty) \to K^i$ be the isomorphism induced by $\pi^i$. Let $\zeta^i \coloneqq  (\phi^i_{1,\infty}|)^{-1} \circ \iota^i : H^i \to G^i$. Let $q^i : G^i \to \coker(1-\phi^i)$ be the quotient map and define $\tau^i = \tilde{\pi}^i\circ \tilde{\phi}^i_{1,\infty} \circ q^i : G^i \to K^i$. Then the diagram
		\[\begin{tikzcd}
			0 \arrow[r] & H^1 \arrow[d, "\eta"] \arrow[r, "\zeta^1"] & G^1 \arrow[d, "\psi"] \arrow[r, "1-\phi^1"] & G^1 \arrow[r, "\tau^1"] \arrow[d, "\psi"] & K^1 \arrow[d, "\theta"] \arrow[r] & 0 \\
			0 \arrow[r] & H^2 \arrow[r, "\zeta^2"]                 & G^2 \arrow[r, "1-\phi^2"]                  & G^2 \arrow[r, "\tau^2"]                  & K^2 \arrow[r]                    & 0.
		\end{tikzcd}\]
commutes and has exact rows.
	\end{lemma}
	
	\begin{proof}
			As $\psi_\infty \circ \phi^1_{1,\infty} = \phi^2_{1,\infty}\circ \psi$ by definition of $\psi_\infty$, the diagram
		\[
		\begin{tikzcd}
			&                                                                    & G^1 \arrow[d, "{\phi^1_{1,\infty}}"] \arrow[r, "1-\phi^1"] \arrow[ddd, "\psi"', dotted, bend right]     & G^1 \arrow[d, "{\phi^1_{1,\infty}}"] \arrow[rd, "\tau^1"] \arrow[ddd, "\psi", dotted, bend right] &                                 &    \\
			0 \arrow[r] & H^1 \arrow[r, "\iota^1"] \arrow[d, "\eta"] \arrow[ru, "\zeta^1"] & {\varinjlim(G^1,\phi^1)} \arrow[r, "1-\phi^1_\infty"] \arrow[d, "\psi_\infty"] & {\varinjlim(G^1,\phi^1)} \arrow[r, "\pi^1"] \arrow[d, "\psi_\infty"]             & K^1 \arrow[r] \arrow[d, "\theta"] & 0  \\
			0 \arrow[r] & H^2 \arrow[r, "\iota^2"] \arrow[rd, "\zeta^2"']               & {\varinjlim(G^2,\phi^2)} \arrow[r, "1-\phi^2_\infty"]                       & {\varinjlim(G^2,\phi^2)} \arrow[r, "\pi^2"]                                   & K^2 \arrow[r]                    & 0 \\
			&                                                                    & G^2 \arrow[u, "{\phi^2_{1,\infty}}"'] \arrow[r, "1-\phi^2"]                                          & G^2 \arrow[u, "{\phi^2_{1,\infty}}"'] \arrow[ru, "\tau^2"']                                     &                                 &
		\end{tikzcd}\]
commutes.
		
\end{proof}

\begin{remark}\label{rem:2 squares}
	Whenever we have a commuting diagram of the form
	\[\begin{tikzcd}
		0 \arrow[r] & H^1 \arrow[d, "\eta"] \arrow[r, "\zeta^1"] & G^1 \arrow[d, "\psi"] \arrow[r, "1-\phi^1"] & G^1 \arrow[r, "\tau^1"] \arrow[d, "\psi"] & K^1 \arrow[d, "\theta"] \arrow[r] & 0 \\
		0 \arrow[r] & H^2 \arrow[r, "\zeta^2"]                 & G^2 \arrow[r, "1-\phi^2"]                  & G^2 \arrow[r, "\tau^2"]                  & K^2 \arrow[r]                    & 0,
	\end{tikzcd}\]  with exact rows it is easy to see that $\psi(\ker(1 - \phi^1)) \subseteq \ker(1 - \phi^2)$, there is a homomorphism $\tilde\psi : \coker(1 - \phi^1) \to \coker(1 - \phi^2)$ such that $\tilde\psi(g + \image(1 - \phi^1)) = \psi(g) + \image(1 - \phi^2)$ for all $g \in G^1$, and the diagrams
\[\begin{tikzcd}
	H^1 \arrow[d, "\eta"] \arrow[r, "\zeta^1"] &[2em] \ker(1-\phi^1) \arrow[d, "\psi"] &  & \coker(1-\phi^1) \arrow[d, "\tilde{\psi}"] \arrow[r, "\tau^1"] &[2em] K^1 \arrow[d, "\theta"] \\
	H^2 \arrow[r, "\zeta^2"]                   & {\ker(1-\phi^2),}                &  & \coker(1-\phi^2) \arrow[r, "\tau^2"]                           & {K^2}
\end{tikzcd}\]
commute. We will use this extensively, in some cases when the $G^i$ are themselves direct limits.
\end{remark}

	Lemma~\ref{lem:PR naturality} leads to the following description of the $K$-theory of a $1$-graph algebra. The first part of the statement comes directly from \cite{PaskRaeburn1996, RaeburnSzymanski2004}; the second part is well known, but does not appear to have been recorded in this form.
	
	\begin{theorem}[{\cite{PaskRaeburn1996, RaeburnSzymanski2004}}]\label{thm:1-graph K-th}
		Let $\Lambda$ be a row-finite $1$-graph with no sources, and let $A$ denote its adjacency matrix, regarded as an endomorphism of $\Z \Lb^0$. Let $\phi_0 : \Z\Lambda^0 \to K_0(C^*(\Lambda))$ be the homomorphism such that $\phi_0(\delta_v) = [p_v]$ for all $v \in \Lambda^0$. There exists a homomorphism $\phi_1 : K_1(C^*(\Lambda)) \to \Z\Lambda^0$ that makes the following sequence exact:
		\[
		\begin{tikzcd}
			0 \arrow[r] & K_1(C^*(\Lambda)) \arrow[r,"\phi_1"] & \Z\Lambda^0 \arrow[r, "1 - A^t"] &
			\Z\Lambda^0 \arrow[r, "\phi_0"] & K_0(C^*(\Lambda)) \arrow[r] & 0.
		\end{tikzcd}
		\]
		This sequence is natural in the following sense: if $H \subseteq \Lambda^0$ is hereditary, then $A^t$ restricts to an endomorphism $A^t_H$ of $\Z H$, and the inclusion map $\iota_H : C^*(H\Lambda) \to C^*(\Lambda)$ makes the diagram
		\[
		\begin{tikzcd}
			0 \arrow[r] & K_1(C^*(H\Lambda)) \arrow[r,"\phi_{H,1}"] \arrow[d,"(\iota_H)_*"]
			& \Z H \arrow[r, "1 - A^t_H"] \arrow[d, hook]
			& \Z H \arrow[r, "\phi_{H,0}"] \arrow[d,hook]
			& K_0(C^*(H\Lambda)) \arrow[r] \arrow[d, "(\iota_H)_*"]
			& 0 \\
			0 \arrow[r] & K_1(C^*(\Lambda)) \arrow[r,"\phi_1"]
			& \Z\Lambda^0 \arrow[r, "1 - A^t"]
			& \Z\Lambda^0 \arrow[r, "\phi_0"]
			& K_0(C^*(\Lambda)) \arrow[r]
			& 0
		\end{tikzcd}
		\]
		commute.
	\end{theorem}
	\begin{proof}
		That $K_1(C^*(\Lambda)) \cong \ker(1 - A^t)$ and $K_0(C^*(\Lambda)) \cong \coker(1 - A^t)$ is precisely \cite[Theorem~7.16]{Raeburn2005} (the original result is in \cite{RaeburnSzymanski2004}, but \cite{Raeburn2005} matches our edge-direction conventions). But to justify that the map $\phi_0$ is as claimed and to establish that the second diagram commutes, we must recap the key points of the proof.
		
		The crossed product $C^*(\Lambda) \rtimes_\gamma \T$ by the gauge action $\gamma$ fits into a dual Pimsner--Voiculescu exact sequence
		\[
		\begin{tikzcd}
			0 \arrow[r] &[-1em] K_1(C^*(\Lambda)) \arrow[r]
			& K_0(C^*(\Lambda) \rtimes_\gamma \T) \arrow[r, "1 - \hat{\gamma}_*^{-1}"]
			& K_0(C^*(\Lambda) \rtimes_\gamma \T) \arrow[r, "\Xi"]
			& K_0(C^*(\Lambda)) \arrow[r]
			&[-1em] 0
		\end{tikzcd}
		\]
		(see \cite[page~66]{Raeburn2005}). We need some detail about the map $\Xi$. For this, we write $\chi_n$ for the $n$th spectral projection $z \mapsto z^n$ in $C^*(\T) \cong C_0(\Z)$. The argument of Lemma~4.2.2 of \cite{PaskRaeburn1996} shows that for a projection $p$ in $C^*(\Lambda)$, the homomorphism $\Xi$ carries the class of $i_\T(\chi_n)i_{C^*(\Lambda)}(p)$ to $[p] \in K_0(C^*(\Lambda))$.
		
		By \cite[Lemma~7.10]{Raeburn2005} there is an isomorphism $C^*(\Lambda) \rtimes_\gamma \T \to C^*(\Lambda \times_d \Z)$ that carries $i_\T(\chi_n)i_{C^*(\Lambda)}(s_\lb)$ to $s_{(\lb, n)}$, and this isomorphism intertwines the dual action $\hat{\gamma}$ with the inverse of the action $\lt_*$ on $C^*(\Lambda \times_d \Z)$ induced by left translation in the second coordinate on $\Lambda \times_d \Z$. Hence we obtain an exact sequence
\[
\begin{tikzcd}
	0 \arrow[r] &[-1em] K_1(C^*(\Lambda)) \arrow[r]
	& K_0(C^*(\Lambda \times_d \Z)) \arrow[r, "1 - \lt_*"]
	& K_0(C^*(\Lambda \times_d \Z)) \arrow[r, "P"]
	& K_0(C^*(\Lambda)) \arrow[r]
	&[-1em] 0
\end{tikzcd}
\]
in which the homomorphism $P$ satisfies $P([p_{(v,n)}]) = [p_v]$ for all $v \in \Lambda^0$ and $n \in \Z$.
		
		Fix a hereditary subset $H$ of $\Lambda^0$. As in Lemma~\ref{lem:saturated product}, it is routine to see that $H \times \Z$ is a hereditary subset of $(\Lambda \times_d \Z)^0$, and that $(H \times \Z)(\Lambda \times_d \Z)$ is identical to $H\Lambda \times_d \Z$. The inclusion $\iota_H : C^*(H\Lambda) \to C^*(\Lambda)$ described in Section~\ref{subsec:H<->I} is gauge equivariant, so it induces an inclusion $\iota_H \times 1 : C^*(H\Lambda) \rtimes_{\gamma^H} \T \to C^*(\Lambda) \rtimes_\gamma \T$. Direct calculation on generators shows that the isomorphisms $C^*(H\Lambda) \rtimes_{\gamma^H} \T \to C^*(H\Lambda \times_d \Z)$ and $C^*(\Lambda) \rtimes_\gamma \T \to C^*(\Lambda \times_d \Z)$ intertwine $\iota_H \times 1$ with $\iota_{H \times \Z}$.
		
		As discussed in the first paragraph of Section~3 of \cite{RaeburnSzymanski2004}, the dual Pimsner--Voiculescu sequence is natural with respect to $\T$-equivariant homomorphisms. It therefore follows from the preceding two paragraphs that $\iota_H : C^*(H\Lambda) \to C^*(\Lambda)$ induces a commuting diagram
		\[
		\begin{tikzcd}[nodes={inner sep=3pt}]
			0 \arrow[r] &[-1.25em] K_1(C^*(H\Lambda)) \arrow[r] \arrow[d,"(\iota_H)_*"]
			&[-1em] K_0(C^*(H\Lambda \times_d \Z)) \arrow[r, "1 - \lt_*"] \arrow[d, "(\iota_{H \times \Z})_*"]
			& K_0(C^*(H\Lambda \times_d \Z)) \arrow[r, "P_H"] \arrow[d, "(\iota_{H \times \Z})_*"]
			&[-0.5em] K_0(C^*(H\Lambda)) \arrow[r] \arrow[d, "(\iota_H)_*"]
			&[-1.25em] 0 \\
			0 \arrow[r] & K_1(C^*(\Lambda)) \arrow[r]
			& K_0(C^*(\Lambda \times_d \Z)) \arrow[r, "1 - \lt_*"]
			& K_0(C^*(\Lambda \times_d \Z)) \arrow[r, "P"]
			& K_0(C^*(\Lambda)) \arrow[r]
			& 0
		\end{tikzcd}
		\]
		of exact sequences.
		
		Corollary~7.14 of \cite{Raeburn2005} shows that there is an isomorphism $\xi:\varinjlim(\Z\Lambda^0, A^t) \to K_0(C^*(\Lambda \times_d \Z))$ that carries $(A^t)_{n,\infty}(\delta_v)$ to $[p_{(v,n)}]$, and direct calculation on generators shows that composing $P : K_0(C^*(\Lambda \times_d \Z)) \to K_0(C^*(\Lambda))$ with this isomorphism yields a homomorphism $(\phi_0)_\infty : \varinjlim(\Z\Lambda^0, A^t) \to K_0(C^*(\Lambda))$ satisfying $(\phi_0)_\infty((A^t)_{n,\infty}(\delta_v)) = [p_v]$. (We denote the corresponding map for $H\Lambda$ by $(\phi_{H,0})_\infty$.) The argument preceding \cite[Corollary~7.14]{Raeburn2005} also shows that the isomorphism $\xi:\varinjlim(\Z\Lambda^0, A^t) \to K_0(C^*(\Lambda \times_d \Z))$ intertwines $(A^t)_\infty$ with $\lt_*$; that is,
		\begin{equation}\label{eq:1-graph translation}\xi \circ (A^t)_\infty =  \lt_* \circ \,\xi.
		\end{equation}
		The inclusion $\psi : \Z H \hookrightarrow \Z\Lambda^0$ intertwines $A_H^t$ with $A^t$. So we obtain a commuting diagram
		\[
		\begin{tikzcd}
			0 \arrow[r] &[-1em] K_1(C^*(H\Lambda)) \arrow[r] \arrow[d,"(\iota_H)_*"]
			& \varinjlim(\Z H, A^t_H) \arrow[r, "1 - (A^t_H)_\infty"] \arrow[d, "\psi_\infty"]
			&[1em] \varinjlim(\Z H, A^t_H) \arrow[r, "(\phi_{H,0})_\infty"] \arrow[d, "\psi_\infty"]
			& K_0(C^*(H\Lambda)) \arrow[r] \arrow[d, "(\iota_H)_*"]
			&[-1em] 0 \\
			0 \arrow[r] & K_1(C^*(\Lambda)) \arrow[r]
			& \varinjlim(\Z\Lambda^0, A^t) \arrow[r, "1 - (A^t)_\infty"]
			& \varinjlim(\Z\Lambda^0, A^t) \arrow[r, "(\phi_{0})_\infty"]
			& K_0(C^*(\Lambda)) \arrow[r]
			& 0
		\end{tikzcd}
		\]
		of exact sequences.
		
		The result now follows from an application of Lemma~\ref{lem:PR naturality}.
	\end{proof}

	\section{\texorpdfstring{$K$}{K}-theory computations for \texorpdfstring{$2$}{2}-graph algebras}\label{sec: k theory main}

This section is devoted to computing the $K$-theory of a $2$-graph $C^*$-algebra $C^*(\Lambda)$ in such a way that we can keep track of the homomorphism of $K_0$-groups induced by the inclusion of the subalgebra $C^*(H\Lambda)$ associated to a given hereditary set $H \subseteq \Lambda^0$. This requires a few steps.

	\subsection{Commutative diagrams involving the \texorpdfstring{$K$}{K}-theory of \texorpdfstring{$C^*(\Lambda\times_{d_1}\Z)$}{the skew product}}\label{subsec:1-dim skew}

	Our first step is to compute the $K$-theory of the graph algebra $C^*(\Lambda\times_{d_1}\Z)$ via an exact sequence and describe morphisms between instances of this exact sequence induced by two maps between copies of $C^*(\Lambda \times_{d_1} \Z)$. We compile the main results of this section in Theorem~\ref{thm:main k theory theorem} and Corollary~\ref{cor: main k theory theorem}. Recall that given a group homomorphism $\phi: G\to G$, the notation $[g]_n$ refers to $\phi_{n,\infty}(g) \in \varinjlim(G,\phi)$. We begin with two preliminary results about the inclusion induced by the hereditary subgraph $H\Lb$ and direct limits respectively.
	
	\begin{lemma}\label{lem:hereditary direct limit}
	Let $\Lambda$ be a row-finite $2$-graph with no sources. Let $H$ be a hereditary subset of $\Lambda^0$. Then there is an injective homomorphism $\varinjlim(\Z H, A_{H,1}^t) \hookrightarrow \varinjlim(\Z\Lambda^0, A_1^t)$ that carries $[\delta_v]_n$ to $[\delta_v]_n$ for $v \in H$.
	\end{lemma}
\begin{proof}
	Because $H$ is a hereditary subset, both $A_{H,1}^t$ and $A_1^t$ map $\Z H$ to $\Z H$.
	In other words, the following diagram commutes:
	\[
	\begin{tikzcd}
		\Z H \arrow[d, "{A_{H,1}^t}"] \arrow[r, hook] & \Z\Lb^0 \arrow[d, "A_1^t"] \\
		\Z H \arrow[r, hook]                          & \Z\Lb^0.
	\end{tikzcd}\]
	Since $\Z H \hookrightarrow \Z\Lb^0$ is an injection, Lemma~\ref{lem:direct limits} implies that the homomorphism $\varinjlim(\Z H, A_{H,1}^t) \to \varinjlim(\Z \Lb^0,A_1^t)$ induced by $i$ is injective.
\end{proof}

The next result shows that the second adjacency matrix $A_2^t$ induces an endomorphism of the direct limits appearing in Lemma~\ref{lem:hereditary direct limit}, and describes its kernel and cokernel.

	\begin{proposition}\label{prop:basic direct limit}
		Let $\Lb$ be a row-finite $2$-graph with no sources. Then there is a homomorphism $A_2^{t,\infty} : \varinjlim(\Z \Lambda^0, A_1^t) \to \varinjlim(\Z \Lambda^0,
			A_1^t)$
			that satisfies $A_2^{t,\infty}[g]_n = [A_2^t g]_n$ for $g \in
			\Z\Lambda^0$. The subgroup $\varinjlim(\ker(1-A_2^t),A_1^t)$ of $\varinjlim(\Z\Lb^0, A_1^t)$ is equal to $\ker(1-A_2^{t,\infty})$. Writing $\tilde{A}_1^t$ for the endomorphism of $\coker(1-A_2^t)$ induced by $A_1^t$, there is an isomorphism \[i_2 : \coker(1-A_2^{t,\infty}) \to \varinjlim(\coker(1-A_2^t),\tilde{A}_1^t)\] such that $i_2([g]_n+\image(1-A^{t,\infty}_2)) = [g+\image(1-A_2^t)]_n$ for all $g \in \Z\Lb^0$.
	\end{proposition}

	\begin{proof}
		The formula for $A_2^{t,\infty}$ is well-defined because $A_1^t$ and $A_2^t$ commute. It clearly defines a homomorphism.
		
		 Fix $g \in \ker(1-A_2^t) \subset \Z\Lambda^0$. Then
     \[
        (1-A_2^{t,\infty})([g]_n) = [(1-A_2^{t})g]_n = [0]_n = 0.
     \]
     Hence $\varinjlim(\ker(1-A_2^t),A_1^t) \subseteq \varinjlim(\Z\Lb^0, A_1^t)$. Fix $h \in \ker(1 - A_2^{t,\infty})$. Then $h = [g]_n$ for some $g \in \Z\Lb^0$ and $n \in \N$. We have $0 = (1 - A_2^{t,\infty})h = [(1-A_2^{t})g]_n$, so there exists $k\ge 0$ such that $0 = (A_1^t)^k(1-A_2^{t})g = (1-A_2^t)(A_1^t)^k g$. Hence $(A_1^t)^k g \in \ker(1-A_2^t)$, and therefore $h = [g]_n = [(A_1^t)^k g]_{n+k} \in \varinjlim(\ker(1-A_2^t),A_1^t)$. This proves that $\varinjlim(\ker(1-A_2^t),A_1^t) = \ker(1-A_2^{t,\infty})$.
		
		Now fix $(1-A_2^{t,\infty})[g]_n \in \image(1-A_2^{t,\infty})$. Then
    \[
        (1-A_2^{t,\infty})[g]_n = [(1-A_2^t)g]_n \in A_1^{t,\infty}(\image(1-A_2^t)).
    \]
    Hence $i_2$ is well defined. To see that $i_2$ is injective, suppose that $[g+\image(1-A_2^t)]_n = 0$. Then there exists $k\ge 0$ such that $(A_1^t)^k g+\image(1-A_2^t) = (\tilde{A}_1^t)^k (g+\image(1-A_2^t)) = 0$. In particular, there exists $h\in \Z\Lambda^0$ such that $(A_1^t)^k g = (1-A_2^t)h$. Hence $[g]_n = [(A_1^t)^k g]_{n+k} = [(1-A_2^t)h]_{n+k} = (1-A_2^{t,\infty})[h]_{n+k} \in \image(1-A_2^{t,\infty})$.
    For surjectivity, just note that each $[g+\image(1-A_2^t)]_n = i_2([g]_n+\image(1-A^{t,\infty}_2))$.
	\end{proof}
	
\begin{theorem}\label{thm:main k theory theorem}
Let $\Lambda$ be a row-finite $2$-graph with no sources. Let $H$ be a
hereditary subset of $\Lambda^0$. With the notation of Lemmas \ref{lem:saturated product}~and~\ref{lem:hereditary direct limit} and Proposition~\ref{prop:basic direct limit},
\begin{enumerate}
\item There exist group homomorphisms $\phi_0 : \varinjlim(\Z\Lambda^0, A_1^t) \to
    K_0(C^*(\Lambda \times_{d_1} \Z))$ and $\phi_{H, 0} : \varinjlim(\Z H, A_{H,1}^t) \to
    K_0(C^*(H\Lambda \times_{d_1} \Z))$ such that $\phi_0([\delta_{v}]_n) = [p_{(v,n)}]$ for all
    $v \in \Lambda^0$ and $\phi_{H, 0}([\delta_{v}]_n) = [p_{(v,n)}]$ for all $v \in H$, and
    there are homomorphisms $\phi_1:K_1(C^*(\Lambda \times_{d_1} \Z)) \to \varinjlim(\Z
    \Lambda^0, A_1^t)$ and $\phi_{H,1} : K_1(C^*(H\Lambda \times_{d_1} \Z)) \to \varinjlim(\Z H,
    A_{H,1}^t)$ such that the diagram
	\[
	\begin{tikzcd}[nodes={inner sep=2pt}]
		0 \arrow[r]
		&[-1.25em] K_1(C^*(H\Lambda \times_{d_1} \Z)) \arrow[r, "\phi_{H,1}"] \arrow[d, "(\iota_{H \times \Z})_*"]
			&[-1.75ex] \varinjlim(\Z H, A_{H,1}^t) \arrow[r, "1 - A_{H,2}^{t,\infty}"] \arrow[d, hook]
			& \varinjlim(\Z H, A_{H,1}^t) \arrow[r, "\phi_{H,0}"] \arrow[d, hook]
			&[-1.75ex] K_0(C^*(H\Lambda \times_{d_1} \Z)) \arrow[r] \arrow[d, "(\iota_{H \times \Z})_*"]
			&[-1.25em] 0\\
			0 \arrow[r]
			& K_1(C^*(\Lambda \times_{d_1} \Z)) \arrow[r, "\phi_{1}"]
			& \varinjlim(\Z \Lambda^0, A_1^t) \arrow[r, "1 - A_2^{t,\infty}"]
			& \varinjlim(\Z \Lambda^0, A_1^t) \arrow[r, "\phi_{0}"]
			& K_0(C^*(\Lambda \times_{d_1} \Z)) \arrow[r]
			& 0
	\end{tikzcd}
	\]
	commutes and has exact rows.
\item Let $A_1^{t,\infty}$ be the map on $\varinjlim(\Z\Lb^0,A_1^t)$ given by $[g]_n\mapsto [A_1^tg]_n$. Using the same maps $\phi_1$ and $\phi_0$ from above, the diagram
	\[\begin{tikzcd}[nodes={inner sep=2pt}]
		0 \arrow[r] &[-1.5em] K_1(C^*(\Lambda\times_{d_1}\Z)) \arrow[r, "{\phi_{1}}"] \arrow[d, "1-(\lt_1)_*"] & {\varinjlim(\Z \Lambda^0,A_{1}^t)} \arrow[r, "{1-A_{2}^{t,\infty}}"] \arrow[d,"1-A_1^{t,\infty}"] & {\varinjlim(\Z \Lambda^0,A_{1}^t)} \arrow[r, "{\phi_{0}}"] \arrow[d,"1-A_1^{t,\infty}"] & K_0(C^*(\Lambda\times_{d_1}\Z)) \arrow[r] \arrow[d, "1-(\lt_1)_*"] &[-1.5em] 0 \\
			0 \arrow[r] &[-1.5em] K_1(C^*(\Lambda\times_{d_1}\Z)) \arrow[r, "\phi_1"]                                       & {\varinjlim(\Z\Lambda^0,A_{1}^t)} \arrow[r, "{1-A_2^{t,\infty}}"]                & {\varinjlim(\Z\Lambda^0,A_{1}^t)} \arrow[r, "\phi_0"]                  & K_0(C^*(\Lambda\times_{d_1}\Z)) \arrow[r]                                 &[-1.5em] 0
	\end{tikzcd}\]
commutes and has exact rows.
\end{enumerate}
\end{theorem}

\begin{remark}\label{rem: thm via naturality}
	Theorem~\ref{thm:main k theory theorem} essentially says that the exact sequence
	\[\begin{tikzcd}
		0 \arrow[r] &[-1.5em] K_1(C^*(\Lambda\times_{d_1}\Z)) \arrow[r, "\phi_1"] & {\varinjlim(\Z\Lambda^0,A_{1}^t)} \arrow[r, "{1-A_2^{t,\infty}}"] & {\varinjlim(\Z\Lambda^0,A_{1}^t)} \arrow[r, "\phi_0"] & K_0(C^*(\Lambda\times_{d_1}\Z)) \arrow[r] &[-1.5em] 0
	\end{tikzcd}\]
is natural with respect to both the inclusion of hereditary subgraph algebras and the action of left-translation.
\end{remark}

We will prove Theorem~\ref{thm:main k theory theorem} later. But first we state the following corollary which explicitly states the main ingredients needed for our second main $K$-theory result, namely Theorem~\ref{thm:inc}.

\begin{corollary}\label{cor: main k theory theorem}
	Let $\Lambda$ be a row-finite $2$-graph with no sources. Let $H$ be a hereditary subset of $\Lambda^0$. Resume the notation of Lemmas \ref{lem:saturated product}~and~\ref{lem:hereditary direct limit}, Proposition~\ref{prop:basic direct limit} and Theorem~\ref{thm:main k theory theorem}.
	\begin{enumerate}
		\item Let $\phi_{1}| : \varinjlim(\ker(1-A_2^t),A_1^t) \to K_1(C^*(\Lambda\times_{d_1}\Z))$ be the restriction of $\phi_1$ to $\ker(1-A_2^{t,\infty}) = \varinjlim(\ker(1-A_2^t),A_1^t)$. Let $\tilde{\phi}_0 : \coker(1-A_2^{t,\infty}) \cong \varinjlim(\coker(1-A_2^t),\tilde{A}_1^t) \to K_0(C^*(\Lambda \times_{d_1} \Z))$ be the map induced by $\phi_0$. The diagrams
		\[
		\begin{tikzcd}
			\varinjlim(\ker(1-A_{H,2}^t),A_{H,1}^t) \arrow[r,"\iota"] \arrow[d,"(\phi_{H,1}|)^{-1}"] & \varinjlim(\ker(1-A_2^t),A_1^t) \arrow[d,"(\phi_1|)^{-1}"] \\
			K_1(C^*(H\Lb \times_{d_1}\Z)) \arrow[r,"\iota_*"] & K_1(C^*(\Lb \times_{d_1}\Z))
		\end{tikzcd}
		\]
		and
		\[
		\begin{tikzcd}
			\varinjlim(\coker(1-A_{H,2}^t),\tilde{A}_{H,1}^t) \arrow[r,"\iota"] \arrow[d,"\tilde{\phi}_{H,0}"] & \varinjlim(\coker(1-A_2^t),\tilde{A}_1^t) \arrow[d,"\tilde{\phi}_0"] \\ K_0(C^*(H\Lb \times_{d_1}\Z)) \arrow[r,"\iota_*"] & K_0(C^*(\Lb \times_{d_1}\Z))
		\end{tikzcd}
		\] commute. In both diagrams, the vertical maps are isomorphisms. We have $\tilde{\phi}_0([\delta_v + \image(1-A_2^t)]_n) = [p_{(v,n)}]$ and likewise for $\tilde{\phi}_{0,H}$.
		\item The diagrams
		\[
		\begin{tikzcd}
			\varinjlim(\ker(1-A_2^t),A_1^t) \arrow[r,"1-A_1^{t,\infty}"] \arrow[d,"(\phi_1|)^{-1}"] & \varinjlim(\ker(1-A_2^t),A_1^t) \arrow[d,"(\phi_1|)^{-1}"] \\
			K_1(C^*(\Lb \times_{d_1}\Z)) \arrow[r,"1-(\lt_1)_*"] & K_1(C^*(\Lb \times_{d_1}\Z))
		\end{tikzcd}
		\]
		and
		\[
		\begin{tikzcd}
			\varinjlim(\coker(1-A_2^t),\tilde{A}_1^t) \arrow[r,"1-\tilde{A}_1^{t,\infty}"] \arrow[d,"\tilde{\phi}_0"] & \varinjlim(\coker(1-A_2^t),\tilde{A}_1^t) \arrow[d,"\tilde{\phi}_0"] \\ K_0(C^*(\Lb \times_{d_1}\Z)) \arrow[r,"1-(\lt_1)_*"] & K_0(C^*(\Lb \times_{d_1}\Z))
		\end{tikzcd}
		\] commute, as do the analogous diagrams for $H\Lb$.
	\end{enumerate}
\end{corollary}

\begin{proof}
	This is a simple application of Remark~\ref{rem:2 squares} and Proposition~\ref{prop:basic direct limit} to the diagrams of Theorem~\ref{thm:main k theory theorem}. The explicit mapping for $\tilde{\phi}_0$ follows from that of $\phi_0$.
\end{proof}
	
	We will now focus on proving Theorem~\ref{thm:main k theory theorem} by splitting it into three main labelled steps that mirror Remark~\ref{rem: thm via naturality}: first we obtain what is needed for the long exact sequence in Remark~\ref{rem: thm via naturality}, then we demonstrate that this sequence is natural with respect to inclusion of hereditary subgraph algebras, and finally we show that it is natural with respect to left-translation. Before we proceed to these three steps, we first need some preliminary work that will allow us to break our steps down to results about smaller pieces of the $C^*$-algebras involved corresponding to levels in various skew-product graphs.
	
	\subsection{Smaller slices: the subalgebras \texorpdfstring{$B^\Lambda_n$}{Bn}, \texorpdfstring{$C^\Lambda_{n,m}$}{Cn,m} and \texorpdfstring{$C^\Lambda_n$}{Cn}}\label{subsec:slices}

Let $\Lambda$ be a row-finite $2$-graph with no sources. For $n\in\Z$, let
\begin{align*}B^\Lambda_n & \coloneqq  \overline{\spann}\{s_\alpha s_\beta^* : \alpha,\beta \in \Lambda\times_{d_1}\Z, s(\alpha) = s(\beta) \in \Lambda^0\times \{n\}\} \\
	& = \overline{\spann}\{s_{(\mu,n-d_1(\mu))}s_{(\nu,n-d_1(\nu))}^* : \mu,\nu \in \Lb\} \subseteq C^*(\Lb\times_{d_1}\Z)
\end{align*}
and
\begin{align*}C^\Lambda_n & \coloneqq  \overline{\spann}\{s_\alpha s_\beta^* : \alpha,\beta \in \Lambda\times_{d}\Z^2, s(\alpha) = s(\beta) \in \Lambda^0\times \{n\}\times\Z\} \\
	& = \overline{\spann}\{s_{(\mu,p-d(\mu))}s_{(\nu,p-d(\nu))}^* : \mu,\nu \in \Lb, p \in \{n\} \times \Z\} \subseteq C^*(\Lb\times_d \Z^2).
\end{align*}

For $m,n\in\Z$, let
\begin{align*}C^\Lambda_{n,m} & \coloneqq  \overline{\spann}\{s_\alpha s_\beta^* : \alpha,\beta \in \Lambda\times_{d}\Z^2, s(\alpha) = s(\beta) \in \Lambda^0\times \{n\}\times\{m\}\} \\
	& = \overline{\spann}\{s_{(\mu,(n,m)-d(\mu))}s_{(\nu,(n,m)-d(\nu))}^* : \mu,\nu \in \Lb\} \subseteq C^*(\Lb\times_d \Z^2).
\end{align*}

The following general lemma will allow us to describe each of the $B^\Lambda_n$, the $C^\Lambda_n$ and the $C^\Lambda_{n,m}$ and the relationships between them. For the statement of the result,  given $J \subseteq \{1, \dots, k\}$, write $\N^J$ for the submonoid of $\N^k$ generated by $\{e_j : j \in J\}$ and let $\Z^J = \N^J - \N^J$ be the enveloping subgroup of $\Z^k$. Let $\pi_J : \Z^k \to \Z^J$ be the projection map.

\begin{lemma}\label{lem:skew subalgebras}
Let $\Lambda$ be a row-finite $k$-graph with no sources. Fix $K \subseteq J \subseteq \{1, \dots, k\}$. Let $d_J \coloneqq  \pi_J \circ d : \Lambda \to \N^J$. Consider the skew-product graph $\Lambda \times_{d_J} \Z^J$. For $X \subseteq \Z^J$, let
\[
A_X \coloneqq  \overline{\spann}\{s_{(\mu, n - d_J(\mu))} s^*_{(\nu, n-d_J(\nu))} : s(\mu) = s(\nu)\text{ and }n \in X\} \subseteq C^*(\Lb\times_{d_J} \Z^J).
\]
For $\mu,\nu \in \Lambda$, $n \in X$ and $p \in \N^k$, we have
\begin{equation}\label{eq:CK in skew prod}
s_{(\mu, n - d_J(\mu))} s^*_{(\nu, n-d_J(\nu))} = \sum_{\lb \in (s(\mu), n-d_J(\mu))(\Lambda \times_{d_J} \Z^J)^p} s_{(\mu\lb, n - d_J(\mu))} s^*_{(\nu\lb, n-d_J(\nu))} \in A_{X + \pi_J(p)};
\end{equation}
in particular $A_X \subseteq A_{X + \pi_J(p)}$. For $q \in \Z^J$ the space $A_{q + \Z^K}$ is a $C^*$-subalgebra of $C^*(\Lambda\times_{d_J}\Z^J)$. For $q \in \Z^J$ and $p \in \N^J$, we have $A_{q + \Z^K} \subseteq A_{q + p + \Z^K}$. Finally, each $A_{q + \Z^K} = \overline{\bigcup_{p \in \Z^K} A_{\{q + p\}}}$.
\end{lemma}
\begin{proof}
The equality in~\eqref{eq:CK in skew prod} is a straight application of the Cuntz--Krieger relation~(4) (see Subsection~\ref{ss:cuntzkrieger}). Membership of $A_{X + \pi_J(p)}$ is because each $n - d_J(\mu) = n + d_J(\lb) - d_J(\mu\lb) = (n + \pi_J(p)) - d_J(\mu\lb)$.

Fix $\alpha,\beta,\mu,\nu \in \Lb \times_{d_J}\Z^J$ with $s(\alpha)=s(\beta), s(\mu)=s(\nu) \in \Lb^0\times (q + \Z^K)$. Then by the Cuntz--Krieger relations, for $m = d(\beta)\vee d(\mu)$ we have
	\[
    s_\alpha s_\beta^* s_\mu s_\nu^* =  s_\alpha\left(\sum_{\beta\tau = \mu\rho \in \Lb^m}s_{\tau} s_{\rho}^*\right)s_\nu^* = \sum_{\beta\tau = \mu\rho \in \Lb^m}s_{\beta\tau} s_{\nu\rho}^*.
     \]
So it suffices to fix $\beta\tau = \mu\rho \in \Lambda^m$ and show that $s(\tau) \in \Lb^0 \times (q + \Z^K)$. We have $r(\beta) = r(\mu)$ and $s(\beta),s(\mu) \in \Lb^0\times (q + \Z^K)$. So $d_{J\setminus K}(\beta) = d_{J\setminus K}(\mu)$. Hence $d_J(\beta)\vee d_J(\mu) = d_{J\setminus K} (\beta) + (d_K(\beta)\vee d_K(\mu))$. Since $d_J(\tau) = (d_J(\beta) \vee d_J(\mu)) - d_J(\beta)$, we deduce that $d_{J \setminus K}(\tau) = 0$. Since $r(\tau) = s(\beta) \in \Lb^0\times (q + \Z^K)$, it follows that $s(\tau) \in \Lb^0 \times (q + \Z^K)$. Hence $A_{q + \Z^K}$ is closed under multiplication. It is clearly closed under adjoints, and it is a closed linear subspace by definition, so it is a $C^*$-subalgebra. We have $A_{q + \Z^K} = \overline{\bigcup_{p \in \Z^K} A_{q + p}}$ because both are densely spanned by $\{s_\alpha s^*_\beta : s(\alpha) = s(\beta) \in \Lb^0 \times (q+\Z^K)\}$.
\end{proof}

\begin{corollary}\label{cor:subalgebras} Let $\Lambda$ be a row-finite $2$-graph with no sources. Then
	\begin{enumerate}
		\item each $B^\Lambda_n \subseteq B^\Lambda_{n+1}$ is a $C^*$-subalgebra of $C^*(\Lambda \times_{d_1} \Z)$, and $C^*(\Lb\times_{d_1}\Z) = \overline{\bigcup_{n \in \Z} B^\Lambda_n}$;
		\item each $C^\Lambda_n \subseteq C^\Lambda_{n+1}$ is a $C^*$-subalgebra of $C^*(\Lb\times_{d}\Z^2)$, and $C^*(\Lb\times_{d}\Z^2)  = \overline{\bigcup_{n \in \Z} C^\Lambda_n}$; and
		\item each $C^\Lambda_{n,m}$ is a $C^*$-subalgebra of $C^\Lambda_{n+1,m} \cap C^\Lambda_{n,m+1}$, and $C^\Lambda_n = \overline{\bigcup_{m \in \Z} C^\Lambda_{n,m}}$.
	\end{enumerate}
\end{corollary}
\begin{proof}
In the notation of Lemma~\ref{lem:skew subalgebras}, putting $J = \{1\}$ and $K = \emptyset$ gives $B^\Lambda_n = A_{n + \Z^K} = A_{\{n\}}$; putting $J = \{1,2\}$ and $K = \{2\}$ gives $C^\Lambda_n = A_{(n,0) + \Z^K}$; and putting $J = \{1,2\}$ and $K = \emptyset$ gives $C^\Lambda_{m,n} = A_{(n,m) + \Z^K} = A_{\{(n,m)\}}$. So all three assertions follow from Lemma~\ref{lem:skew subalgebras}.
\end{proof}

Continuity of $K$-theory gives us the following corollary.

\begin{corollary}\label{cor:inductivelimitsubalgebras} Let $\Lambda$ be a row-finite $2$-graph with no sources.
	\begin{enumerate}
		\item For $n \in \Z$, let $\iota_n : B^\Lambda_n \to B^\Lambda_{n+1}$ be the inclusion of Corollary~\ref{cor:subalgebras}. Then \[K_*(C^*(\Lb \times_{d_1}\Z)) \cong \varinjlim(K_*(B^\Lambda_n),(\iota_n)_*).\]
		\item For $m,n\in\Z$, let $\iota_m : C^\Lambda_{n,m} \to C^\Lambda_{n,m+1}$ be the inclusion of Corollary~\ref{cor:subalgebras}. Then \[K_*(C^\Lambda_n) \cong \varinjlim(K_*(C^\Lambda_{n,m}),(\iota_m)_*).\]
	\end{enumerate}
\end{corollary}

\subsection{Step 1}\label{subsec:step1}

We now state our first step, but the proof must wait until we have the tools needed to prove it.

\begin{proposition}\label{prop: Bn step 1}
Let $\Lambda$ be a row-finite $2$-graph with no sources. Resume the notation of Corollary~\ref{cor:subalgebras}. For each $n \in \Z$ let $\phi_{n,0} : \Z\Lb^0 \to K_0(B^\Lambda_n)$ be the homomorphism such that $\phi_{n,0}(\delta_v) = [p_{(v,n)}]$ for all $v \in \Lb^0$. Then there exist homomorphisms $\phi_{n,1} : K_1(B^\Lambda_n) \to \Z\Lb^0$, $n \in \Z$ such that the diagrams
\[
	\begin{tikzcd}
		0 \arrow[r] & K_1(B^\Lambda_n) \arrow[r,"\phi_{n,1}"] \arrow[d, "(\iota_n)_*"] & \Z\Lambda^0 \arrow[r, "1-A_2^t"] \arrow[d, "A_1^t"] & \Z\Lambda^0 \arrow[r,"\phi_{n,0}"] \arrow[d, "A_1^t"] & K_0(B^\Lambda_n) \arrow[r] \arrow[d, "(\iota_n)_*"] & 0 \\
		0 \arrow[r] & K_1(B^\Lambda_{n+1}) \arrow[r,"\phi_{n+1,1}"]                      & \Z\Lambda^0 \arrow[r, "1-A_2^t"]                    & \Z\Lambda^0 \arrow[r,"\phi_{n+1,0}"]                    & K_0(B^\Lambda_{n+1}) \arrow[r]                      & 0
	\end{tikzcd}\]
commute.
\end{proposition}

Now to develop our toolkit. While our application will be to $2$-graphs, the following result is easy to state for general $k$-graphs. We adopt again the notation used in Lemma~\ref{lem:skew subalgebras}. In addition, for $K \subseteq J \subseteq \{1, \dots, k\}$, we write $\lt_K$ for the action of $\Z^k$ on $C^*(\Lb\times_{d_J}\Z^J)$ such that $\lt_K(m)(s_{(\mu,p)}) = s_{(\mu,p-m)}$ for $\mu \in \Lb$, $p \in \Z^J$ and $m \in \Z^K \subseteq \Z^J$.

\begin{lemma}\label{lem:fullcorners}
Let $\Lb$ be a row-finite $k$-graph with no sources. Suppose that $I_1$ and $I_2$ are disjoint subsets of $\{1,\dots,k\}$ and let $I \coloneqq  I_1\cup I_2$. Let $\iota_1 : C^*(\Lb\times_{d_I}\Z^I) \to C^*(\Lb \times_{d_{I}}\Z^{I})\rtimes_{\lt_{I_2}}\Z^{I_2}$ be the canonical inclusion. Let $Q_0 \in \mathcal{M}(C^*(\Lb \times_{d_{I}}\Z^{I})\rtimes_{\lt_{I_2}}\Z^{I_2})$ be the projection
\[
    \sum_{v\in\Lb^0,\, p \in \Z^{I_1}} \iota_1(s_{(v,p,0)}).
\]
Let $\iota_2 : \Z^{I_2} \to \mathcal{UM}(C^*(\Lb \times_{d_{I}}\Z^{I})\rtimes_{\lt_{I_2}}\Z^{I_2})$ be the canonical unitary representation. Then there is an isomorphism
\[
    \psi: C^*(\Lb\times_{d_{I_1}}\Z^{I_1}) \to Q_0 (C^*(\Lb \times_{d_{I}}\Z^{I})\rtimes_{\lt_{I_2}}\Z^{I_2}) Q_0
\]
such that $\psi(s_{(\mu,n)}) = \iota_1(s_{(\mu,n,0)})\iota_2(d_{I_2}(\mu))$ for all $\mu \in \Lb,n \in \Z^{I_1}$.
\end{lemma}

\begin{proof}
The map $t : (\mu,n)\mapsto \iota_1(s_{(\mu,n,0)})\iota_2(d_{I_2}(\mu))$ is a Cuntz--Krieger $(\Lb\times_{d_{I_1}}\Z^{I_1})$-family in $C^*(\Lb \times_{d_{I}}\Z^{I})\rtimes_{\lt_{I_2}}\Z^{I_2}$. The dual of $\lt_{I_2}$ is an action $\alpha$ of $\mathbb{T}^{I_2}$ on $C^*(\Lb \times_{d_{I}}\Z^{I})\rtimes_{\lt_{I_2}}\Z^{I_2}$ such that
\[
    \alpha_z(\iota_1(s_{(\mu,p,q)})\iota_2(m)) = z^m\iota_1(s_{(\mu,p,q)})\iota_2(m), \quad \text{ for all $z \in \mathbb{T}^{I_2}$.}
\]
So by the gauge invariant uniqueness theorem \cite[Theorem 3.4]{KP2000}, $s_{(\mu,n)}\mapsto t_{(\mu,n)}$ is an injection from $C^*(\Lb\times_{d_{I_1}}\Z^{I_1})$ into $C^*(\Lb \times_{d_{I}}\Z^{I})\rtimes_{\lt_{I_2}}\Z^{I_2}$. We have
	\begin{align*}
		C^*(\Lb \times_{d_{I}}\Z^{I})& \rtimes_{\lt_{I_2}}\Z^{I_2} \\ & = \overline{\spann}\{\iota_1(s_{(\mu,p,q)}s_{(\nu,r,s)}^*)\iota_2(u) : \mu,\nu \in \Lb, p,r, \in \Z^{I_1}, q,s,u \in \Z^{I_2}\}
	\end{align*} so that
	\begin{align*}
		Q_0 \big(C^*(\Lb&{} \times_{d_{I}}\Z^{I}) \rtimes_{\lt_{I_2}}\Z^{I_2}\big) Q_0 \\ & = \overline{\spann}\big\{\iota_1(s_{(\mu,p,0)})\iota_2(d_{I_2}(\mu)-d_{I_2}(\nu))\iota_1(s_{(\nu,r,0)}^*) : \mu,\nu \in \Lb, p,r, \in \Z^{I_1}\big\} \\ & = C^*(\{t_{(\mu,n)} : (\mu,n) \in \Lb\times_{d_{I_1}}\Z^{I_1}\}).\qedhere
	\end{align*}
\end{proof}

Applying Lemma~\ref{lem:fullcorners} when $\Lb$ is a $2$-graph with $I_1 = \{1\}$ and $I_2 = \{2\}$ gives most of the following corollary.

\begin{corollary}\label{cor:Bn k theory iso Cncrossedproduct}
Let $\Lambda$ be a row-finite $2$-graph with no sources. Let $\iota_1 : C^*(\Lambda \times_d \Z) \to C^*(\Lambda \times_d \Z) \times_{\lt_2} \Z$ be the canonical inclusion. Let $Q_0 \in \mathcal{M}(C^*(\Lb \times_{d}\Z^{2})\rtimes_{\lt_{2}}\Z)$ be the projection
\[
    Q_0 = \sum_{v\in\Lb^0,\, p \in \Z} \iota_1(s_{(v,(p,0))}).
\]
Then there is an isomorphism
\[
\psi : C^*(\Lb\times_{d_1}\Z) \to Q_0(C^*(\Lb\times_d \Z^2)\rtimes_{\lt_2}\Z)Q_0
\]
such that $\psi(s_{(\mu,n)}) = \iota_1(s_{(\mu,(n,0))})\iota_2(d_2(\mu))$ for all $\mu \in \Lambda$ and $n \in \Z$. For each $n \in \Z$, this isomorphism $\psi$ restricts to an isomorphism $\psi_n : B^\Lambda_n \to Q_0(C^\Lambda_n\rtimes_{\lt_2}\Z)Q_0$, and induces an isomorphism $(\psi_n)_* : K_*(B^\Lambda_n) \to K_*(C^\Lambda_n\rtimes_{\lt_2}\Z)$.
\end{corollary}

\begin{proof} Everything but the statements about $\psi_n$ follows immediately from Lemma~\ref{lem:fullcorners}. Recall that the standard spanning family for $B^\Lambda_n$ is $\{s_{(\mu, n - d_1(\mu))} s^*_{(\nu, n - d_1(\nu))} : \mu,\nu \in \Lb\}$. Since each
\begin{align*}
\psi(s_{(\mu, n - d_1(\mu))} s^*_{(\nu, n - d_1(\nu))})
    = \iota_1(s_{(\mu, (n - d_1(\mu),0))})\iota_2(d_2(\mu) - d_2(\nu))\iota_1( & s^*_{(\nu, (n - d_1(\nu), 0))}) \\ & \in Q_0(C^\Lambda_n\rtimes_{\lt_2}\Z)Q_0,
\end{align*}
we see that $\psi$ restricts to a map $\psi_n : B^\Lambda_n \to Q_0(C^\Lambda_n\rtimes_{\lt_2}\Z)Q_0$. To see that it is surjective, observe that $C^\Lambda_n\rtimes_{\lt_2}\Z$ is spanned by the elements $S_{\mu, \nu, p} \coloneqq  \iota_1(s_{(\mu,p - d(\mu))} s_{(\nu,p - d(\nu))}^*)\iota_2(d_2(\mu) - d_2(\nu))$ where $p$ ranges over $\{n\} \times \Z \subseteq \Z^2$. For a given $\mu,\nu, p$ we have
\begin{align*}
Q_0 \iota_1(S_{\mu, \nu, p}) Q_0
    &= \begin{cases}
        \iota_1(S_{\mu, \nu, p}) &\text{ if $r((\mu,p-d(\mu))),r((\nu,p-d(\nu))) \in \Lb^0 \times \Z \times \{0\}$}\\
        0 &\text{ otherwise}\\
    \end{cases}\\
    &= \begin{cases}
        \iota_1(S_{\mu, \nu, p}) &\text{ if $p - d(\mu), p - d(\nu) \in \Z \times \{0\}$}\\
        0 &\text{ otherwise.}\\
    \end{cases}
\end{align*}
So $Q_0(C^\Lambda_n\rtimes_{\lt_2}\Z)Q_0$ is spanned by the elements $s_{(\mu,(n-d_1(\mu), 0))} s^*_{(\nu,(n-d_1(\nu),0))})\iota_2(d_2(\mu) - d_2(\nu))$, which are all in the image of $\psi_n$.

To see that $(C^\Lambda_n\rtimes_{\lt_2}\Z)Q_0(C^\Lambda_n\rtimes_{\lt_2}\Z) = C^\Lambda_n\rtimes_{\lt_2}\Z$, note that for each $(v, (n,m)) \in \Lb^0 \times \{n\} \times \{\Z\}$, we have $s_{(v,(n,m))} = \lt_2(-m)(s_{(v, (n, 0))})$, and so by covariance, $\iota_1(s_{(v,(n,m))})$ belongs to $i_2(\Z) Q_0 i_2(\Z) \cap (C^\Lambda_n\rtimes_{\lt_2}\Z) \subseteq (C^\Lambda_n\rtimes_{\lt_2}\Z)Q_0(C^\Lambda_n\rtimes_{\lt_2}\Z)$. So the ideal $(C^\Lambda_n\rtimes_{\lt_2}\Z)Q_0(C^\Lambda_n\rtimes_{\lt_2}\Z)$ contains an approximate identity for, and hence is all of, $C^\Lb_n\rtimes_{\lt_2}\Z$.
\end{proof}

Note that the inclusion $\iota_n : C^\Lb_n \to C^\Lb_{n+1}$ is equivariant under $\lt_2$, so it induces a homomorphism $\iota_n \times 1 : C^\Lb_n \rtimes_{\lt_2}\Z \to C^\Lb_{n+1}\rtimes_{\lt_2}\Z$.

\begin{lemma}\label{lem:Bn Cncrossedproduct inc}
Let $\Lambda$ be a row-finite $2$-graph with no sources. With $\iota_n$ and $\iota_n \times 1$ as above, the diagrams
	\[
	\begin{tikzcd}
		B^\Lb_n \arrow[r,"\psi_n"] \arrow[d,"\iota_n"] & C^\Lb_n\rtimes_{\lt_2}\Z \arrow[d,"\iota_n \times 1"]\\
		B^\Lb_{n+1} \arrow[r,"\psi_{n+1}"] & C^\Lb_{n+1}\rtimes_{\lt_2}\Z
	\end{tikzcd}
	\]
commute.
\end{lemma}

\begin{proof}
Let $S \coloneqq (s(\mu),n)(\Lb\times_{d_1}\Z)^{e_1} $ and $S' \coloneqq  (s(\mu),(n,0))(\Lb\times_d\Z^2)^{e_1}$. Then
	
	\begin{align*}
		\psi_{n+1}\left(\iota_n(s_{(\mu,n-d_1(\mu))}s_{(\nu,n-d_1(\nu))}^*)\right) & = \psi_{n+1}\left(\sum_{(\lb,n)\in S} s_{(\mu\lb,n-d_1(\mu))}s_{(\nu\lb,n-d_1(\nu))}^*\right) \\
		& = \sum_{(\lb,n,0)\in S'} \iota_1\left(s_{(\mu\lb,(n-d_1(\mu),0))})\iota_2(d_2(\mu)-d_2(\nu))\iota_1(s_{(\nu\lb,(n-d_1(\nu),0))}^*\right) \\
		& = (\iota_n\times 1) \left(\iota_1(s_{(\mu,(n-d_1(\mu),0))})\iota_2(d_2(\mu)-d_2(\nu))\iota_1(s_{(\nu,(n-d_1(\nu),0))}^*)\right) \\
		& = (\iota_n \times 1)\left(\psi_n (s_{(\mu,n-d_1(\mu))}s_{(\nu,n-d_1(\nu))}^*)\right).\qedhere
	\end{align*}
\end{proof}

\begin{lemma}\label{lem:Cnm K theory}
Let $\Lambda$ be a row-finite $2$-graph with no sources. For $n,m \in \Z$, there is an isomorphism $C^\Lb_{n,m} \cong \bigoplus_{v\in\Lb^0} \mathcal{K}_{\Lb v}$ that carries $s_{(\mu, (n,m)-d(\mu))} s^*_{(\nu, (n,m)-d(\nu))}$ to the matrix unit $\Theta_{\mu,\nu}$ in $\mathcal{K}_{\Lb s(\mu)}$. So there is an isomorphism $K_0(C^\Lb_{n,m}) \cong \Z\Lb^0$ that carries $[s_{(\mu, (n,m)-d(\mu))} s^*_{(\mu, (n,m)-d(\mu))}]$ to $\delta_{s(\mu)}$ for all $\mu\in \Lb$, and we have $K_1(C^\Lb_{n,m}) =0$. Writing $(\imath_{n,m})_*$ for the map in $K$-theory induced by the inclusion $C^\Lb_{n,m} \subseteq C^\Lb_{n+1,m}$, and $(\jmath_{n,m})_*$ for the one induced by $C^\Lb_{n,m} \subseteq C^\Lb_{n,m+1}$, the diagrams
	\[
	\begin{tikzcd}
		K_0(C^\Lb_{n,m})  \arrow[r,"(\imath_{n,m})_*"] \arrow[d,"\cong"] & K_0(C^\Lb_{n+1,m}) \arrow[d,"\cong"] & K_0(C^\Lb_{n,m}) \arrow[r,"(\jmath_{n,m})_*"] \arrow[d,"\cong"] & K_0(C^\Lb_{n,m+1}) \arrow[d,"\cong"] \\
		\Z \Lb^0 \arrow[r,"A_1^t"] & \Z\Lb^0, & \Z\Lb^0 \arrow[r,"A_2^t"] & \Z\Lb^0
	\end{tikzcd}
	\] commute.
	For all $n\in\Z$, we have $K_1(C^\Lambda_n) = 0$ and there is an isomorphism
	\[
    \varinjlim (\Z\Lb^0,A_2^t) \cong K_0(C^\Lb_n)
    \]
    that carries $[\delta_v]_m$ to $[s_{(v,n,m)}]$.
\end{lemma}

\begin{proof}
For a given $v \in \Lb^0$, the set $\{s_{(\mu, (n,m)-d(\mu))} s^*_{(\nu, (n,m)-d(\nu))} : s(\mu) = s(\nu)\}$ is a family of matrix units over $\Lb v$ and these families are mutually orthogonal for distinct $v, w$. The isomorphism in $K$-theory then follows because $K_0(\Kk) \cong \Z$ via $[\theta_{i,i}] \mapsto 1$ and $K_1(\Kk) = 0$, and because $K$-theory respects direct sums.

Temporarily writing $h_{n,m} : K_0(C^\Lb_{n,m}) \to \Z\Lb^0$ for the isomorphisms established in the preceding paragraph, the formula~\eqref{eq:CK in skew prod} shows that
\begin{align*}
h_{n+1, m}&\big(\imath_{n,m}([s_{(\mu, (n,m)-d(\mu))} s^*_{(\mu, (n,m)-d(\mu))}])\big)\\
    &= h_{n+1,m}\Big(\Big[\sum_{\lambda \in s(\mu)\Lambda^{e_1}} s_{(\mu\lambda, (n,m)-d(\mu))} s^*_{(\mu\lambda, (n,m)-d(\mu))}\Big]\Big)
    = \sum_{\lambda \in s(\mu)\Lambda^{e_1}} \delta_{s(\lambda)}
    = A_1^t \delta_{s(\mu)},
\end{align*}
so the left-hand diagram commutes, and the argument for the right-hand diagram is identical.

The final statement follows from continuity of $K$-theory.
\end{proof}

Looking at the maps induced by the inclusion $C^\Lb_{n,m} \subseteq C^\Lb_{n+1,m}$ for $m\in \Z$ of the direct limit $C^\Lb_n = \varinjlim(C^\Lb_{n,m},\iota_m)$ gives us the following corollary.

\begin{corollary}\label{cor:Cn k theory inc}
Let $\Lambda$ be a row-finite $2$-graph with no sources. Let $\mathcal{A}_1 : \varinjlim (\Z\Lb^0,A_2^t) \to \varinjlim (\Z\Lb^0,A_2^t)$ be the map induced by $A_1^t : \Z\Lb^0 \to \Z\Lb^0$ at each level of the direct limit. Then the following diagram commutes:
	\[
	\begin{tikzcd}
		\varinjlim (\Z\Lb^0,A_2^t)  \arrow[r,"\cong"] \arrow[d,"\mathcal{A}_1"] & K_0(C^\Lb_n) \arrow[d, "(\iota_n)*"] \\
		\varinjlim (\Z\Lb^0,A_2^t)  \arrow[r,"\cong"] & K_0(C^\Lb_{n+1}).
	\end{tikzcd}
	\]
\end{corollary}

\begin{lemma}\label{lem:Cnm lt2}
Let $\Lambda$ be a row-finite $2$-graph with no sources, and fix $m,n \in \Z$. The diagram 	
\[
	\begin{tikzcd}
		\Z\Lb^0 \arrow[r,"I"] \arrow[d,"\cong"] & \Z\Lb^0 \arrow[r,"A_2^t"] \arrow[d,"\cong"] & \Z\Lb^0 \arrow[d,"\cong"] \\
		K_0(C^\Lb_{n,m+1}) \arrow[r,"(\lt_2)_*"] & K_0(C^\Lb_{n,m}) \arrow[r,"(\iota_m)*"] & K_0(C^\Lb_{n,m+1}).
	\end{tikzcd}
	\]
    commutes. In particular, if $\mathcal{A}_2 : \varinjlim (\Z\Lb^0,A_2^t) \to \varinjlim (\Z\Lb^0,A_2^t)$ denotes the map induced by $A_2^t : \Z\Lb^0 \to \Z\Lb^0$ at each level of the direct limit, then the diagram
	\[
	\begin{tikzcd}
		\varinjlim (\Z\Lb^0,A_2^t) \arrow[r,"\mathcal{A}_2"] \arrow[d,"\cong"] & \varinjlim (\Z\Lb^0,A_2^t) \arrow[d,"\cong"] \\
		K_0(C^\Lb_{n}) \arrow[r,"(\lt_2)_*"]  & K_0(C^\Lb_{n})
	\end{tikzcd}
	\]
commutes.
\end{lemma}

\begin{proof}
	The proof is very similar to one of the steps in our computation of $K$-theory for $1$-graph $C^*$-algebras. The action $\lt_2$ of $\Z$ on the graph $C^*$-algebra $C^*(\Lb\times_d \Z^2)$ is given by $\lt_2 \left(s_{(\mu,n-d_1(\mu),m+1-d_2(\mu))}\right) = s_{(\mu,n-d_1(\mu),m-d_2(\mu))}$, so in $K$-theory, $(\lt_2)_*: [s_{(v,n,m+1)}] \mapsto [s_{(v,n,m)}]$.
	
	The second diagram commutes because the map $\mathcal{A}_2$ is induced by the above action at each level $m\in \Z$ of the direct limit.
\end{proof}

We have used $\mathcal{A}_2$ to denote the induced map on the direct limit $\varinjlim(\Z\Lb^0,A_2^t)$ in order to distinguish it from the induced map $A_2^{t,\infty}$ on $\varinjlim(\Z\Lb^0,A_1^t)$ as the bonding maps are not the same.

\begin{proof}[Proof of Proposition~\ref{prop: Bn step 1}]
	Naturality of Pimnser--Voiculescu sequences gives us the following commutative diagram in which the rows are exact and the vertical arrows are induced from the natural inclusions:
	\[
	\begin{tikzcd}
		0 \arrow[r] & K_1(C^\Lb_{n}\rtimes_{\lt_2}\Z) \arrow[r] \arrow[d,"(\iota_n \times 1)_*"] & K_0(C^\Lb_n) \arrow[r,"1-(\lt_2)_*"] \arrow[d,"(\iota_n)_*"] & K_0(C^\Lb_n)  \arrow[r,"(\iota_1)_*"] \arrow[d,"(\iota_n)_*"] & K_0(C^\Lb_{n}\rtimes_{\lt_2}\Z) \arrow[r] \arrow[d,"(\iota_n \times 1)_*"] & 0 \\
		0 \arrow[r] & K_1(C^\Lb_{n+1}\rtimes_{\lt_2}\Z) \arrow[r] & K_0(C^\Lb_{n+1}) \arrow[r,"1-(\lt_2)_*"] & K_0(C^\Lb_{n+1})  \arrow[r,"(\iota_1)_*"] & K_0(C^\Lb_{n+1}\rtimes_{\lt_2}\Z) \arrow[r] & 0.
	\end{tikzcd}
	\]
	
We claim that there are unique homomorphisms $\zeta_{n},\zeta_{n+1}, \eta_n, \eta_{n+1}$  such that the tube-like diagram
	\[
	\begin{tikzcd}[nodes={inner sep=3pt}]
		0 \arrow[r] &[-1em] K_1(B^\Lb_n) \arrow[r,"\zeta_n"] \arrow[d,"(\psi_n)_*"] \arrow[ddd,"(\iota_n)_*",dotted,bend right=49] & \varinjlim (\Z\Lb^0,A_2^t) \arrow[r,"1-\mathcal{A}_2"] \arrow[d,"\cong"] \arrow[ddd,"\mathcal{A}_1",dotted,bend right=49] & \varinjlim (\Z\Lb^0,A_2^t)  \arrow[r,"\eta_n"] \arrow[d,"\cong"] \arrow[ddd,dotted,"\mathcal{A}_1",bend right=49] & K_0(B^\Lb_n) \arrow[r] \arrow[d,"(\psi_n)_*"] \arrow[ddd,dotted,"(\iota_n)_*",bend right=49] &[-1em] 0 \\
		0 \arrow[r] & K_1(C^\Lb_{n}\rtimes_{\lt_2}\Z) \arrow[r] \arrow[d,"(\iota_n\times 1)_*"] & K_0(C^\Lb_n) \arrow[r,"1-(\lt_2)_*"] \arrow[d,"(\iota_n)_*"] & K_0(C^\Lb_n)  \arrow[r,"(\iota_1)_*"] \arrow[d,"(\iota_n)_*"] & K_0(C^\Lb_{n}\rtimes_{\lt_2}\Z) \arrow[r] \arrow[d,"(\iota_n \times 1)_*"] & 0 \\
		0 \arrow[r] & K_1(C^\Lb_{n+1}\rtimes_{\lt_2}\Z) \arrow[r] & K_0(C^\Lb_{n+1}) \arrow[r,"1- (\lt_2)_*"] & K_0(C^\Lb_{n+1})  \arrow[r,"(\iota_1)_*"] & K_0(C^\Lb_{n+1}\rtimes_{\lt_2}\Z) \arrow[r] & 0\\
		0 \arrow[r] & K_1(B^\Lb_{n+1}) \arrow[r,"\zeta_{n+1}"] \arrow[u,"(\psi_{n+1})_*"] & \varinjlim (\Z\Lb^0,A_2^t) \arrow[r,"1-\mathcal{A}_2"] \arrow[u,"\cong"] &\varinjlim (\Z\Lb^0,A_2^t)  \arrow[r,"\eta_{n+1}"] \arrow[u,"\cong"] & K_0(B^\Lb_{n+1}) \arrow[r] \arrow[u,"(\psi_{n+1})_*"] & 0
	\end{tikzcd}
	\]
commutes, and that for those homomorphisms the rows of the diagram are exact.
	
	To see this, we examine each part of the diagram in turn. The two squares formed by the first and fourth columns commute because they are induced by the commuting diagram of $C^*$-homomorphisms from Proposition~\ref{lem:Bn Cncrossedproduct inc}.
	
	The two squares formed by the second and third columns commute by Corollary~\ref{cor:Cn k theory inc}.
	
	There are unique maps $\zeta_{n},\zeta_{n+1}, \eta_n, \eta_{n+1}$ such that the corner squares commute because the first and third rows of vertical arrows are all isomorphisms by Corollary~\ref{cor:Bn k theory iso Cncrossedproduct} and Lemma~\ref{lem:Cnm K theory}.
	
	The subdiagram formed by the second and third sequence-rows commutes by naturality of the Pimsner--Voiculescu sequences as mentioned before.
	
	Finally, the top and bottom squares of the middle column formed by the second and third sequence-columns commute by Lemma~\ref{lem:Cnm lt2}. Thus the whole darned thing commutes and the rows are exact because the second and fourth rows of vertical arrows are isomorphisms.
	
We now obtain the diagram in the statement of Proposition~\ref{prop: Bn step 1} from Lemma~\ref{lem:PR naturality} for
\begin{gather*}
H^1 = K_1(B^\Lb_n),\quad H^2 = K_1(B^\Lb_{n+1}),\quad G^1 = G^2 = \Z\Lb^0,\quad K^1 = K_0(B^\Lb_n),\quad K^2 = K_0(B^\Lb_{n+1}),\\
\phi^1 = \phi^2 = A_2^t, \quad \psi = A_1^t, \quad \eta = (\iota_n)_*,\quad\text{ and }\quad \theta = (\iota_n)_*,
\end{gather*} except that we have $\tau^1, \tau^2$ instead of $\phi_{n,0}$ and $\phi_{n+1,0}$ respectively. To finish, we need to show that $\tau^1 = \phi_{n,0} : \delta_v \mapsto [p_{(v,n)}]$ (and similarly for $\tau^2 = \phi_{n+1,0}$).

To that end, we claim that $\eta_n([\delta_v]_m) = [p_{(v,n)}]$ for all $m\ge 0$. Following the left and bottom arrows of that square, $[\delta_v]_m$ gets mapped first to $[s_{(v,n,m)}] \in K_0(C^\Lb_n)$ and then to $[\iota_1(s_{(v,n,m)})] \in K_0(C^\Lb_n\rtimes_{\lt_2}\Z)$. Now $\iota_1(s_{(v,n,0)}) = \iota_2(m)\iota_1(s_{(v,n,m)})\iota_2(-m)$ since the action is left-translation in the second coordinate, so that  $\iota_1(s_{(v,n,0)})$ is unitarily equivalent to $\iota_1(s_{(v,n,m)})$. In particular, $(\psi_n)_*([p_{(v,n)}]) = [\iota_1(s_{(v,n,m)})] = [\iota_1(s_{(v,n,0)})]$. Thus $(\psi_n)_*([p_{(v,n)}]) = (\psi_n)_*(\eta_n([\delta_v]_n)$ so that by injectivity of $(\psi_n)_*$, $\eta_n$ is as claimed. Now $\eta_n$ corresponds to $\pi^1$ in Lemma \ref{lem:PR naturality}, so if $\tilde{\eta}_n$ is the map induced by $\eta_n$ on $\coker(1-\mathcal{A}_2)$, $\tilde{\phi}_{1,\infty}$ is the isomorphism from Lemma \ref{lem:paskraeburn}, and $q:\Z\Lb^0 \to \coker(1-A_2^t)$ is the quotient map, then $\tau^1 = \tilde{\eta}_n\circ \tilde{\phi}_{1,\infty} \circ q$. Thus
\[
    \tau^1(\delta_v) = \tilde{\eta}_n\circ \tilde{\phi}_{1,\infty}(\delta_v + \image(1-A_2^t)) = \tilde{\eta}_n([\delta_v]_1+\image(1-\mathcal{A}_2)) = \eta_n([\delta_v]_1) = [p_{(v,n)}].
\]
Similarly, $\tau^2(\delta_v) = [p_{(v,n+1)}]$, and we are done.
\end{proof}

\subsection{Step 2}\label{subsec:step2}

The following is a $2$-graph analogue of Theorem~\ref{thm:1-graph K-th} for the $B^\Lb_n$.

\begin{proposition}\label{prop: hereditary inclusions step 2}
Let $\Lb$ be a row-finite $2$-graph with no sources, and let $H \subseteq \Lb^0$ be a hereditary subset. Let $\phi_{n,i}$ and $\phi^{H}_{n,i}$ be the maps of Proposition~\ref{prop: Bn step 1} applied to the $2$-graphs $\Lambda$ and $H\Lambda$ respectively. Then the diagram
	\[
	\begin{tikzcd}
		0 \arrow[r] & K_1(B_n^{H\Lb}) \arrow[r, "{\phi_{n,1}^H}"] \arrow[d, "(\iota_{H\times \Z})_*"] & \Z H \arrow[r, "{1-A_{H,2}^t}"] \arrow[d, hook] & \Z H \arrow[r, "{\phi_{n,0}^H}"] \arrow[d, hook] & K_0(B_n^{H\Lb}) \arrow[r] \arrow[d, "(\iota_{H\times \Z})_*"] & 0 \\
		0 \arrow[r] & K_1(B^\Lb_{n}) \arrow[r, "\phi_{n,1}"]                                         & \Z\Lambda^0 \arrow[r, "1-A_2^t"]                & \Z\Lambda^0 \arrow[r, "\phi_{n,0}"]                & K_0(B^\Lb_{n}) \arrow[r]                                 & 0.
	\end{tikzcd}\]
commutes.
\end{proposition}
\begin{proof}
In this proof, we let $\iota \coloneqq  \iota_{H\times\Z^2}$ denote the inclusion of $C_n^{H\Lb}$ into $C^\Lb_n$ analogous to $\iota_{H\times \Z} : B_n^{H\Lb} \to B^\Lb_n$, but we will drop the extra subscripts to avoid notational clutter. Likewise, $\iota \coloneqq  \iota_{H\times\Z^2}\times 1 :C_n^{H\Lb}\rtimes_{\lt_2}\Z \to C^\Lb_n\rtimes_{\lt_2}\Z $, denuded of its adornments, will also denote the homomorphism induced by the equivariant inclusion of $C_n^{H\Lb} \hookrightarrow C^\Lb_n$.
	
Naturality of Pimnser--Voiculescu sequences yields the following commutative diagram in which the rows are exact and the vertical arrows are induced from the natural inclusions:
	\[
	\begin{tikzcd}
		0 \arrow[r] & K_1(C_{n}^{H\Lb}\rtimes_{\lt_2}\Z) \arrow[r] \arrow[d,"\iota_*"] & K_0(C_n^{H\Lb}) \arrow[r,"1-(\lt_2)_*"] \arrow[d,"\iota_*"] & K_0(C_n^{H\Lb})  \arrow[r,"(\iota_1)_*"] \arrow[d,"\iota_*"] & K_0(C_{n}^{H\Lb}\rtimes_{\lt_2}\Z) \arrow[r] \arrow[d,"\iota_*"] & 0 \\
		0 \arrow[r] & K_1(C^\Lb_{n}\rtimes_{\lt_2}\Z) \arrow[r] & K_0(C^\Lb_{n}) \arrow[r,"1-(\lt_2)_*"] & K_0(C^\Lb_{n})  \arrow[r,"(\iota_1)_*"] & K_0(C^\Lb_{n}\rtimes_{\lt_2}\Z) \arrow[r] & 0.
	\end{tikzcd}
	\]
	
Just as in the proof of Proposition \ref{prop: Bn step 1}, we claim that there are unique homomorphisms $\zeta_{n},\zeta_{n}^H, \eta_n, \eta_{n}^H$ that make the tube-like diagram
	\[
	\begin{tikzcd}[nodes={inner sep=3pt}]
		0 \arrow[r] &[-1em] K_1(B_n^{H\Lb}) \arrow[r,"\zeta_n^H"] \arrow[d,"(\psi_n^H)_*"] \arrow[ddd,"(\iota_{H\times\Z})_*",dotted,bend right = 55] & \varinjlim (\Z H,A_{H,2}^t) \arrow[r,"1-\mathcal{A}_{H,2}"] \arrow[d,"\cong"] \arrow[ddd,"\iota",dotted,bend right = 55] & \varinjlim (\Z H,A_{H,2}^t)  \arrow[r,"\eta_n^H"] \arrow[d,"\cong"] \arrow[ddd,"\iota",dotted, bend right = 55] & K_0(B_n^{H\Lb}) \arrow[r] \arrow[d,"(\psi_n^H)_*"] \arrow[ddd,"(\iota_{H\times\Z})_*",dotted, bend right=55] &[-1em] 0 \\
		0 \arrow[r] &[-1em] K_1(C_{n}^{H\Lb}\rtimes_{\lt_2}\Z) \arrow[r] \arrow[d,"\iota_*"] & K_0(C_n^{H\Lb}) \arrow[r,"1-(\lt_2)_*"] \arrow[d,"\iota_*"] & K_0(C_n^{H\Lb})  \arrow[r,"(\iota_1)_*"] \arrow[d,"\iota_*"] & K_0(C_{n}^{H\Lb}\rtimes_{\lt_2}\Z) \arrow[r] \arrow[d,"\iota_*"] &[-1em] 0 \\
		0 \arrow[r] &[-1em] K_1(C^\Lb_{n}\rtimes_{\lt_2}\Z) \arrow[r] & K_0(C^\Lb_{n}) \arrow[r,"1- (\lt_2)_*"] & K_0(C^\Lb_{n})  \arrow[r,"(\iota_1)_*"] & K_0(C^\Lb_{n}\rtimes_{\lt_2}\Z) \arrow[r] &[-1em] 0\\
		0 \arrow[r] &[-1em] K_1(B^\Lb_{n}) \arrow[r,"\zeta_{n}"] \arrow[u,"(\psi_{n})_*"] & \varinjlim (\Z\Lb^0,A_2^t) \arrow[r,"1-\mathcal{A}_2"] \arrow[u,"\cong"] &\varinjlim (\Z\Lb^0,A_2^t)  \arrow[r,"\eta_{n}"] \arrow[u,"\cong"] & K_0(B^\Lb_{n}) \arrow[r] \arrow[u,"(\psi_{n})_*"] &[-1em] 0
	\end{tikzcd}
	\]
commute, and that with these homomorphisms, the rows of the diagram are exact.
	
The two squares formed by the first and fourth columns commute because they are induced from the diagram \[
	\begin{tikzcd}
		B_n^{H\Lb} \arrow[r,"\psi_n^H"] \arrow[d,"\iota_{H\times \Z}"] & C_n^{H\Lb}\rtimes_{\lt_2}\Z \arrow[d,"\iota"] \\
		B^\Lb_n \arrow[r,"\psi_n"] & C^\Lb_n \rtimes_{\lt_2}\Z
	\end{tikzcd}
	\] which commutes because $\psi_n(s_{(\mu,n)}) = \iota_1(s_{(\mu,n,0)})\iota_2(d_2(\mu))$ and likewise for $\psi_n^H$ (see Corollary~\ref{cor:Bn k theory iso Cncrossedproduct}).
	
	The second and third columns both just give the diagram
	\[
	\begin{tikzcd}
		\varinjlim (\Z H,A_{H,2}^t) \arrow[r,"\cong"] \arrow[d,"\iota"] & K_0(C_n^{H\Lb}) \arrow[d,"\iota_*"] \\
		\varinjlim (\Z \Lb^0,A_{2}^t) \arrow[r,"\cong"] & K_0(C^\Lb_n),
	\end{tikzcd}
	\] which commutes by Lemma~\ref{lem:Cnm K theory}.
	
	The subdiagram formed by the second and third rows commutes by naturality of the Pimsner--Voiculescu sequence.
	
	For the top and bottom bands of squares, see the top band of the big diagram in the proof of Proposition~\ref{prop: Bn step 1} applied to the $2$-graphs $H\Lambda$ and $\Lb$.
	
	Thus the whole diagram commutes. The rows are exact because the second and fourth rows of vertical arrows are isomorphisms.
	
	Now the result follows from Lemma~\ref{lem:PR naturality} for
\begin{gather*}
H^1 = K_1(B_n^{H\Lb}),\quad H^2 = K_1(B^\Lb_n),\quad G^1 = \Z H,\quad G^2 = \Z\Lb^0, \quad K^1 = K_0(B_n^{H\Lb}),\\
K^2 = K_0(B^\Lb_n), \quad \phi^1 = A_{H,2}^t,\quad \phi^2 = A_2^t,\quad  \psi = \iota,\quad \eta = (\iota_{H\times\Z})_*,\quad\text{ and }\quad \theta = (\iota_{H\times\Z})_*. \qedhere
\end{gather*}
\end{proof}

\subsection{Step 3 and proof of Theorem\texorpdfstring{~\ref{thm:main k theory theorem}}{ 4.3}}\label{subsec:step3}

The last main piece of the puzzle we shall later assemble is what's needed for an analogue for $B^\Lb_n$ analogue of equation~\eqref{eq:1-graph translation}.

\begin{proposition}\label{prop:2-graph skewed lt}
Let $\Lambda$ be a row-finite $2$-graph. Let $\phi_{n,i}$ be the maps of Proposition~\ref{prop: Bn step 1}. Then the diagram
	\[
	\begin{tikzcd}
		0 \arrow[r] & K_1(B^\Lb_n) \arrow[r, "{\phi_{n,1}}"] \arrow[d, "(\lt_1)_*"] & \Z\Lambda^0 \arrow[r, "1-A_2^t"] \arrow[d, "A_1^t"] & \Z\Lambda^0 \arrow[r, "{\phi_{n,0}}"] \arrow[d, "A_1^t"] & K_0(B^\Lb_n) \arrow[r] \arrow[d, "(\lt_1)_*"] & 0 \\
		0 \arrow[r] & K_1(B^\Lb_{n}) \arrow[r, "{\phi_{n,1}}"]                     & \Z\Lambda^0 \arrow[r, "1-A_2^t"]                    & \Z\Lambda^0 \arrow[r, "{\phi_{n,0}}"]                    & K_0(B^\Lb_{n}) \arrow[r]                     & 0
	\end{tikzcd}\]
commutes.
\end{proposition}
\begin{proof}
Since the homomorphism $\lt_1:C^\Lb_n \to C^\Lb_n$ commutes with $\lt_2$, naturality of the Pimnser--Voiculescu sequence implies that the diagram
	\[
	\begin{tikzcd}
		0 \arrow[r] & K_1(C^\Lb_{n}\rtimes_{\lt_2}\Z) \arrow[r] \arrow[d,"(\lt_1\times 1)_*"] & K_0(C^\Lb_n) \arrow[r,"1-(\lt_2)_*"] \arrow[d,"(\lt_1)_*"] & K_0(C^\Lb_n)  \arrow[r,"(\iota_1)_*"] \arrow[d,"(\lt_1)_*"] & K_0(C^\Lb_{n}\rtimes_{\lt_2}\Z) \arrow[r] \arrow[d,"(\lt_1 \times 1)_*"] & 0 \\
		0 \arrow[r] & K_1(C^\Lb_{n}\rtimes_{\lt_2}\Z) \arrow[r] & K_0(C^\Lb_{n}) \arrow[r,"1-(\lt_2)_*"] & K_0(C^\Lb_{n})  \arrow[r,"(\iota_1)_*"] & K_0(C^\Lb_{n}\rtimes_{\lt_2}\Z) \arrow[r] & 0
	\end{tikzcd}
	\]
commutes and has exact rows.
	
	Just as in the proofs of Propositions \ref{prop: Bn step 1}~and~\ref{prop: hereditary inclusions step 2}, there are unique homomorphisms $\zeta_n$ and $\eta_n$ such that the corner squares of the 24-term tube-like diagram
	\[
	\begin{tikzcd}
		0 \arrow[r] &[-1em] K_1(B^\Lb_n) \arrow[r,"\zeta_n"] \arrow[d,"(\psi_n)_*"] \arrow[ddd,"(\lt_1)_*",dotted,bend right=49] & \varinjlim (\Z\Lb^0,A_2^t) \arrow[r,"1-\mathcal{A}_2"] \arrow[d,"\cong"] \arrow[ddd,"\mathcal{A}_1",dotted,bend right=49] & \varinjlim (\Z\Lb^0,A_2^t)  \arrow[r,"\eta_n"] \arrow[d,"\cong"] \arrow[ddd,"\mathcal{A}_1",dotted,bend right=49] & K_0(B^\Lb_n) \arrow[r] \arrow[d,"(\psi_n)_*"] \arrow[ddd,"(\lt_1)_*",dotted,bend right=49] &[-1em] 0 \\
		0 \arrow[r] & K_1(C^\Lb_{n}\rtimes_{\lt_2}\Z) \arrow[r] \arrow[d,"(\lt_1 \times 1)_*"] & K_0(C^\Lb_n) \arrow[r,"1-(\lt_2)_*"] \arrow[d,"(\lt_1)_*"] & K_0(C^\Lb_n)  \arrow[r,"(\iota_1)_*"] \arrow[d,"(\lt_1)_*"] & K_0(C^\Lb_{n}\rtimes_{\lt_2}\Z) \arrow[r] \arrow[d,"(\lt_1 \times 1)_*"] & 0 \\
		0 \arrow[r] & K_1(C^\Lb_{n}\rtimes_{\lt_2}\Z) \arrow[r] & K_0(C^\Lb_{n}) \arrow[r,"1- (\lt_2)_*"] & K_0(C^\Lb_{n})  \arrow[r,"(\iota_1)_*"] & K_0(C^\Lb_{n}\rtimes_{\lt_2}\Z) \arrow[r] & 0\\
		0 \arrow[r] & K_1(B^\Lb_{n}) \arrow[r,"\zeta_{n}"] \arrow[u,"(\psi_{n})_*"] & \varinjlim (\Z\Lb^0,A_2^t) \arrow[r,"1-\mathcal{A}_2"] \arrow[u,"\cong"] &\varinjlim (\Z\Lb^0,A_2^t)  \arrow[r,"\eta_{n}"] \arrow[u,"\cong"] & K_0(B^\Lb_{n}) \arrow[r] \arrow[u,"(\psi_{n})_*"] & 0
	\end{tikzcd}
	\]
commute. We claim that the whole diagram commutes and that its rows are exact.
	
Consider the diagram
	\begin{equation}\label{eq:C*-diagram}
	\begin{tikzcd}
		B^\Lb_{n+1} \arrow[r,"\lt_1"] \arrow[d,"\psi_{n+1}"] & B^\Lb_n \arrow[r,"\iota_n"] \arrow[d,"\psi_n"] & B^\Lb_{n+1} \arrow[d,"\psi_{n+1}"] \\
		C^\Lb_{n+1}\rtimes_{\lt_2}\Z \arrow[r,"\lt_1 \times 1"] & C^\Lb_n \rtimes_{\lt_2}\Z \arrow[r,"\iota_n"] & C^\Lb_{n+1} \rtimes_{\lt_2}\Z
	\end{tikzcd}
	\end{equation}
of $C^*$-homomorphisms. The left-hand square commutes because $\lt_1$ commutes with $\lt_2$. The right-hand square commutes by Lemma~\ref{lem:Bn Cncrossedproduct inc}. So the whole diagram commutes. Since the two squares formed by the first and fourth columns of the 24-term diagram above are induced by~\eqref{eq:C*-diagram}, it follows that they commute.
	
	The diagram
	\[
	\begin{tikzcd}
		\Z\Lambda^0 \arrow[r,"I"] \arrow[d] & \Z\Lambda^0 	\arrow[r,"A_1^t"] \arrow[d] & \Z\Lambda^0 \arrow[d] \\
		K_0(C^\Lb_{n+1,m}) \arrow[r,"(\lt_1)_*"] & K_0(C^\Lb_{n,m}) 	\arrow[r,"(\iota_n)_*"] & K_0(C^\Lb_{n+1,m})
	\end{tikzcd}
	\]
is analogous to that in Lemma~\ref{lem:Cnm lt2} and the proof that it commutes is more or less exactly the same after remembering that the induced map is $\mathcal{A}_1$ (not the same as $A_1^{t,\infty}$ in our notation); we won't repeat the argument.

	Thus the squares formed by the second and third columns of the 24-term diagram commute by Corollary~\ref{cor:inductivelimitsubalgebras}.
	
	The subdiagram of the 24-term diagram formed by the second and third rows commutes by naturality of the Pimsner--Voiculescu sequences.
	
	Finally, the top and bottom bands of squares in the 24-term diagram are exactly the same as the top band in the 24-term commuting diagram in the proof of Proposition~\ref{prop: Bn step 1}. So the whole diagram commutes, and the rows are exact because the second and fourth rows of vertical arrows are isomorphisms.
	
	Now the result follows from Lemma~\ref{lem:PR naturality} applied with
    \begin{gather*}
    H^1 = H^2 = K_1(B^\Lb_n),\quad G^1 = G^2 = \Z\Lb^0, \quad K^1 = K^2 = K_0(B^\Lb_n),\\
    \phi^1 = \phi^2 = A_2^t,\quad  \psi = A_1^t, \quad \eta = (\lt_1)_*, \quad\text{ and }\quad \theta = (\lt_1)_*. \qedhere
    \end{gather*}
\end{proof}

We are finally able to prove Theorem~\ref{thm:main k theory theorem}.

\begin{proof}[Proof of Theorem~\ref{thm:main k theory theorem}]
	Proposition~\ref{prop: Bn step 1} implies that the diagrams
	\[
	\begin{tikzcd}
		K_1(B^\Lb_n) \arrow[d, "(\iota_n)_*"] \arrow[r, "{\phi_{n,1}}"] & \Z\Lb^0 \arrow[d, "A_1^t"] &  & \Z\Lb^0 \arrow[d, "A_1^t"] \arrow[r, "{\phi_{n,0}}"] & K_0(B^\Lb_n) \arrow[d, "(\iota_n)_*"] \\
		K_1(B^\Lb_{n+1}) \arrow[r, "{\phi_{n+1,1}}"]                    & {\Z\Lb^0,}                 &  & \Z\Lb^0 \arrow[r, "{\phi_{n+1,0}}"]                  & K_0(B^\Lb_{n+1})
	\end{tikzcd}\]
commute. Since direct limits of groups preserve exact sequences (Lemma~\ref{lem:direct limits}), we can take a direct limit of the entire commutative diagram in Proposition~\ref{prop: hereditary inclusions step 2} to obtain the commuting diagram
	\[
	\begin{tikzcd}[nodes={inner sep=3pt}]
		0 \arrow[r] &[-1.5em] {\varinjlim\big(K_1(B_n^{H\Lb}),(\iota_n)_*\big)} \arrow[d] \arrow[r] &[-1.25em] {\varinjlim(\Z H,A_{H,1}^t)} \arrow[d, hook] \arrow[r, "{1-A_{H,2}^{t,\infty}}"] & {\varinjlim(\Z H,A_{H,1}^t)} \arrow[d, hook] \arrow[r] &[-1.25em] {\varinjlim\big(K_0(B_n^{H\Lb}),(\iota_n)_*\big)} \arrow[d] \arrow[r] &[-1.5em] 0 \\
		0 \arrow[r] &[-1.5em] {\varinjlim\big(K_1(B^\Lb_n),(\iota_n)_*\big)} \arrow[r]             &[-0.25em] {\varinjlim(\Z\Lambda^0,A_1^t)} \arrow[r, "{1-A_2^{t,\infty}}"]              & {\varinjlim(\Z\Lambda^0,A_1^t)} \arrow[r]          &[-0.25em] {\varinjlim\big(K_0(B^\Lb_n),(\iota_n)_*\big)} \arrow[r]             &[-1.25em] 0
	\end{tikzcd}\]
of exact sequences. Theorem~\ref{thm:main k theory theorem}(1) then follows from an application of Corollary~\ref{cor:inductivelimitsubalgebras}(1).
	
Similarly, we can take a direct limit of the entire commutative diagram in Proposition~\ref{prop:2-graph skewed lt} to obtain the commuting diagram
	 \[
	 \begin{tikzcd}[nodes={inner sep=3pt}]
	 	0 \arrow[r] &[-1.25em] {\varinjlim\big(K_1(B^\Lb_n),(\iota_n)_*\big)} \arrow[r] \arrow[d, "1-(\lt_1)_*"] &[-0.25em] {\varinjlim(\Z\Lambda^0,A_1^t)} \arrow[r, "{1-A_2^{t,\infty}}"] \arrow[d, "{1-A_1^{t,\infty}}"] & {\varinjlim(\Z\Lambda^0,A_1^t)} \arrow[r] \arrow[d, "{1-A_2^{t,\infty}}"] &[-0.25em] {\varinjlim\big(K_0(B^\Lb_n),(\iota_n)_*\big)} \arrow[r] \arrow[d, "1-(\lt_1)_*"] &[-1.25em] 0 \\
	 	0 \arrow[r] &[-1.25em] {\varinjlim\big(K_1(B^\Lb_n),(\iota_n)_*\big)} \arrow[r]                                   & {\varinjlim(\Z\Lambda^0,A_1^t)} \arrow[r, "{1-A_1^{t,\infty}}"]                                 & {\varinjlim(\Z\Lambda^0,A_1^t)} \arrow[r]                                 & {\varinjlim\big(K_0(B^\Lb_n),(\iota_n)_*\big)} \arrow[r]                         &[-1.25em] 0
	 \end{tikzcd}\]
of exact sequences.	Theorem~\ref{thm:main k theory theorem}(2) then follows as before from an application of Corollary~\ref{cor:inductivelimitsubalgebras}(1).
\end{proof}

\subsection{\texorpdfstring{$K$}{K}-theory of \texorpdfstring{$C^*(H\Lb)$}{C*(H Lambda)} and \texorpdfstring{$C^*(\Lb)$}{C*(Lambda)}}

We will now go from our results about the $K$-theory of skew product $2$-graph algebras to the $K$-theory of the original $2$-graph algebras themselves. Let $\Lb$ be a row-finite $2$-graph with no sources and $H$ be a saturated hereditary subset of $\Lb^0$. Denote the inclusion of $C^*(H\Lb)$ into $C^*(\Lb)$ by $\iota_H$ and recall the map $\iota_{H\times \Z}$ from Lemma~\ref{lem:saturated product}. Since $\iota_{H\times\Z}$ commutes with $\lt_1$, it induces a homomorphism $\iota_{H\times\Z}\times 1: C^*(H\Lb\times_{d_1}\Z)\rtimes_{\lt_1}\Z \to C^*(\Lb\times_{d_1}\Z)\rtimes_{\lt_1}\Z$.


\begin{lemma}\label{lem: takai duality}
Let $\Lambda$ be a row-finite $2$-graph with no sources. Let $P_0 \in \mathcal{M}(C^*(\Lb\times_{d_1}\Z)\rtimes_{\lt_1}\Z)$ be the projection
	$$P_0 = \sum_{v\in \Lb^0}\iota_1(p_{(v,0)}).$$
Then there is an isomorphism $\rho: C^*(\Lb) \to P_0 (C^*(\Lb\times_{d_1}\Z)\rtimes_{\lt_1}\Z) P_0$ such that
\begin{equation}\label{eq:rho}
    \rho(s_{\mu}) = \iota_1(s_{(\mu,0)})\iota_2(d_1(\mu))\quad\text{ for all $\mu \in \Lb$.}
\end{equation}
The induced homomorphism $\rho_* : K_*(C^*(\Lb)) \to K_*(C^*(\Lb\times_{d_1}\Z)\rtimes_{\lt_1}\Z)$ is an isomorphism.
Let $H \subseteq \Lambda^0$ be a hereditary set, and write $\rho_H : C^*(H\Lb) \to C^*(H\Lb\times_{d_1}\Z)\rtimes_{\lt_1}\Z$ for the homomorphism obtained from the first statement for the $2$-graph $H\Lambda$. Then the diagram
	\[
	\begin{tikzcd}
		C^*(H\Lb) \arrow[r,"\rho_{H}"] \arrow[d,"\iota_H"] & C^*(H\Lb \times_{d_1}\Z)\rtimes_{\lt_1} \Z \arrow[d,"\iota_{H\times\Z}\times 1"] \\
		C^*(\Lb) \arrow[r,"\rho"] & C^*(\Lb \times_{d_1}\Z)\rtimes_{\lt_1} \Z
	\end{tikzcd}
	\]
commutes.
\end{lemma}

\begin{proof}
Statement~(1) is an application of Lemma~\ref{lem:fullcorners}---the induced homomorphism $\rho_*$ is an isomorphism because $P_0$ is full. Statement~(2) follows directly from the formula~\eqref{eq:rho}.
\end{proof}

\begin{theorem}\label{thm:inc}
Let $\Lambda$ be a row-finite $2$-graph with no sources, and let $H \subseteq \Lambda^0$ be a hereditary set. Let $\tilde{A}_1^t \in \operatorname{End}(\coker(1 - A_2^t))$ be the map induced by $A_1^t \in \operatorname{End}(\Z\Lb^0)$, and let $A^t_1| \in \operatorname{End}(\ker(1 - A_2^t))$ be the restriction of $A_1^t$. Write $\tilde{A}^t_{H, 1}$ and $A^t_{H,1}|$ for the corresponding maps for $H\Lambda$. There is a homomorphism $j : \coker(1-\tilde{A}_1^{t}) \to K_0(C^*(\Lb))$ such that \[j\big(\delta_v+\image(1-A_2^t) + \image(1-\tilde{A}_1^t)\big) = [p_v]\] for all $v \in \Lb^0$, and a corresponding homomorphism $j_H$ for $H\Lambda$. There exist homomorphisms $\tau_H, \tau$ such that the diagram
	\[
	\begin{tikzcd}
		0 \arrow[r] & \coker(1-\tilde{A}_{H,1}^{t}) \arrow[r,"j_H"] \arrow[d,"\tilde{\iota}"] & K_0(C^*(H\Lb)) \arrow[r,"\tau_H"] \arrow[d,"(\iota_H)_*"] & \ker(1-A_{H,1}^{t}|) \arrow[r] \arrow[d,"\iota"] & 0 \\
		0 \arrow[r] & \coker(1-\tilde{A}_{1}^{t}) \arrow[r,"j"] & K_0(C^*(\Lb)) \arrow[r,"\tau"] & \ker(1-A_{1}^{t}|) \arrow[r] & 0
	\end{tikzcd}
	\]
commutes and has exact rows. There exist homomorphisms from $\coker(1-A_{1}^{t}|)$ to $K_1(C^*(\Lb))$ and from $K_1(C^*(\Lb))$ to $\ker(1-\tilde{A}_{1}^{t})$, and similarly for $H\Lambda$, such that the diagram
	\[
	\begin{tikzcd}
		0 \arrow[r] & \coker(1-A_{H,1}^{t}|) \arrow[r] \arrow[d,"\tilde{\iota}"] & K_1(C^*(H\Lb)) \arrow[r] \arrow[d,"(\iota_H)*"] & \ker(1-\tilde{A}_{H,1}^{t}) \arrow[r] \arrow[d,"\iota"] & 0 \\
		0 \arrow[r] & \coker(1-A_{1}^{t}|) \arrow[r] & K_1(C^*(\Lb)) \arrow[r] & \ker(1-\tilde{A}_{1}^{t}) \arrow[r] & 0
	\end{tikzcd}
	\]
commutes and has exact rows.
\end{theorem}
	\begin{proof}
	In this proof, we will drop the adornments on the map $\iota_{H\times\Z}\times 1$ to avoid cluttering the diagrams and just write $\iota$. Naturality of Pimsner--Voiculescu sequences gives us the following commutative diagram in which the six term sequence forming the inner cycle and the six term sequence forming the outer cycle are both exact:
	\begin{equation}\label{eq:nested cyclic} 
	\begin{tikzcd}[nodes={inner sep=3pt}]
		K_0(C^*(H\Lambda\times_{d_1}\Z)) \arrow[r, "1-(\lt_1)_*"] \arrow[d, "(\iota_{H\times\Z})_*", dashed]                  & K_0(C^*(H\Lambda\times_{d_1}\Z)) \arrow[d, "(\iota_{H\times\Z})_*", dashed] \arrow[r, "(\iota_1)_*"] & K_0(C^*(H\Lambda\times_{d_1}\Z)\rtimes_{\lt_1}\Z) \arrow[ddd, bend left=75] \arrow[d, "\iota_*", dashed]                        \\
		K_0(C^*(\Lambda\times_{d_1}\Z)) \arrow[r, "1-(\lt_1)_*"]                                                & K_0(C^*(\Lambda\times_{d_1}\Z)) \arrow[r, "(\iota_1)_*"]                                           & K_0(C^*(\Lambda\times_{d_1}\Z)\rtimes_{\lt_1}\Z) \arrow[d,bend left =60]              \\
		K_1(C^*(\Lambda\times_{d_1}\Z)\rtimes_{\lt_1}\Z) \arrow[u,bend left=60]                                               & K_1(C^*(\Lambda\times_{d_1}\Z)) \arrow[l, "(\iota_1)_*"]                                           & K_0(C^*(\Lambda\times_{d_1}\Z)) \arrow[l, "1-(\lt_1)_*"]                                           \\
		K_1(C^*(H\Lambda\times_{d_1}\Z)\rtimes_{\lt_1}\Z) \arrow[uuu, bend left=75] \arrow[u, "\iota_*", dashed] & K_1(C^*(H\Lambda\times_{d_1}\Z)) \arrow[u, "(\iota_{H\times\Z})_*", dashed] \arrow[l, "(\iota_1)_*"]             & K_1(C^*(H\Lambda\times_{d_1}\Z)). \arrow[l, "1-(\lt_1)_*"] \arrow[u, "(\iota_{H\times\Z})_*", dashed]
	\end{tikzcd}\end{equation}
	
	We want to add an outer cycle and an inner cycle just as we have added rows to previous diagrams; to avoid a graphical bloodbath, we stretch the cycles in the previous diagram into columns so that the newly-added cyclic sequences appear as the leftmost and rightmost columns. So we claim that there are unique homomorphisms $\theta_H, \theta, \bar{\theta}_H, \bar{\theta}, \xi^H, \xi, \bar{\xi}_H,\bar{\xi}$ such that the diagram
\[
	\begin{tikzpicture}[baseline=(a)]
		\node[xscale=0.925] (a) at (0,0){
			\begin{tikzcd}[nodes={inner sep=0.5pt}, row sep=3em]
				\varinjlim(\coker(1-A_{H,2}^t),\tilde{A}_{H,1}^t) \arrow[r,swap,pos=0.4,"\cong"] \arrow[d,"1-\tilde{A}_{H,1}^{t,\infty}"] \arrow[rrr,"\iota",dashed,bend left = 8] &[-2.75em]  K_0(C^*(H\Lb\times_{d_1}\Z)) \arrow[r,"(\iota_{H\times\Z})_*"] \arrow[d,"1-(\lt_1)_*"] &[-1.25em] K_0(C^*(\Lb\times_{d_1}\Z)) \arrow[d,"1-(\lt_1)_*"] &[-2.75em] \varinjlim(\coker(1-A_{2}^t),\tilde{A}_{1}^t) \arrow[l,pos=0.4,"\cong"] \arrow[d,"1-\tilde{A}_{1}^{t,\infty}"] \\
				\varinjlim(\coker(1-A_{H,2}^t),\tilde{A}_{H,1}^t) \arrow[r,swap,pos=0.4,"\cong"] \arrow[d,"\theta_H"] \arrow[rrr,"\iota",dashed,bend left = 8] & K_0(C^*(H\Lb\times_{d_1}\Z)) \arrow[r,"(\iota_{H\times\Z})_*"] \arrow[d,"(\iota_1)_*"] & K_0(C^*(\Lb\times_{d_1}\Z)) \arrow[d,"(\iota_1)_*"] & \varinjlim(\coker(1-A_{2}^t),\tilde{A}_{1}^t) \arrow[l,pos=0.4,"\cong"] \arrow[d,"\theta"] \\
				K_0(C^*(H\Lb)) \arrow[r,swap,"(\rho_H)_*"] \arrow[d,"\bar{\theta}_H"] \arrow[rrr,"(\iota_H)_*",dashed,bend left = 8] & K_0(C^*(H\Lb \times_{d_1}\Z)\rtimes_{\lt_1}\Z) \arrow[r,"\iota_*"] \arrow[d] & K_0(C^*(\Lb \times_{d_1}\Z)\rtimes_{\lt_1}\Z) \arrow[d] & K_0(C^*(\Lb)) \arrow[l,"\rho_*"] \arrow[d,"\bar{\theta}"] \\
				\varinjlim(\ker(1-A_{H,2}^t),A_{H,1}^t|) \arrow[r,"\cong"] \arrow[d,"1-A_{H,1}^{t,\infty}"] \arrow[rrr,swap,"\iota",dashed,bend right = 8] & K_1(C^*(H\Lb\times_{d_1}\Z)) \arrow[r,"(\iota_{H\times\Z})_*"] \arrow[d,"1-(\lt_1)_*"] & K_1(C^*(\Lb\times_{d_1}\Z)) \arrow[d,"1-(\lt_1)_*"] & \varinjlim(\ker(1-A_{2}^t),A_{1}^t|) \arrow[l,swap,"\cong"] \arrow[d,"1-A_1^{t,\infty}"] \\
				\varinjlim(\ker(1-A_{H,2}^t),A_{H,1}^t|) \arrow[r,"\cong"] \arrow[d,"\xi_H"] \arrow[rrr,swap,"\iota",dashed,bend right = 8] & K_1(C^*(H\Lb\times_{d_1}\Z)) \arrow[r,"(\iota_{H\times\Z})_*"] \arrow[d,"(\iota_1)_*"] & K_1(C^*(\Lb\times_{d_1}\Z)) \arrow[d,"(\iota_1)_*"] & \varinjlim(\ker(1-A_{2}^t),A_{1}^t|) \arrow[l,swap,"\cong"] \arrow[d,"\xi"]\\
				K_1(C^*(H\Lb)) \arrow[r,"(\rho_H)_*"] \arrow[uuuuu,"\bar{\xi}_H",dotted,bend left = 40] \arrow[rrr,swap,"(\iota_H)_*",dashed,bend right = 8] & K_1(C^*(H\Lb \times_{d_1}\Z)\rtimes_{\lt_1}\Z) \arrow[r,"\iota_*"] \arrow[uuuuu,dotted,bend left = 60] & K_1(C^*(\Lb \times_{d_1}\Z)\rtimes_{\lt_1}\Z) \arrow[uuuuu,dotted,bend right = 60] & K_1(C^*(\Lb)) \arrow[l,swap,"\rho_*"] \arrow[uuuuu,"\bar{\xi}",dotted,bend right = 40]
			\end{tikzcd}
		};
	\end{tikzpicture}
\]
commutes and has exact columns.

There are unique homomorphisms $\theta_H, \theta, \bar{\theta}_H, \bar{\theta}, \xi^H, \xi, \bar{\xi}_H, \bar{\xi}$ that make the squares in which they appear commute because the adjoining horizontal maps are isomorphisms. Since the middle two columns are Pimsner--Voiculescu sequences, and hence exact, it follows that all four columns are exact as well. The first, second, fourth, and fifth rows commute by Corollary~\ref{cor: main k theory theorem}(1), and the third and sixth rows commute as they are induced from the commuting diagram of Lemma~\ref{lem: takai duality}. The middle vertical band of squares is~\eqref{eq:nested cyclic}, which commutes by naturality of the Pimsner--Voiculescu sequence. The remaining squares commute either by definition or by Corollary~\ref{cor: main k theory theorem}(2).
	
Deleting the middle two columns and retaining the outer two columns and dashed connecting maps yields the commuting diagram	
	\[
	\begin{tikzcd}[nodes={inner sep=2.5pt}]
		\varinjlim(\coker(1-A_{H,2}^t),\tilde{A}_{H,1}^t) \arrow[r, "1-\tilde{A}_{H,1}^{t,\infty}"] \arrow[d, "\iota", dashed]                  & \varinjlim(\coker(1-A_{H,2}^t),\tilde{A}_{H,1}^t) \arrow[d, "\iota", dashed] \arrow[r, "\theta_H"] & K_0(C^*(H\Lambda)) \arrow[ddd,swap, "\bar{\theta}_H", bend left=75] \arrow[d, "(\iota_H)_*", dashed]                        \\
		\varinjlim(\coker(1-A_{2}^t),\tilde{A}_{1}^t) \arrow[r, "1-\tilde{A}_1^{t,\infty}"]                                                & \varinjlim(\coker(1-A_{2}^t),\tilde{A}_{1}^t) \arrow[r, "\theta"]                                           & K_0(C^*(\Lambda)) \arrow[d,"\bar{\theta}",bend left =60]              \\
		K_1(C^*(\Lambda)) \arrow[u,pos=0.6,"\bar{\xi}", bend left=60]                                               & \varinjlim(\ker(1-A_{2}^t),A_{1}^t|)\arrow[l, "\xi"]                                           & \varinjlim(\ker(1-A_{2}^t),A_{1}^t|)\arrow[l, "1-A_1^{t,\infty}"]                                           \\
		K_1(C^*(H\Lb)) \arrow[uuu,swap,pos=0.47,"\bar{\xi}_H", bend left=75] \arrow[u, "(\iota_H)_*", dashed] & \varinjlim(\ker(1-A_{H,2}^t),A_{H,1}^t|) \arrow[u, "\iota", dashed] \arrow[l, "\xi_H"]             & \varinjlim(\ker(1-A_{H,2}^t),A_{H,1}^t|) \arrow[l, "1-A_{H,1}^{t,\infty}"] \arrow[u, "\iota", dashed]
\end{tikzcd}\]
of exact sequences.
	
Applying the first isomorphism theorem to the commuting diagram of five-term exact sequences centred on $K_0(C^*(H\Lb))$ and on $K_0(C^*(\Lb))$ yields the commuting diagram
	\[
	\begin{tikzcd}
		0 \arrow[r] & \coker(1-\tilde{A}_{H,1}^{t,\infty}) \arrow[r, "j_H"] \arrow[d,"\tilde{\iota}"] & K_0(C^*(H\Lb)) \arrow[r] \arrow[d,"(\iota_H)_*"] & \ker(1-A_{H,1}^{t,\infty}) \arrow[r] \arrow[d,"\iota"] & 0 \\
		0 \arrow[r] & \coker(1-\tilde{A}_{1}^{t,\infty}) \arrow[r, "j"] & K_0(C^*(\Lb)) \arrow[r] & \ker(1-A_{1}^{t,\infty}) \arrow[r] & 0
	\end{tikzcd}
	\]
with exact rows. (In the bottom row of this diagram, the cokernel is that of $1-\tilde{A}_{1}^{t,\infty}$ regarded as an endomorphism of $\varinjlim(\coker(1-A_{2}^t),\tilde{A}_{1}^t)$, and the kernel is that of $1-A_{1}^{t,\infty}$ regarded as an endomorphism of $\varinjlim(\ker(1-A_{2}^t),A_{1}^t|)$; and analogously in the top row.)
	
Likewise, applying the first isomorphism theorem to the commuting diagram of five term exact sequences centred on $K_1(C^*(H\Lb))$ and on $K_1(C^*(\Lb))$ gives the commuting diagram
	\[
	\begin{tikzcd}
		0 \arrow[r] & \coker(1-A_{H,1}^{t,\infty}) \arrow[r] \arrow[d,"\tilde{\iota}"] & K_1(C^*(H\Lb)) \arrow[r] \arrow[d,"(\iota_H)_*"] & \ker(1-\tilde{A}_{H,1}^{t,\infty}) \arrow[r] \arrow[d,"\iota"] & 0 \\
		0 \arrow[r] & \coker(1-A_{1}^{t,\infty}) \arrow[r] & K_1(C^*(\Lb)) \arrow[r] & \ker(1-\tilde{A}_{1}^{t,\infty}) \arrow[r] & 0
	\end{tikzcd}
	\]
with exact rows. (Again, here the kernel and cokernel in the bottom row are with respect to the endomorphisms $1-A_{1}^{t,\infty} \in \operatorname{End}\big(\varinjlim(\ker(1-A_{2}^t),A_{H,1}^t|)\big)$ and $1-\tilde{A}_{1}^{t,\infty} \in \operatorname{End}\big(\varinjlim(\coker(1-A_{2}^t),\tilde{A}_{1}^t)\big)$; and analogously in the top row.)
	
Lemma~\ref{lem:paskraeburn} now yields the commuting diagrams in the statement of Theorem~\ref{thm:inc}.
	
It remains to verify that $j$ and $j_H$ satisfy the desired formulas; the arguments are identical so we just argue for $j$. By definition, $j$ is induced by $\theta_H$ via the first isomorphism theorem, so
\[
j\big(\delta_v+\image(1-A_2^t) + \image(1-\tilde{A}_1^t)\big)
    = \theta([\delta_v + \image(1-A_2^t)]_1),
\]
so we must show that $\theta([\delta_v + \image(1-A_2^t)]_1) = [p_v]$.

Let $\tilde{\phi}_0 : \varinjlim(\coker(1 - A^t_2), \tilde{A}^t_1) \to K_0(C^*(\Lambda \times_{d_1} \Z)$ be as in Corollary~\ref{cor: main k theory theorem}(1). Then $\theta$ is, by definition, the homomorphism for which the diagram
\[
\begin{tikzcd}
	K_0(C^*(\Lb\times_{d_1}\Z)) \arrow[d, "(\iota_1)_*"] & {\varinjlim(\coker(1-A_2^t),\tilde{A}_1^t)} \arrow[d, "\theta"] \arrow[l, "\tilde{\phi}_0"] \\
	K_0(C^*(\Lb\times_{d_1}\Z)\rtimes_{\lt_1}\Z)           & K_0(C^*(\Lb)) \arrow[l, "\phi_*"]
\end{tikzcd}\]
commutes. Hence
\[
    \phi_* \circ \theta([\delta_v + \image(1-A_2^t)]_1) = (\iota_1)_* \circ \tilde{\phi}_0([\delta_v+\image(1-A_2^t)]_1) = (\iota_1)_*[p_{(v,1)}] =[\iota_1(p_{(v,1)})].
\]
The projection $\iota_1(p_{(v,1)})$ is equivalent to $\iota_1(p_{(v,0)})$ in $C^*(\Lb\times_{d_1}\Z)\rtimes_{\lt_1}\Z$ as in the proof of Proposition~\ref{prop: Bn step 1}. Hence $[\iota_1(p_{(v,1)})] =[\iota_1(p_{(v,0)})] = \phi_*[p_v]$. Hence $\theta[\delta_v+\image(1-A_2^t)]_1 = [p_v]$ as required, because $\phi_*$ is injective.
\end{proof}

We can re-interpret the diagram describing $K_0$-groups in Theorem~\ref{thm:inc} in terms of block matrices as in \cite{Evans2008}. Regard the block matrices
\[
(1-A_1^t,1-A_2^t) \in M_{\Lb^0 \times (\Lb^0 \sqcup \Lb^0)}(\Z)
\quad\text{ and }\quad
\begin{pmatrix}1-A_1^t \\ 1- A_2^t\end{pmatrix} \in M_{(\Lb^0 \sqcup \Lb^0) \times \Lb^0}(\Z)
\]
as homomorphisms from $\Z\Lb^0 \oplus \Z\Lb^0$ to $\Z\Lb^0$ and from $\Z\Lb^0$ to $\Z\Lb^0 \oplus \Z\Lb^0$ respectively, and similarly for $H\Lambda$.

Since $H$ is invariant under the $A^t_i$, the inclusion of $\Z H$ in $\Z\Lambda$ induces homomorphisms
\[
\tilde{\iota} : \coker(1-A_{H,1}^t,1-A_{H,2}^t) \to \coker(1-A_1^t,1-A_2^t)
    \quad\text{and} \quad
\iota : \ker\begin{pmatrix}1-A_{H,1}^t \\ 1- A_{H,2}^t\end{pmatrix} \to \ker\begin{pmatrix}1-A_1^t \\ 1- A_2^t\end{pmatrix}.
\]

\begin{theorem}\label{thm:inc nice}
Let $\Lb$ be a row-finite $2$-graph with no sources, and let $H$ be a hereditary subset of $\Lambda^0$. Let $\tilde{\iota}$ and $\iota$ be as above, and let $(\iota_H)_* : K_0(C^*(H\Lambda)) \to K_0(C^*(\Lambda))$ be the homomorphism induced by the inclusion of $C^*$-algebras. There is a homomorphism $j : \coker(1-A_1^t,1-A_2^t) \to K_0(C^*(\Lambda))$ such that $j(\delta_v + \image(1-A_1^t,1-A_2^t)) = [p_v]$, and a corresponding homomorphism $j_H$ for $H\Lambda$, and there are homomorphisms $\tau_H, \tau$ such that the diagram	
\[
	\begin{tikzcd}
		0 \arrow[r] & \coker(1-A_{H,1}^t,1-A_{H,2}^t) \arrow[r,"j_H"] \arrow[d,"\tilde{\iota}"] & K_0(C^*(H\Lambda)) \arrow[r,"\tau_H"] \arrow[d,"(\iota_H)_*"] & \ker\begin{pmatrix}1-A_{H,1}^t \\ 1- A_{H,2}^t\end{pmatrix} \arrow[r] \arrow[d,"\iota"] & 0 \\
		0 \arrow[r] & \coker(1-A_1^t,1-A_2^t) \arrow[r,"j"] & K_0(C^*(\Lambda)) \arrow[r,"\tau"] & \ker\begin{pmatrix}1-A_1^t \\ 1- A_2^t\end{pmatrix} \arrow[r] & 0
	\end{tikzcd}
	\]
commutes and has exact rows.
\end{theorem}

\begin{proof}
We have $\ker([1-A_1^t]|_{\ker(1-A_2^t)}) = \ker(1-A_1^t) \cap \ker(1-A_2^t) = \ker\big(\begin{smallmatrix}1-A_1^t \\ 1- A_2^t\end{smallmatrix}\big)$. So the result will follow once we show that there is an isomorphism $\coker(1 - \tilde{A}^t_1) \cong \coker(1-A_1^t,1-A_2^t)$ that carries $y + \image(1 - A^t_2)$ to $y + \image(1-A_1^t,1-A_2^t)$, and similarly for $H$. For $y \in \Z\Lb^0$, we have
\begin{align*}
y + \image(1-A_2^t) &{}\in \image(1-\tilde{A}_1^t)\\
    &\iff y-(1-A_1^t)\tilde{y} \in \image(1-A_2^t)\; \text{ for some $\tilde{y}\in \Z\Lb^0$}\\
    &\iff y = (1-A_1^t)\tilde{y} + (1-A_2^t)\bar{y}\; \text{ for some $\tilde{y},\bar{y} \in \Z\Lb^0$}\\
    &\iff y \in \image(1-A_1^t,1-A_2^t).
\end{align*}
The argument for $H$ is identical.
\end{proof}

\begin{remark}
	In \cite[Proposition 3.16]{Evans2008}, Evans shows that $$K_0(C^*(\Lb))\cong \coker(1-A_1^t,1-A_2^t) \oplus \ker\begin{pmatrix}A_2^t-1 \\ 1- A_1^t\end{pmatrix}$$ and $$K_1(C^*(\Lb)) \cong \ker(1-A_1^t,1-A_2^t)/\image\begin{pmatrix}1-A_1^t \\ 1- A_2^t\end{pmatrix}.$$ The advantage of Evans' method using spectral sequences is the neat description of the $K$-groups, and in particular the explicit description of $K_1(C^*(\Lambda))$. The drawback is that we were unable to deduce from Evans' approach whether the sequence
	\[
	\begin{tikzcd}
		0 \arrow[r] & {\coker(1-A_1^t,1-A_2^t)} \arrow[r] & K_0(C^*(\Lb)) \arrow[r] & \ker\begin{pmatrix}A_2^t-1 \\ 1- A_1^t\end{pmatrix} \arrow[r] & 0
	\end{tikzcd}\] was natural with respect to the inclusion of $C^*(H\Lb)$ into $C^*(\Lb)$. The advantage of our approach using Pimsner--Voiculescu sequences alone was that we were able to obtain Theorem~\ref{thm:inc nice}. The drawback of our method is that we do not obtain the splitting for $K_0$, nor a nice description of $K_1$.
	
	Note that in Evans' result, there is a switch of signs and in the block matrices ordering compared to our description. Evans' conventions were chosen to match the standard form of the maps in the Koszul resolution that he uses to compute the homology groups appearing in the Kasparov spectral sequence (see \cite[pp. 11--12]{Evans2008} and \cite[Corollary~4.5.5]{Wiebel}). However, the kernels and images of his maps are identical to those of ours, so that the difference is just one of presentation.
\end{remark}

\section{Stably finite extensions}\label{sec:stablyfiniteextensions}

Let $\Lambda$ be a row-finite $2$-graph with no sources. The \emph{matrix condition}
    	\begin{gather}
        	\{(1-A_1^t)f + (1-A_2^t)g : f,g \in \Z\Lb^0\}\cap \mathbb{N}\Lambda^0 = \{0\}, \tag{M}\label{matrixcond}
        \end{gather}
was introduced in \cite{CaHS}. Theorem~1.1 of \cite{CaHS} implies that if $\Lb$ is cofinal in the sense that it admits no saturated hereditary subsets other than $\emptyset$ and $\Lb^0$, then $C^*(\Lb)$ is stably finite if and only if $\Lb$ satisfies~\eqref{matrixcond}. When $H$ is a hereditary subset, the extra condition that $H$ is also saturated guarantees that $\Lambda \setminus \Lambda H$ has no sources. Since we defined the matrix condition~\eqref{matrixcond} specifically for row-finite $2$-graphs with no sources, we require that $H$ is saturated in the following lemma.

	\begin{lemma}\label{lem:scl1}
		Let $\Lambda$ be a row-finite $2$-graph with no sources and let $H$ be a saturated hereditary subset. There is a homomorphism $\tilde{\iota} : \coker(1-A_{H,1}^t,1-A_{H,2}^t) \to \coker(1-A_1^t,1-A_2^t)$ such that $\tilde{\iota}(f + \image(1-A_{H,1}^t,1-A_{H,2}^t)) = f + \image(1-A_1^t,1-A_2^t)$ for all $f \in \Z H$. If $\Lb$ satisfies~\eqref{matrixcond}, then $H\Lb$ satisfies \eqref{matrixcond} and
		\begin{equation}\label{eq:lem1}
			\ker(\tilde{\iota})\cap  [\N H + \image(1-A_{H,1}^t,1-A_{H,2}^t)] = \{0\}.
		\end{equation}
		Conversely, if $\Lambda \setminus \Lambda H$ and $H\Lambda$ both satisfy~\eqref{matrixcond} and~\eqref{eq:lem1} holds, then $\Lambda$ satisfies~\eqref{matrixcond}.
	\end{lemma}	
\begin{proof}
Let $\iota$ be the inclusion $\Z H \hookrightarrow \Z \Lb^0$. Lemma~\ref{lem:fundamental} implies that $\iota(\image(1 - A^t_{H,1}, 1 - A^t_{H_2})) \subseteq \image(1 - A^t_1, 1 - A^t_2)$, which shows that there exists a homomorphism $\tilde\iota$ as described. It also implies that if $\Lambda$ satisfies~\eqref{matrixcond}, then
\begin{align*}
\{(1-A_{H,1}^t)f + (1-A_{H,2}^t)g &{}: f,g \in \Z H\}\cap \mathbb{N} H\\
    &\subseteq \{(1-A_1^t)f + (1-A_2^t)g : f,g \in \Z\Lb^0\}\cap \mathbb{N} \Lambda^0 = \{0\},
\end{align*}
so that $H\Lambda$ also satisfies~\eqref{matrixcond}. Moreover, if $m \in \N H$ and $\tilde{\iota}(m + \image(1-A_{H,1}^t,1-A_{H,2}^t)) = 0$, then $\iota(m) \in \image(1-A_1^t,1-A_2^t) \cap \N\Lb^0 = \{0\}$, and since $\iota$ is injective, we deduce that $m = 0$. So~\eqref{eq:lem1} is satisfied.

For the other direction, suppose that $H\Lb$ and $\Lb \setminus \Lb H$ satisfy (\ref{matrixcond}), and that
\[
    \ker(\tilde{\iota})\cap  [\N H + \image(1-A_{H,1}^t,1-A_{H,2}^t)] = \{0\}.
\]
Fix $f \in \image(1-A_{1}^t,1-A_2^t)\cap \N\Lb^0$ and write $f = \begin{pmatrix} f_1 \\ f_2	\end{pmatrix}$ where $f_1 \in \N (\Lb^0\setminus H)$ and $f_2 \in \N H$, and fix $g, h \in \Z\Lb^0$ such that $f = (1-A_1^t)g + (1-A_2^t)h$. To avoid notational clutter, in this proof we write $B_i^t \coloneqq  A_{\Lb^0\setminus H,i}^t$ for $i=1,2$. Then writing
$g = \begin{pmatrix} g_1 \\ g_2 \end{pmatrix}$ and $h = \begin{pmatrix}	h_1 \\ h_2
\end{pmatrix}$ for $g_1, h_1 \in \Z (\Lb^0 \setminus H)$ and $g_2, h_2 \in \Z H$, we have
\[
	\begin{pmatrix}
		f_1 \\ f_2
	\end{pmatrix} = \begin{pmatrix}
		1 - B_1^t & 0 \\ * & 1- A_{H,1}^t
	\end{pmatrix} \begin{pmatrix} g_1 \\ g_2
	\end{pmatrix} + \begin{pmatrix}
		1 - B_2^t & 0 \\ * & 1- A_{H,2}^t
	\end{pmatrix} \begin{pmatrix}
		h_1 \\ h_2
	\end{pmatrix}.
\]
Thus $f_1 = (1-B_1^t)g_1 + (1-B_2^t)h_1 \in \image(1-B_1^t,1-B_2^t) \cap \N(\Lb^0\setminus H)$. Note that the $B_i$ are precisely the $1$-graph adjacency matrices for $(\Lb \setminus \Lb H)_i$ and that $(\Lb \setminus \Lb H)^0 = \Lb^0 \setminus H$. Thus the matrix condition (\ref{matrixcond}) for $\Lb\setminus \Lb H$ implies that $f_1 = 0$. In particular,
\[
    0 = f + \image(1-A_1^t,1-A_2^t) = \tilde{\iota}(f_2 + \image(1-A_{H,1}^t,1-A_{H,2}^t)),
\]
so that $f_2 \in \ker(\tilde{\iota})\cap[\N H + \image(1-A_{H,1}^t,1-A_{H,2}^t)] = \{0\}$. So $f=0$.
\end{proof}
	
\begin{remark}
For the purposes of this paper, we only need the first direction of Lemma~\ref{lem:scl1}; we include the converse for completeness and in case it is useful in future work.
\end{remark}
	
Let $j :\coker(1-A_1^t,1-A_2^t) \to K_0(C^*(\Lb))$ be as in Theorem~\ref{thm:inc nice}. Then
\begin{align*}
    j(\N \Lb^0 + \image(1-A_1^t,1-A_2^t)) = \N \{[p_v] : v \in \Lb^0\} \subset \Z \{[p_v] : v \in \Lb^0\} & \cap K_0(C^*(\Lb))_+ \\ & = \image(j) \cap K_0(C^*(\Lb))_+.
\end{align*}
It is not known when the other inclusion holds. It does hold for many examples, including all $2$-graphs that either have no red cycles or no blue cycles by the combination of \cite[Lemma~1.5]{Spielberg1988} and \cite[Theorem~5.1]{CaHS}. So we take it as an assumption in our main theorem later.
	
We say that a row-finite $2$-graph $\Lambda$ with no sources satisfies the \textit{positive-elements condition} if
\begin{gather}\label{apecond}\tag{P}
	j[\N \Lb^0 + \image(1-A_1^t,1-A_2^t)] = \image(j)\cap K_0(C^*(\Lb))_+.
\end{gather}
	
There are various equivalent ways of describing condition \eqref{apecond}.
	
\begin{proposition}\label{prop:equivalent conditions}
Let $\Lb$ be a row-finite $2$-graph with no sources. The following conditions are equivalent.
\begin{enumerate}
	\item $\Lambda$ satisfies condition \eqref{apecond}.
	\item $\Lambda$ satisfies the condition 	
    \begin{gather}\label{nicecond}\tag{N}
		\Z \{[p_v] : v \in \Lb^0\} \cap K_0(C^*(\Lb))_+ = \N \{[p_v] : v \in \Lb^0\}.
	\end{gather}
	\item Writing $\iota_1 : C^*(\Lb\times_{d_1}\Z) \to C^*(\Lb\times_{d_1}\Z)\rtimes_{\lt_1}\Z$
    for the natural inclusion of the skew-graph algebra into its crossed product, we have
    \begin{gather}\label{KCP}\tag{K}
		(\iota_1)_*\big(K_0(C^*(\Lb\times_{d_1}\Z))_+\big)
            = \image((\iota_1)_*)\cap K_0(C^*(\Lb\times_{d_1}\Z)\rtimes_{\lt_1}\Z)_+.
	\end{gather}
\end{enumerate}
\end{proposition}

	
\begin{remark}
Of the conditions in Proposition \ref{prop:equivalent conditions}, \eqref{nicecond} is the most natural and conceptually easy to grasp. In the proof of our main result, we will use the condition in the form \eqref{apecond} for the graph $H\Lb$ because the map $j_H$ appears in a commuting diagram. We include the third condition becase it is a familiar sticking point for crossed-product $C^*$-algebras, so it provides context for what may look, on the face of it, like an innocuous hypothesis.
\end{remark}
	
To prove the above proposition, we first need the following lemma.
	
\begin{lemma}\label{lem: positive cone of skew-graph K0}
Let $\Lambda$ be a row-finite $2$-graph with no sources. Then
\[
    K_0(C^*(\Lb \times_{d_1}\Z))_+ = \N\{[p_{(v,n)}] : v\in\Lb^0, n\in \N\}.
\]
\end{lemma}
\begin{proof}
Taking the $B^\Lb_n$ and $\iota_n:B^\Lb_n \to B^\Lb_{n+1}$ from Section~\ref{sec: k theory main} for which $K_0(C^*(\Lb\times_{d_1}\Z)) = \varinjlim(B^\Lb_n,\iota_n)$, we know that $$\varinjlim(K_0(B^\Lb_n),(\iota_n)_*) \cong K_0(C^*(\Lb\times_{d_1} \Z)).$$
		
Let $\Lb_2$ be the sub-$1$-graph $\{\lb \in \Lb : d_1(\lb) = 0\}$ and let $Q_n \in \mathcal{M}(B^\Lb_n)$ be the projection $\sum_{v\in \Lb^0}p_{(v,n)}$. In the spirit of Section \ref{sec: k theory main}, it is not hard to see that $C^*(\Lb_2) \cong Q_n B^\Lb_n Q_n$ and that $B^\Lb_n = B^\Lb_n Q_n B^\Lb_n$. Thus, $C^*(\Lb_2)$ is isomorphic to a full corner in $B_n$. Hence, for each $n\in \Z$, there is an isomorphism $K_0(C^*(\Lb_2)) \cong K_0(B^\Lb_n)$ that carries a class $[p_v]$ to the class $[p_{(v,n)}]$.
		
Since $\Lambda_2$ is a $1$-graph, \cite[Theorem 7.1]{AMP2007} implies that $K_0(C^*(\Lb_2))_+ = \N\{[p_v] : v \in \Lb^0\}$ (this is because, in their notation, for a directed graph $E$ the positive cone $K_0(C^*(E))_+$ is the monoid $V(C^*(E))$ and $\N\{[p_v] : v \in E^0\}$ is $M_E$). In particular, $K_0(B^\Lb_n)_+ = \N\{[p_{(v,n)}] : v \in \Lb^0\}$ and $(\iota_n)_*[K_0(B^\Lb_n)_+] \subset K_0(B^\Lb_{n+1})_+$ since the inclusion map $\iota_n$ corresponds to the nonnegative-integer matrix $A_1^t$ (Proposition~\ref{prop: Bn step 1}). Thus, writing $\iota_{n,\infty}: B^\Lb_n \to C^*(\Lb\times_{d_1}\Z)$ for the inclusion maps,
\[
    K_0(C^*(\Lb\times_{d_1}\Z))_+ = \bigcup_{n \in \N} \iota_{n,\infty} (K_0(B^\Lb_n)_+) = \N\{[p_{(v,n)}] : v \in \Lb^0, n \in \N\}.\qedhere
\]
\end{proof}
	
\begin{proof}[Proof of Proposition~\ref{prop:equivalent conditions}]
Conditions \eqref{apecond}~and~\eqref{nicecond} are clearly equivalent. We prove that conditions \eqref{apecond}~and~\eqref{KCP} are also equivalent.
		
Consider the commuting diagram
		\[ \begin{tikzcd}
			{\varinjlim(\coker(1-A_2^t),\tilde{A}_1^t)} \arrow[d, "\tilde{\phi}_0"] \arrow[r, "{1-\tilde{A}_1^{t,\infty}}"] & {\varinjlim(\coker(1-A_2^t),\tilde{A}_1^t)} \arrow[d, "\tilde{\phi}_0"] \arrow[r, "\theta"] & K_0(C^*(\Lb)) \arrow[d, "\rho_*"]          \\
			K_0(C^*(\Lb\times_{d_1}\Z)) \arrow[r, "1-(\lt_1)_*"]                                                                  & K_0(C^*(\Lb\times_{d_1}\Z)) \arrow[r, "(\iota_1)_*"]                                                     & K_0(C^*(\Lb\times_{d_1}\Z)\rtimes_{\lt_1}\Z)
		\end{tikzcd} \]
comprised of the top two squares of the right-most column of the 24-term diagram in the proof of Theorem~\ref{thm:inc}. Let $\pi$ denote the quotient map from $\varinjlim(\coker(1-A_2^t),\tilde{A}_1^t)$ onto $\coker(1-\tilde{A}_1^{t,\infty})$ and let $\tilde{\theta}$ be the induced map such that $\tilde{\theta}\circ \pi = \theta$. Let $R: \coker(1-\tilde{A}_1^{t,\infty})\to \coker(1-\tilde{A}_1^t), \; [\delta_v + \image(1-A_2^t)]_n + \image(1-\tilde{A}_1^{t,\infty}) \mapsto \delta_v + \image(1-A_2^t) + \image(1-\tilde{A}_1^t)$ be the induced isomorphism from Lemma~\ref{lem:paskraeburn} and let $j$ be the map in Theorem~\ref{thm:inc}. Then the diagram
		\[\begin{tikzcd}
			{\varinjlim(\coker(1-A_2^t),\tilde{A}_1^t)} \arrow[r, "\pi"] \arrow[rd, "\theta"] & {\coker(1-\tilde{A}_1^{t,\infty})} \arrow[r, "R"] \arrow[d, "\tilde{\theta}"] & \coker(1-\tilde{A}_1^t) \arrow[ld, "j"] \\
			& K_0(C^*(\Lambda))                                                             &
		\end{tikzcd}\]
commutes.
Now
\[
    \N\Lb^0 + \image(1-A_2^t) + \image(1-\tilde{A}_1^t) = R\big(\N\{[\delta_v+\image(1-A_2^t)]_n + \image(1-\tilde{A}_1^t) : v\in \Lb^0, n\in \N\}\big)
\]
and $\image(j) = \image(\tilde{\theta})$. Hence~\eqref{apecond} holds if and only if
\begin{equation}\label{apecond2}
\tilde{\theta}\big(\N\{[\delta_v+\image(1-A_2^t)]_n
        + \image(1-\tilde{A}_1^t) : v\in \Lb^0, n\in \N\}\big)
    = \image(\tilde{\theta})\cap K_0(C^*(\Lb))_+.
\end{equation}
We have
\[
\N\{[\delta_v+\image(1-A_2^t)]_n + \image(1-\tilde{A}_1^t) : v\in \Lb^0, n\in \N\}
    = \pi\left(\N\{[\delta_v+\image(1-A_2^t)]_n : v\in \Lb^0, n\in \N\}\right),
\]
and since $\pi$ is surjective, $\image(\tilde{\theta}) = \image(\theta)$. Hence~\eqref{apecond2} holds if and only if
\begin{equation}\label{apecond3}
\theta[\N\{[\delta_v+\image(1-A_2^t)]_n : v\in \Lb^0, n\in \N\}]
    = \image(\theta)\cap K_0(C^*(\Lb))_+.
\end{equation}
		
Since $\phi_* \circ\theta = (\iota_1)_*\circ\tilde{\phi}_0$ and $\tilde{\phi}_0$ is an isomorphism,
\[
\phi_*(\image(\theta))
    = \image(\phi_* \circ \theta)
    = \image((\iota_1)_*\circ\tilde{\phi}_0)
    = \image((\iota_1)_*).
\]
By Takai duality, $\phi:C^*(\Lb) \to C^*(\Lb\times_{d_1}\Z)\rtimes_{\lt_1}\Z$ induces an order isomorphism in $K$-theory, so
\[
    \phi_*(K_0(C^*(\Lb))_+) = K_0(C^*(\Lb\times_{d_1}\Z)\rtimes_{\lt_1}\Z)_+.
\]
Applying $\phi_*$ to both sides of~\eqref{apecond3}, we see that~\eqref{apecond3} holds if and only if
\begin{equation}\label{apecond4}
(\iota_1)_*\circ \tilde{\phi}_0 \left([\N\{[\delta_v+\image(1-A_2^t)]_n : v\in \Lb^0, n\in \N\}\right)
    = \image((\iota_1)_*)\cap K_0(C^*(\Lb\times_{d_1}\Z)\rtimes_{\lt_1}\Z)_+.
\end{equation}

Finally,
\[
\tilde{\phi}_0 \left(\N\{[\delta_v+\image(1-A_2^t)]_n : v\in \Lb^0, n\in \N\}\right)
    = \N\{[p_{(v,n)}] : v\in\Lb^0, n\in \N\} = K_0(C^*(\Lb\times_{d_1}\Z))_+.
\]
Hence~\eqref{apecond4} is equivalent to~\eqref{KCP}.
\end{proof}
	
We can now state and prove our main result.
	
\begin{theorem}\label{thm:simplified stably finite restated}
Let $\Lb$ be a row-finite $2$-graph with no sources and $H \subset \Lb^0$ be a saturated hereditary subset. Assume that
\begin{gather}
  	\{(1-A_1^t)f + (1-A_2^t)g : f,g \in \Z\Lb^0\}\cap \mathbb{N}\Lambda^0 = \{0\},\text{ and}
    \tag{\ref{matrixcond}}\\
	\Z \{[p_v] : v \in H\} \cap K_0(C^*(H\Lb))_+ = \N \{[p_v] : v \in H\}.
    \tag{\ref{nicecond}}
\end{gather}
If $C^*(H\Lb)$ and $C^*(\Lb\setminus \Lb H)$ are stably finite then $C^*(\Lb)$ is stably finite.
\end{theorem}
\begin{proof}
Let $I_H \coloneqq  \overline{\spann}\{s_\mu s_\nu^* : s(\mu)= s(\nu) \in H\}$. Consider the following short exact sequence:
\[
\begin{tikzcd}
	0 \arrow[r] & I_H \arrow[r,"i"] & C^*(\Lb) \arrow[r] & C^*(\Lb)/I_H \arrow[r] & 0.
\end{tikzcd}
\]
Then $C^*(H\Lb)$ is stably isomorphic to $I_H$ and $C^*(\Lb\setminus \Lb H)$ is isomorphic to $C^*(\Lb)/I_H$ \cite[Theorem 5.2]{RaeburnSimsYeend2003}, so $I_H$ and $C^*(\Lb)/I_H$ are stably finite because $C^*(H\Lb)$ and $C^*(\Lb\setminus \Lb H)$ are.
		
By Spielberg's result Theorem~\ref{thm:spielberg}, $C^*(\Lb)$ is stably finite if and only if \begin{gather}\label{jack}\tag{S}
	\ker{i_*}\cap K_0(I_H)_+ = \{0\}.
\end{gather}
Let $\iota_H : C^*(H\Lb) \to C^*(\Lb)$ be the canonical inclusion. Since $C^*(H\Lb)$ embeds into $I_H$ which embeds into $C^*(\Lb)$, and since $K_0(C^*(H\Lb)) \cong K_0(I_H)$, Equation~\eqref{jack} holds if and only if
\begin{equation}\label{jack2}
	\ker{(\iota_H)_*}\cap K_0(C^*(H\Lb))_+ = \{0\}.
\end{equation}

We are going to use Theorem \ref{thm:inc nice}, so we adopt its notation. The map
\[
\iota : \ker\begin{pmatrix}1-A_{H,1}^t \\ 1- A_{H,2}^t\end{pmatrix} \subset \Z H
    \to \ker\begin{pmatrix}1-A_{1}^t \\ 1- A_{2}^t\end{pmatrix} \subset \Z \Lb^0
\]
is injective. So if $x \notin \image(j_H) = \ker(\tau_H)$, then $x \notin \ker(\iota \circ \tau_H) = \ker(\tau \circ (\iota_H)_*) \supset \ker ((\iota_H)_*)$. Hence~\eqref{jack2} is equivalent to
\begin{equation}\label{jack3}
	\ker((\iota_H)_*)\cap \left(\image(j_H)\cap K_0(C^*(H\Lb))_+\right) = \{0\}.
\end{equation}

Since $H\Lb$ satisfies~\eqref{nicecond}, and hence~\eqref{apecond} by Proposition~\ref{prop:equivalent conditions},
\[
    \image(j_H)\cap K_0(C^*(H\Lb))_+ = j_H\left(\N H + \image(1-A_{H,1}^t,1-A_{H,2}^t)\right).
\]
Thus~\eqref{jack3} is equivalent to
\[
	\ker((\iota_H)_*) \cap j_H\left(\N H + \image(1-A_{H,1}^t,1-A_{H,2}^t)\right) = \{0\},
\]
and hence to
\begin{equation}\label{jack5}
	\ker((\iota_H)_* \circ j_H)\cap \left(\N H + \image(1-A_{H,1}^t,1-A_{H,2}^t)\right) = \{0\}.
\end{equation}
		
By exactness of the bottom row of the diagram in Theorem~\ref{thm:inc nice}, the map $j : \coker(1 - A^t_1, 1- A^t_2) \to K_0(C^*(\Lambda))$ is injective. So $\ker((\iota_H)_* \circ j_H) = \ker(j\circ \tilde{\iota}) = \ker(\tilde{\iota})$. The first part of Lemma~\ref{lem:scl1} gives
\[
    \ker(\tilde{\iota})\cap \left(\N H + \image(1-A_{H,1}^t,1-A_{H,2}^t)\right) = \{0\},
\]
giving~\ref{jack5} and hence~\eqref{jack}, whence $C^*(\Lb)$ is stably finite.
\end{proof}

To conclude, we show that we can bootstrap from Theorem~\ref{thm:simplified stably finite restated} to obtain a result for $2$-graphs with finite chains of saturated hereditary subsets. Observe that this includes all $2$-graphs with finitely many vertices. We need the following definition. A \textit{cycle in a coordinate graph $\Lb_i$} of a $2$-graph $\Lb$ is a path $\lb \in \Lb \setminus \Lb^0$ such that $r(\lb) = s(\lb)$ and $d(\lb) \in \N e_i$.

\begin{corollary}\label{cor: maximal chain}
Let $\Lambda$ be a row-finite $2$-graph with no sources. Suppose that $\Lambda$ has a finite chain
\begin{equation}\label{maximal chain}
	H_0\coloneqq \emptyset\subsetneq H_1\subsetneq \dots \subsetneq  H_{n}\subsetneq H_{n+1}\coloneqq \Lambda^0
\end{equation}
of saturated hereditary subsets that is maximal in the sense that for each $m \ge 1$ there is no saturated hereditary set $H$ such that $H_{m-1} \subsetneq H \subsetneq H_m$. Suppose that $\Lambda\setminus \Lambda H_m$ satisfies \eqref{matrixcond} for each $m=0,\dots,n$ and that $H_m\Lambda$ satisfies \eqref{nicecond} for $m=1,\dots,n$. Then $C^*(\Lambda)$ is stably finite.  If there are no cycles in the coordinate graph $\Lb_1$, or if there are no cycles in the coordinate graph $\Lb_2$, then $C^*(\Lambda)$ is AF-embeddable.
\end{corollary}
\begin{proof}
Since $H_1$ has no proper saturated hereditary subsets, $H_1\Lambda$ is cofinal by \cite[Lemma~5.2]{PSS2020}. As $\Lb = \Lb \setminus \Lb H_0$ satisfies \eqref{matrixcond} so does $H_1\Lb$ by Lemma~\ref{lem:scl1}. Hence $C^*(H_1\Lambda)$ is stably finite by \cite[Theorem~1.1(1)]{CaHS}.
		
Fix $m\in \{1,\dots, n\}$ and assume that $C^*(H_m\Lambda)$ is  stably finite. We need to show that $C^*(H_{m+1}\Lambda)$ is stably finite by applying Theorem~\ref{thm:simplified stably finite restated} to the $2$-graph $H_{m+1}\Lb$ and the saturated hereditary subset $H_m \subset H_{m+1} = (H_{m+1}\Lb)^0$. Again $H_{m+1}\Lb$ satisfies \eqref{matrixcond} by Lemma~\ref{lem:scl1} as $\Lb$ does. Furthermore, as $H_m \Lb$ satisfies \eqref{nicecond} and $C^*(H_m \Lb)$ is stably finite by assumption, it remains to show that $C^*(H_{m+1}\Lb \setminus ((H_{m+1}\Lb)H_m))$ is stably finite.
		
Since the chain \eqref{maximal chain} is maximal, the set of vertices $H_{m+1}\setminus H_m$ of $(H_{m+1}\Lambda)\setminus ((H_{m+1} \Lambda) H_{m})$ has no proper saturated hereditary subsets, and hence $H_{m+1}\Lambda\setminus ((H_{m+1} \Lambda) H_{m})$  is cofinal by \cite[Lemma~5.2]{PSS2020}. As $\Lambda\setminus \Lambda H_m$ satisfies \eqref{matrixcond}, so does $H_{m+1}\Lambda\setminus ((H_{m+1} \Lambda) H_{m}) = (H_{m+1}\setminus H_m)(\Lb \setminus \Lb H_m)$ by Lemma~\ref{lem:scl1}. Thus the $C^*$-algebra of $H_{m+1}\Lambda \setminus ((H_{m+1} \Lambda) H_{m})$ is stably finite by \cite[Theorem~1.1(1)]{CaHS}.
		
So Theorem~\ref{thm:simplified stably finite restated} implies that $C^*(H_{m+1}\Lambda)$ is stably finite. By induction, $C^*(\Lb) = C^*(H_{n+1}\Lb)$ is stably finite.
		
The final statement about AF-embeddability follows from \cite[Theorem 5.1]{CaHS}.
\end{proof}

\section{Examples}\label{sec:examples}

We illustrate our results by discussing two 2-graphs. In both examples, the bulk of the work is in verifying that the $2$-graph satisfies the matrix condition (M). We call the first coordinate graph $\Lambda^{e_1}$ the \emph{blue subgraph} and the second coordinate graph  $\Lambda^{e_2}$ the \emph{red subgraph}. Cycles in the blue and red subgraphs are called blue and red cycles respectively.

\subsection{}\label{subsec:spine example} Throughout this subsection, let $\Lambda$ be the $2$-graph with skeleton shown in Figure~\ref{figure-skeleton}.

\begin{figure}[ht]\label{figure-skeleton}
	\[
	\begin{tikzpicture}[>=stealth, decoration={markings, mark=at position 0.5 with {\arrow{>}}}]
		\node[circle, inner sep=0pt] (v1) at (0,0) {$v_1$};
		\node[circle, inner sep=0pt] (v2) at (2,0) {$v_2$};
		\node[circle, inner sep=0pt] (v3) at (4,0) {$v_3$};
		\node[circle, inner sep=0pt] (v4) at (6,0) {$v_4$};
		\node[circle, inner sep=0pt] at (7.5,0) {$\dots$};
		\node[circle, inner sep=0pt] (w1) at (0,2) {$w_1$};
		\node[circle, inner sep=0pt] (w2) at (2,2) {$w_2$};
		\node[circle, inner sep=0pt] (w3) at (4,2) {$w_3$};
		\node[circle, inner sep=0pt] (w4) at (6,2) {$w_4$};
		\node[circle, inner sep=0pt] at (7.5,2) {$\dots$};
		\node[circle, inner sep=0pt] (x1) at (4.5,-2) {$x_1$};
		\node[circle, inner sep=0pt] (x2) at (5.5,-2) {$x_2$};
		\node[circle, inner sep=0pt] (x3) at (6.5,-2) {$x_3$};
		\node[circle, inner sep=0pt] (x4) at (7.5,-2) {$x_4$};
		\node[circle, inner sep=0pt] at (8.5,-2) {$\dots$};
		\node[circle, inner sep=0pt] (y1) at (1.5,-2) {$y_1$};
		\node[circle, inner sep=0pt] (y2) at (0.5,-2) {$y_2$};
		\node[circle, inner sep=0pt] (y3) at (-0.5,-2) {$y_3$};
		\node[circle, inner sep=0pt] (y4) at (-1.5,-2) {$y_4$};
		\node[circle, inner sep=0pt] at (-2.5,-2) {$\dots$};
		\draw[blue, postaction=decorate, out=150, in=30] (v2) to (v1);
		\draw[blue, postaction=decorate, out=150, in=30] (v3) to (v2);
		\draw[blue, postaction=decorate, out=150, in=30] (v4) to (v3);
		\draw[red, dashed, postaction=decorate, out=210, in=330] (v2) to (v1);
		\draw[red, dashed, postaction=decorate, out=210, in=330] (v3) to (v2);
		\draw[red, dashed, postaction=decorate, out=210, in=330] (v4) to (v3);
		\draw[blue, postaction=decorate, out=240, in=120] (w1) to (v1);
		\draw[blue, postaction=decorate, out=240, in=120] (w2) to (v2);
		\draw[blue, postaction=decorate, out=240, in=120] (w3) to (v3);
		\draw[blue, postaction=decorate, out=240, in=120] (w4) to (v4);
		\draw[red, dashed, postaction=decorate, out=300, in=60] (w1) to (v1);
		\draw[red, dashed, postaction=decorate, out=300, in=60] (w2) to (v2);
		\draw[red, dashed, postaction=decorate, out=300, in=60] (w3) to (v3);
		\draw[red, dashed, postaction=decorate, out=300, in=60] (w4) to (v4);
		\draw[blue, postaction=decorate] (w1).. controls +(-0.75, 0.75) and +(-0.75, -0.75) .. (w1);
		\draw[blue, postaction=decorate] (w2).. controls +(-0.75, 0.75) and +(-0.75, -0.75) .. (w2);
		\draw[blue, postaction=decorate] (w3).. controls +(-0.75, 0.75) and +(-0.75, -0.75) .. (w3);
		\draw[blue, postaction=decorate] (w4).. controls +(-0.75, 0.75) and +(-0.75, -0.75) .. (w4);
		\draw[red, dashed, postaction=decorate] (w1).. controls +(0.75, 0.75) and +(0.75, -0.75) .. (w1);
		\draw[red, dashed, postaction=decorate] (w2).. controls +(0.75, 0.75) and +(0.75, -0.75) .. (w2);
		\draw[red, dashed, postaction=decorate] (w3).. controls +(0.75, 0.75) and +(0.75, -0.75) .. (w3);
		\draw[red, dashed, postaction=decorate] (w4).. controls +(0.75, 0.75) and +(0.75, -0.75) .. (w4);
		\draw[blue, postaction=decorate] (x1) to (v1);
		\draw[blue, postaction=decorate] (x1) to (v2);
		\draw[blue, postaction=decorate] (x1) to (v3);
		\draw[blue, postaction=decorate] (x1) to (v4);
		\draw[blue, postaction=decorate] (x2) to (x1);
		\draw[blue, postaction=decorate] (x3) to (x2);
		\draw[blue, postaction=decorate] (x4) to (x3);
		\draw[red, dashed, postaction=decorate] (x1).. controls +(0.75, -0.75) and +(-0.75, -0.75) .. (x1);
		\draw[red, dashed, postaction=decorate] (x2).. controls +(0.75, -0.75) and +(-0.75, -0.75) .. (x2);
		\draw[red, dashed, postaction=decorate] (x3).. controls +(0.75, -0.75) and +(-0.75, -0.75) .. (x3);
		\draw[red, dashed, postaction=decorate] (x4).. controls +(0.75, -0.75) and +(-0.75, -0.75) .. (x4);
		\draw[red, dashed, postaction=decorate] (y1) to (v1);
		\draw[red, dashed, postaction=decorate] (y1) to (v2);
		\draw[red, dashed, postaction=decorate] (y1) to (v3);
		\draw[red, dashed, postaction=decorate] (y1) to (v4);
		\draw[red, dashed, postaction=decorate] (y2) to (y1);
		\draw[red, dashed, postaction=decorate] (y3) to (y2);
		\draw[red, dashed, postaction=decorate] (y4) to (y3);
		\draw[blue, postaction=decorate] (y1).. controls +(0.75, -0.75) and +(-0.75, -0.75) .. (y1);
		\draw[blue, postaction=decorate] (y2).. controls +(0.75, -0.75) and +(-0.75, -0.75) .. (y2);
		\draw[blue, postaction=decorate] (y3).. controls +(0.75, -0.75) and +(-0.75, -0.75) .. (y3);
		\draw[blue, postaction=decorate] (y4).. controls +(0.75, -0.75) and +(-0.75, -0.75) .. (y4);
		\node[anchor=north] at (3, -3) {Figure~\ref{figure-skeleton}};
	\end{tikzpicture}
	\]
\end{figure}

Notice that $\Lambda$ is not cofinal and hence  \cite[Theorem~1.1]{CaHS} does not apply; also $\Lambda$ has blue and red cycles  and hence  \cite[Theorem~5.1]{CaHS} does not apply.  We will show that  $C^*(\Lambda)$ is stably finite using Theorem~\ref{thm:simplified stably finite restated}.

To do this, we will apply the theorem to  the saturated hereditary subset \[H\coloneqq \{x_n, y_n\colon n\geq 1\},\] which is the disjoint union of the saturated hereditary subsets $H_1\coloneqq \{x_n \colon n\geq 1\}$ and $H_2\coloneqq\{y_n \colon n\geq 1\}$.  Notice that the constant function $1$ is a faithful graph trace on $H_i\Lambda$ for $i=1,2$.   Since each $H_i\Lambda$ is cofinal,   each $C^*(H_i\Lambda)$ is AF-embeddable by \cite[Theorem~1.1]{CaHS}.  So the direct sum  $C^*(H\Lambda)$  is AF-embeddable as well.

Let $E$ be the blue subgraph of the skeleton of the quotient graph $\Lambda\setminus H\Lambda$. Since no cycle of the directed graph $E$ has an entrance, $C^*(E)$ is AF-embeddable by \cite{Schafhauser}. Notice that $\Lambda\setminus H\Lambda$  is the pullback of  $E$  by $f \colon \N^2 \to \N$ given by $f(m) = m_1 + m_2$. By \cite[Corollary~3.5(iii)]{KP2000},  $C^*(\Lambda\setminus H\Lambda)$ is
isomorphic to  $C^*(E) \otimes C(\T)$. So if we choose an embedding $i \colon C^*(E) \to A$ of
$C^*(E)$ in an AF algebra $A$ and let $\pi \colon \cantor \to \T$ be a continuous surjection of
the Cantor set $\cantor$ onto $\T$, then $i \otimes \pi_*$ is an embedding of $C^*(E) \otimes
C(\T)$ in the AF-algebra $A \otimes C(\cantor)$.  Thus $C^*(\Lambda\setminus H\Lambda)$ is AF-embeddable.

It remains to verify that $\Lambda$ satisfies the matrix condition (M) and that $H\Lambda$ satisfies condition (N).

\begin{lemma}\label{lemma matrix condition}
The graph  $\Lambda$ satisfies the matrix condition (M).
 \end{lemma}
 \begin{proof}
Let $A_1, A_2$ be the coordinate matrices for  $\Lambda$. We need to show that \[
\Big(\sum_{i=1}^2\image(1-A^t_i) \Big) \cap  \N{\Lambda^0} = \{0\}.
 \]
Let $g,h\in \Z \Lambda^0$ and suppose that
\[
f\coloneqq(1-A_1^t)g+(1-A_2^t)h\in\N \Lambda^0;
\]
we will  show that $f=0$.

 We write $v, w, x, y$ for the ordered sets of vertices $\{v_n\colon n\geq 1\}$, $\{w_n\colon n\geq 1\}$, $\{x_n\colon n\geq 1\}$ and $\{y_n\colon n\geq 1\}$, respectively. Then  with respect to the ordering $v\cup w\cup x\cup y$ we have
\[
 1-A_1^t
=   \left[
    \begin{array}{c|c|c|c}
      * & 0 &0 &0 \\
      \hline
      -1  & 0  &0 &0 \\
       \hline
      * & 0  &*&0 \\
 \hline
     0  & 0 &0 &0 \\
    \end{array}
    \right]
\quad\text{and}\quad
 1-A_2^t
=   \left[
    \begin{array}{c|c|c|c}
      * & 0 &0 &0 \\
      \hline
      -1  & 0  &0 &0 \\
       \hline
     0& 0  &0&0 \\
 \hline
     *  & 0 &0 &* \\
    \end{array}
    \right]
\]
where the blocks labelled $*$ are given by:
\begin{align*}
( 1-A_1^t)(v_m, v_n)=( 1-A_2^t)(v_m, v_n)&=( 1-A_1^t)(x_m, x_n)=(1-A_2^t)(y_m, y_n)\\
& =\begin{cases}1&\text{if $n=m$}\\
 -1&\text{if $m+1=n$}\\
 0&\text{else;}
 \end{cases}\\
 (1-A_1^t)(x_m, v_n)=(1-A_2^t)(y_m, v_n)&=\begin{cases}-1&\text{if $m=1$}\\
 0&\text{else.}
  \end{cases}
\end{align*}
We write $f=f^v+f^w+f^x+f^y$ where $f^v=f|_v$ and $A_1^v$ for the $(v,v)$ block of $A_1$ (and similarly for $w,x,y$ and $A_2$, and for $g,h \in \Z\Lb^0$ below).
We will show that $f^v$, $f^w$, $f^x$ and $f^y$ are all zero, and then $f=0$ as well.

We have
\[
f^v=(1-(A_1^v)^t)g^v+(1-(A_2^v)^t)h^v\in\N (v\Lambda v)^0.
\]
In $v\Lambda v$, we have $v_n\Lambda^{e_1}v_{n+1}=v_n\Lambda^{e_2}v_{n+1}$ for all $n$ and hence $C^*(v\Lambda v)$  is AF-embeddable by  \cite[Proposition~7.2]{CaHS}. Since $v\Lambda v$ is cofinal, $v\Lambda v$ satisfies the matrix condition by \cite[Theorem~1.1]{CaHS}. Thus $f^v=0$.

Since $1-(A_1^v)^t=1-(A_2^v)^t$ we get $g^v+h^v\in\ker(1-(A_1^v)^t)=\{0\}$. Thus $g^v=-h^v$.
This gives
\begin{align*}
f(w_n)&=((1-A_1^t)g)(w_n)+((1-A_2^t)h)(w_n)\\
&=\sum_{u\in\Lambda^0}(1-A_1^t)(w_n, u)g(u)+\sum_{u\in\Lambda^0}(1-A_2^t)(w_n, u)h(u)\\
&=-g(v_n)+(-h(v_n))=0
\end{align*}
for all $n\ge 1$.
Thus $f^w=0$.

Since the entries of the  rows of $1-A_2^t$ corresponding to  any $x_n$ are all $0$ we have $((1-A_2^t)h)^x=0$.
Now another calculation gives
\[f(x_n)=((1-A_1^t)g)(x_n)=\begin{cases}
-\sum_{i=1}^\infty g(v_i)+g(x_1)&\text{if $n=1$}\\
g(x_n)-g(x_{n-1})&\text{if $n>1$}.
\end{cases}
\]
Thus $f^x= (g(x_1)-\sum_{i=1}^\infty g(v_i))\delta_{x_1}+\sum_{n=2}^\infty (g(x_n)-g(x_{n-1}))\delta_{x_n}$.
Since $f^x\geq 0$ we get the chain
\[
\cdots\geq g(x_n)\geq\cdots\geq g(x_2)\geq g(x_1)\geq \sum_{i=1}^\infty g(v_i).
\]
Similarly,
\[
0\leq f(y_n)=0+((1-A_2^t)h)(y_n)=\begin{cases}
-\sum_{i=1}^\infty h(v_i)+h(y_1)&\text{if $n=1$}\\
h(y_n)-h(y_{n-1})&\text{if $n>1$}
\end{cases}
\]
and we get the chain
\[
\cdots \geq h(y_n)\geq\cdots\geq h(y_2)\geq h(y_1)\geq \sum_{i=1}^\infty h(v_i).
\]
We claim that $\sum_{i=1}^\infty g(v_i)=0$. If $\sum_{i=1}^\infty g(v_i)>0$, then
\[
\cdots\geq  g(x_n)\geq\cdots\geq g(x_2)\geq g(x_1)\geq \sum_{i=1}^\infty g(v_i)>0,
\]
which implies that $g$ has infinite support, a contradiction. If $\sum_{i=1}^\infty g(v_i)<0$, then since $g^v=-h^v$
we have
\[
\geq \cdots h(y_n)\geq\cdots\geq h(y_2)\geq h(y_1)\geq \sum_{i=1}^\infty h(v_i)=-\sum_{i=1}^\infty g(v_i)>0,
\]
which implies that $h$ has infinite support, also a contradiction. Thus  $\sum_{i=1}^\infty g(v_i)=0$ as claimed, and it follows that $\sum_{i=1}^\infty h(v_i)=0$ as well.  Now we have
\begin{gather*}
0\leq f^x = g(x_1)\delta_{x_1}+\big(g(x_2)-g(x_1)\big)\delta_{x_2}+\big(g(x_3)-g(x_2)\big)\delta_{x_3}+\dots=(1-(A_1^x)^t)g^x;\\
0\leq f^y=h(y_1)\delta_{y_1}+\big(h(y_2)-h(y_1)\big)\delta_{y_2}+\big(h(y_3)-h(y_2)\big)\delta_{y_3}+\dots=(1-(A_2^y)^t)h^y.
\end{gather*}
The blue graph $x_1\Lambda^{e_1}$ and the red graph $y_1\Lambda^{e_2}$  have vertex matrices $A_1^x$ and $A_2^y$, respectively, and
their $C^*$-algebras  are the compact operators. Thus $x_1\Lambda^{e_1}$  and   $y_1\Lambda^{e_2}$  satisfy the matrix condition for $1$-graphs by \cite[Lemma~4.2]{CaHS}. Thus $f^x=0=f^y$. We have   shown that $f=0$ and hence that $\Lambda$ satisfies the matrix condition.
\end{proof}

\begin{lemma}  The graph $H\Lambda$ satisfies the positive-elements condition (P).
\end{lemma}

\begin{proof} For condition (P) we have to show that
\[j\big(\N H+ \image(1-A_{H,1}^t,1-A_{H,2}^t)\big) \supset  \image(j)\cap K_0(C^*(H\Lb))_+.\]
Since $H$ is the disjoint union of the saturated hereditary subset $H_1=\{x_n \colon n\geq 1\}$ and $H_2=\{y_n \colon n\geq 1\}$, it suffices to show that this inclusion holds for $H_1\Lambda$.

Let $B_1$ and $B_2$ be the vertex matrices for $H_1\Lambda$. For $n\geq 1$ we have  $(1-B_1^t)\delta_{x_n}=\delta_{x_n}-\delta_{x_{n+1}}$, and $1-B_2^t=0$. In particular, \[\delta_{x_1}+ \image(1-B_1^t,1-B_2^t)=\delta_{x_{n}}+ \image(1-B_1^t,1-B_2^t).\]
Let $x\in  \image(j)\cap K_0(C^*(H_1\Lb))_+$.  Then there exist a finite subset $F$ of $\{1, 2, \dots\}$ and $a_n\in\Z$ for $n\in F$ such that  $x=j\big(\sum_{n\in F} a_n\delta_{x_n} + \image(1-B_1^t,1-B_2^t)\big)$.  Then
\[
x=j \Big(\Big(\sum_{n\in F} a_n \Big)\delta_{x_1} + \image(1-B_1^t,1-B_2^t) \Big)=\Big(\sum_{n\in F} a_n\big) [p_{x_1}]\in K_0(C^*(H_1\Lb))_+.
\]
Thus $\sum_{n\in F} a_n\in\N$ and $x\in j\big(\N H+ \image(1-B_1^t,1-B_2^t)\big)$ as needed.
\end{proof}

Condition (P) is the same as condition (N) by Proposition~\ref{prop:equivalent conditions}.  So we have now verified all the hypotheses of Theorem~\ref{thm:simplified stably finite restated} and we conclude that $C^*(\Lambda)$ is stably finite.

\subsection{}\label{subsec:handlebar example} Throughout this subsection, $\Lambda$ is  the $2$-graph of  \cite[Example~6.3]{ES} with skeleton shown in Figure~\ref{figure-ES}.
\begin{figure}[ht]\label{figure-ES}
\[
\begin{tikzpicture}[>=stealth, decoration={markings, mark=at position 0.5 with {\arrow{>}}}]
    \filldraw[black!10!white] (-5.5,-1.5)--(-1,-1.5)--(-1,7.5)--(-5.5,7.5)--cycle;
    \foreach \x in {-4,-2,0,2,4} {
        \foreach \y in {0,2,4,6} {
            \node[circle, inner sep=1.5pt, fill=black] (v\x\y) at (\x,\y) {};
            \ifnum \y > 0 {
                \ifnum \x = 0 {
                    \draw[postaction=decorate, blue] (v0\y) .. controls +(0.25,-.95) .. +(0,-1.9)
                        \ifnum \y=2 node[pos=0.5, anchor=west, inner sep=1pt]  {$\alpha$}\fi;
                    \draw[postaction=decorate, red] (v0\y) .. controls +(-0.25,-.95) .. +(0,-1.9)
                        \ifnum \y=2 node[pos=0.5, anchor=east, inner sep=1pt] {$\beta$}\fi;
                } \else {
                    \ifnum \x > 0 {
                        \draw[postaction=decorate, red] (v\x\y)--+(0,-1.9);
                    } \else {
                        \draw[postaction=decorate, blue] (v\x\y)--+(0,-1.9);
                    }\fi
                } \fi
            } \fi
            \ifnum \x < 0 {
                \draw[postaction=decorate, red] (v\x\y)--+(1.9,0);
            } \fi
            \ifnum \x > 0 {
                \draw[postaction=decorate, blue] (v\x\y)--+(-1.9,0)
                \ifnum \x=2{\ifnum \y=0 node[pos=0.5, anchor=north, inner sep=2pt] {$\gamma$}\fi}\fi;
            } \fi
        }
        \node at (\x,-1) {$\vdots$};
        \node at (\x,7) {$\vdots$};
    }
    \foreach \y in {0,2,4,6} {
        \node at (-5,\y) {$\dots$};
        \node at (5,\y) {$\dots$};
    }
    \node[anchor=north, inner sep=1pt] at (v00.south) {$v_{0,0}$};
    \node[anchor=north] at (0, -2) {Figure~\ref{figure-ES}};
\end{tikzpicture}
\]
\end{figure}

We pick a vertex of $\Lambda$  that emits blue and red edges with the same range, and name it $v_{0,0}$. Then we put a $\Z^2$ grid on the set of vertices. We call the collection of vertices $\{v_{0,n}\colon n\in\Z\}$ the  spine of the graph.  The hereditary and saturated subset $H$ of vertices to the left of the spine is then $H=\{v_{m,n}\colon  m<0, n\in\Z\}$, and the graph $\Lambda H$ is shown shaded in Figure~\ref{figure-ES}. The quotient graph $\Lambda\setminus \Lambda H$ has the generalised cycle (in the sense of \cite{ES}) $(\alpha,\beta)$ with entrance the blue edge $\gamma$ with range equal to the range of $\alpha$ and $\beta$. Hence the vertex projection at the range of $\gamma$ is infinite. Thus the $C^*$-algebra of $\Lambda\setminus \Lambda H$ is not stably finite.

Nevertheless, we can show that $C^*(\Lambda)$ is AF-embeddable: since $\Lambda$  contains no blue or red cycles, it suffices by \cite[Theorem~5.1]{CaHS} to verify that the matrix condition holds. (We could also apply Theorem~\ref{thm:simplified stably finite restated} to the hereditary and saturated subset $K=\{v_{m,n}\colon m\neq 0\}$ --- it is easy to see that $K\Lb$ satisfies~\eqref{nicecond} --- to conclude that $C^*(\Lambda)$ is stably finite, but using \cite[Theorem~5.1]{CaHS} is obviously more efficient and gives a stronger result.)

\begin{lemma} The graph $\Lambda$ satisfies the matrix condition (M).
\end{lemma}

\begin{proof}
Let $A_1, A_2$ be the vertex matrices for $\Lambda$. Let $0\neq f\in \big(\image(1-A^t_1) +\image(1-A^t_2) \big) \cap \N{\Lambda^0}$ and look for a contradiction. Fix $g,h\in \Z\Lb^0$ such that
\[
f=(1-A^t_1) g+(1-A^t_2)h.
\]
 Our first claim is that $f(v_{0,n})=0$ for all $n\in\Z$. We  have
 \begin{align*}
0\leq f(v_{0,n})&=\big( (1-A^t_1) g+(1-A^t_2)h \big) (v_{0,n})\\
&=g(v_{0,n})-g(v_{0,n-1})+h(v_{0,n})-h(v_{0,n-1}).
\end{align*}
 Let $\Gamma\coloneqq\Lambda\{v_{0,n}\colon n\in\Z\}$ and $M_i=A_i|_{\Gamma^{0}}$ be the vertex matrices for $\Gamma$. Then the $2$-graph $\Gamma$ is the graph with vertices the spine of $\Lambda$ and each vertex $v_{0,n}$ of $\Gamma$ emits one blue and red edge with range $v_{0,n-1}$. The $C^*$-algebra of this graph $\Gamma$ is AF-embeddable by \cite[Proposition~7.2]{CaHS}, and hence $\Gamma$ satisfies the matrix condition by \cite[Theorem~1.1]{CaHS}. We have
\begin{align*}
 f|_{\Gamma^{0}}&=((1-A^t_1) g)|+((1-A^t_2)h)|\\
 &= (1-M^t_1) g|+(1-M^t_2)h|\in \big(\image(1-M^t_1) +\image(1-M^t_2) \big) \cap \N{\Lambda^0}
\end{align*}
 and hence $f|_{\Gamma^{0}}=0$. Thus $f(v_{0,n})=0$ for all $n\in\Z$, as claimed.

 Next we claim that $g(v_{0,n})=-h(v_{0,n})$ for all $n\in\Z$. If $g(v_{0,n})=h(v_{0,n})=0$ for all $n\in\Z$, then we are done.  So suppose there exists $n\in\Z$ such that one of $g(v_{0,n})$ or $h(v_{0,n})$ is nonzero.  Since $g,h$ have finite support, we can let
 \[
 n_0\coloneqq\max\{n\in\Z: \text{one of $g(v_{0,n})$ or $h(v_{0,n})$ is nonzero}\}.
 \]
 Then $g(v_{0,n})=h(v_{0,n})=0$ for $n>n_0$. We have
 \begin{align*}
 0=f(v_{0,n_0+1})&=g(v_{0,n_0+1})-g(v_{0,n_0})+h(v_{0, n_0+1})-h(v_{0,n_0})\\
 &=-g(v_{0,n_0})-h(v_{0,n_0}).
 \end{align*}
Thus $g(v_{0, n_0})=-h(v_{0, n_0})$. Now
\begin{align*}
0=f(v_{0, n_0})&=g(v_{0,n_0})-g(v_{0, n_0-1})+h(v_{0, n_0})-h(v_{0, n_0-1})\\
&=-g(v_{0, n_0-1})-h(v_{0, n_0-1}).
\end{align*}
Thus $g(v_{0, n_0-1})=-h(v_{0, n_0-1})$. By induction,  $g(v_{0,n})=-h(v_{0,n})$ for all $n\in\Z$, as claimed.

 Let $i=1$ and $n\in\Z$. Then
 \begin{align*}
 0&\leq f(v_{1,n})+f(v_{-1,n})\\
 &=g(v_{1,n})-g(v_{0,n})+h(v_{1,n})-h(v_{1,n-1})\\
 &\hskip4cm +g(v_{-1,n})-g(v_{-1,n-1})+h(v_{-1,n})-h(v_{0,n})\\
 &=g(v_{1,n})+h(v_{1,n})-h(v_{1,n-1}) +g(v_{-1,n})-g(v_{-1,n-1})+h(v_{-1,n})
  \end{align*}
using that $g(v_{0,n})=-h(v_{0,n})$ for all $n\in\Z$.
 Summing over $\Z$, and recognising the sum as telescoping, gives
 \[
 0\leq\sum_{n\in\Z}g(v_{1,n})+h(v_{-1,n}).
 \]
 For $i\geq 2$ we have
 \begin{align*}
 0&\leq f(v_{i,n})+f(v_{-i,n})\\
 &=g(v_{i,n})-g(v_{i-1,n})+h(v_{i,n})-h(v_{i,n-1})\\
 &\hskip4cm +g(v_{-i,n})-g(v_{-i,n-1})+h(v_{-i,n})-h(v_{-i+1,n}).
  \end{align*}
This time,  summing over $n$ gives
\[
 \sum_{n\in\Z}g(v_{i,n})+h(v_{-i,n})\geq  \sum_{n\in\Z}g(v_{i-1,n})+h(v_{-(i-1),n}).
\]
It now follows that
\[
 \sum_{n\in\Z}g(v_{i,n})+h(v_{-i,n})\geq  \sum_{n\in\Z}g(v_{i',n})+h(v_{-i',n})
\]
for $1\leq i'\leq i$. But since $f>0$,  there exists $i_0\geq 1$ such that $f(v_{i_0,j})+f(v_{-i_0,j})>0$ for some $j\in\Z$.  Then for all $i\geq i_0+1$ we have
\[
 \sum_{n\in\Z}g(v_{i,n})+h(v_{-i,n})>  \sum_{n\in\Z}g(v_{i_0,n})+h(-i_0,n).
\]
But this implies  that $g+h$ has infinite support, a contradiction. Thus $f=0$, as needed.
\end{proof}

It follows from the lemma and \cite[Theorem~5.1]{CaHS} that $C^*(\Lambda)$ is AF-embeddable.

\begin{remark} The graph $\Lambda$  illustrates that we cannot extend \cite[Theorem~1.1]{CaHS} to include the condition ``there exists a faithful graph trace'' if $\Lambda$ is not cofinal. To see that $\Lambda$ has no faithful graph trace we argue by contradiction. Suppose that $\Lambda$ has a faithful graph trace $\tau$. By rescaling, if necessary, we may assume that $\tau(v_{1,0})=1$.

Let $m\geq 1$. Since
\begin{gather*}
\tau(v_{m,n})=\sum_{\lambda\in v_{m,n}\Lambda^{e_1}}\tau(s(\lambda))=\tau(v_{m,n+1})\\
\tau(v_{m,n})=\sum_{\lambda\in v_{m,n}\Lambda^{e_2}}\tau(s(\lambda))=\tau(v_{m+1,n}),
\end{gather*}
it follows by induction that $\tau(v_{m,n})=1$ for $m\geq 1$ and $n\in\Z$. For each $i\in\N$ we have
\[
\tau(v_{0,0})=\sum_{\lambda\in v_{0,0}\Lambda^{(i,0)}}\tau(s(\lambda)).
\]
But now consider the sources of blue paths with range $v_{0,0}$ that have degree $(i,0)$: these are the vertices  $v_{0,i}, v_{1,i-1}, v_{2,i-2}, \dots,  v_{i,0}$. Thus $\tau(v_{0,0})\geq i-1$ for all $i\in\N$, a contradiction.
Thus $\Lambda$ does not have  a faithful graph trace even though $C^*(\Lambda)$ is AF-embeddable.
\end{remark}
	

\end{document}